\documentclass[preprint,noinfoline,authoryear]{imsart}
\usepackage{fullpage}
\setattribute{thebibliography}{size}{\small}
\renewcommand{\baselinestretch}{1.5}
\RequirePackage[OT1]{fontenc}
\RequirePackage{amsthm,amsmath,amssymb}
\RequirePackage[sort,round]{natbib}
\RequirePackage[breaklinks,colorlinks]{hyperref}
\RequirePackage{hypernat}

\hypersetup{colorlinks=false}

\usepackage{verbatim}
\usepackage{soul}
\usepackage{graphicx}
\usepackage{ctable}
\startlocaldefs
\DeclareMathOperator{\Var}{Var}
\DeclareMathOperator{\Cov}{Cov}

\DeclareMathOperator{\Cum}{Cum}
\newcommand{\LLL}{\mathcal{L}}

\newcommand{\floor}[1]{\lfloor #1 \rfloor}

\newcommand{\field}[1]{\mathbb{#1}}
\newcommand{\R}{\field{R}}

\newcommand{\N}{\field{N}}
\newcommand{\Z}{\field{Z}}

\newcommand{\E}{\field{E}}

\newcommand{\h}[1]{\boldsymbol{#1}}

\newcommand{\PP}{\mathcal{P}}
\newcommand{\FF}{\mathcal{F}}
\newcommand{\HH}{\mathcal{H}}
\newcommand{\CC}{\mathcal{C}}
\newcommand{\RR}{\mathcal{R}}

\def\II{I \negthinspace I}
\def\III{I \negthinspace I \negthinspace I}
\def\IV{I \negthinspace V}
\newcommand{\dd}{{\,\mathrm{d}}}
\theoremstyle{plain}
\newtheorem{theorem}{Theorem}
\newtheorem{corollary}[theorem]{Corollary}
\newtheorem{proposition}[theorem]{Proposition}
\newtheorem{lemma}[theorem]{Lemma}

\theoremstyle{definition}

\theoremstyle{definition}
\newtheorem{remark}{Remark}

\makeatletter

\newcommand{\Rmnum}[1]{\expandafter\@slowromancap\romannumeral #1@}
\makeatother

\newcommand{\comments}[1]{}

\newcounter{regular} 
\newcounter{stable} 

\endlocaldefs

\begin{document}

\begin{frontmatter}

  %\title{Asymetric Linear Regression with Nonstationay Errors}
  %\runtitle{Asymetric Regression with Nonstationay Errors}
%% \thankstext{T1}{Footnote to the title with the `thankstext'
%% command.}

\title{Asymptotic Inference of Autocovariances of Stationary Processes}
\runtitle{Sample Autocovariances of Stationary Processes}
\date{\today}

\begin{aug}
\author{\fnms{Han}
  \snm{Xiao}\ead[label=e1]{xiao@galton.uchicago.edu}}
\and
\author{\fnms{Wei Biao}
  \snm{Wu}\ead[label=e2]{wbwu@galton.uchicago.edu}}
%\snm{Wu}\thanksref{t1}\ead[label=e1]{wbwu@galton.uchicago.edu}}
%\thankstext{t1}{Research partially supported by NSF-DMS-0478704.}
\runauthor{H.~Xiao, W.B.~Wu}

\affiliation{University of Chicago}
\address{
Department of Statistics\\ 5734 S.~University Ave.\\ Chicago, IL 60637\\
\printead{e1} \\
\printead{e2}
}
\end{aug}

%\today

\begin{abstract}
The paper presents a systematic theory for asymptotic inference of
autocovariances of stationary processes. We consider nonparametric
tests for serial correlations based on the maximum (or ${\cal
L}^\infty$) and the quadratic (or ${\cal L}^2$) deviations. For
these two cases, with proper centering and rescaling, the
asymptotic distributions of the deviations are Gumbel and
Gaussian, respectively. To establish such an asymptotic theory, as
byproducts, we develop a normal comparison principle and propose a
sufficient condition for summability of joint cumulants of
stationary processes. We adopt a simulation-based block of blocks
bootstrapping procedure that improves the finite-sample
performance.
\end{abstract}



% history:
% \received{\smonth{1} \syear{0000}}

%\tableofcontents

\end{frontmatter}
\date{}

\renewcommand{\labelenumi}{\arabic{enumi}.}

%%%%%%%%%%%%%%%%%%%%%%%%%%%%%%%%%%%%%%%%%%%%%%%%%%%%%%%%
%%%%%%%%%%%%%%%%%%%%%%%%%%%%%%%%%%%%%%%%%%%%%%%%%%%%%%%%

\section{Introduction}
\label{sec:intro}

If $(X_i)_{i \in \Z}$ is a real-valued stationary process, then
from a second-order inference point of view it is characterized by
its mean $\mu = \E X_i$ and the autocovariance function $\gamma_k
= \E [(X_0-\mu) (X_k-\mu)]$, $k \in \Z$. Assume $\mu = 0$. Given
observations $X_1, \ldots, X_n$, the natural estimates of
$\gamma_k$ and the autocorrelation $r_k = \gamma_k / \gamma_0$ are
\begin{equation}
  \label{eq:autocov}
  \hat{\gamma}_k=(1/n)\sum_{i=|k|+1}^nX_{i-|k|}X_{i}
  \quad \hbox{and} \quad \hat r_k=\hat \gamma_k / \hat \gamma_0,
  \,\, 1-n \le k \le n-1,
\end{equation}
respectively. The estimator $\hat\gamma_k$ plays a crucial role in
almost every aspect of time series analysis. It is well-known that for
linear processes with {\it independent and identically distributed}
(iid) innovations, under suitable conditions, $\sqrt{n}(\hat{\gamma}_k
- \gamma_k) \Rightarrow \mathcal{N}(0,\tau_k^2)$, where $\Rightarrow$
stands for convergence in distribution, $\mathcal{N}(0,\tau_k^2)$
denotes the normal distribution with mean zero and variance
$\tau_k^2$. Here $\tau_k^2$ can be calculated by Bartlett's formula
(see Section 7.2 of \cite{brockwell:1991}). Other contributions on
linear processes include \cite{hannan:1972}, \cite{hosoya:1982},
\cite{anderson:1991} and \cite{phillips:1992} etc. \cite{romano:1996}
and \cite{wu:2009} considered the asymptotic normality of $\hat
\gamma_k$ for nonlinear processes. As a primary goal of the paper, we
shall study asymptotic properties of the quadratic (or ${\cal L}^2$)
and the maximum (or ${\cal L}^\infty$) deviations of $\hat
{\gamma}_k$.

\subsection{The ${\cal L}^2$ Theory}
Testing for serial correlation has been extensively studied in
both statistics and econometrics, and it is a standard diagnostic
procedure after a model is fitted to a time series. Classical
procedures include \cite{durbin:1950,durbin:1951},
\cite{box:1970}, \cite{robinson:1991} and their variants. The
Box-Pierce portmanteau test uses
% the sum of squares of first $K$ sample autocorrelations
% multiplied by the sample size
$Q_K = n \sum_{k=1}^K \hat r_k^2$ as the test statistic, and
rejects if it lies in the upper tail of $\chi^2_K$ distribution.
An arguable deficiency of this test and many of its modified
versions (for a review see for example \cite{escanciano:2009}) is
that the number of lags $K$ included in the test is held as a
constant in the asymptotic theory. As commented by
\cite{robinson:1991}:
\begin{quote}
''{\it ...unless the statistics take account of sample
autocorrelations at long lags there is always the possibility that
relevant information is being neglected...}''
\end{quote}
The problem is particularly relevant if practitioners have no
prior information about the alternatives. The attempt of
incorporating more lags emerged naturally in the spectral domain
analysis; see among others \cite{durlauf:1991}, \cite{hong:1996}
and \cite{deo:2000}. The normalized spectral density $f(\omega) =
(2\pi)^{-1} \sum_{k \in \Z} r_k \cos(k\omega)$ should equal to
$(2\pi)^{-1}$ when the serial correlation is not present. Let
$\hat f(\omega) = \sum_{k=1-n}^{n-1}h(k/s_n)\hat r_k\cos(k\omega)$
be the lag-window estimate of the normalized spectral density,
where $h(\cdot)$ is a kernel function and $s_n$ is the bandwidth
satisfying the natural condition $s_n \to \infty$ and $s_n /n \to
0$. The former aims to include correlations at large lags. A test
for the serial correlation can be obtained by comparing $\hat f$
and the constant function $f(\omega) \equiv (2\pi)^{-1}$ using a
suitable metric. In particular, using the quadratic metric and
rectangle kernel, the resulting test statistic is the Box-Pierce
statistic with unbounded lags. \cite{hong:1996} established the
following result:
\begin{equation}
  \label{eq:hong}
  \frac{1}{\sqrt{2s_n}}\left(n\sum_{k=1}^{s_n}
  (\hat r_k-r_k)^2 -s_n\right) \Rightarrow \mathcal{N}(0,1),
\end{equation}
under the condition that $X_i$ are iid, which implies that all
$r_k$ in the preceding equation are zero. \cite{lee:2001} and
\cite{duchesne:2010} studied similar tests in spectral domain, but
using a wavelet basis instead of trigonometric polynomials in
estimating the spectral density and henceforth working on wavelet
coefficients. \cite{fan:1996} considered a similar problem in a
different context and proposed {\it adapative Neyman} test and
thresholding tests, using $\max_{1\leq k \leq s_n} (Q_k-k) /
\sqrt{2k}$ and $n \sum_{k=1}^{s_n}\hat r_k^2I(|\hat r_k|>\delta)$
as test statistics respectively, where $\delta$ is a threshold
value. \cite{escanciano:2009} proposed to use $Q_{s_n}$ with $s_n$
being selected by AIC or BIC.

It has been an important and difficult question on whether the iid
assumption in \cite{hong:1996} can be relaxed. % One may naturally ask:
% (i) what if the underlying process $(X_i)$ is only uncorrelated but
% not iid in (\ref{eq:hong}) and (ii) what if the serial correlation is
% present in (\ref{eq:hong})? The special case of question (i) where
% $X_i$ are martingale differences has been studied in
% \cite{robinson:1991}, \cite{durlauf:1991} and \cite{deo:2000}.
% \cite{shao:2011} showed that (\ref{eq:hong}) is true when $(X_i)$ is a
% general white noise sequence, under the geometric moment contraction
% condition. To the best of our knowledge, there has been no results in
% the literature for (ii).
Similar problems have been studied by \cite{durlauf:1991},
\cite{deo:2000} and \cite{hong:2003} for the case that $X_i$ are
martingale differences. Recently \cite{shao:2011} showed that
(\ref{eq:hong}) is true when $(X_i)$ is a general white noise
sequence, under the geometric moment contraction (GMC)
condition. Since the GMC condition, which implies that the
autocovariances decay geometrically, is quite strong, the question
arises as to whether it can be replaced by a weaker one. Furthermore,
one may naturally ask: what if the serial correlation is present in
(\ref{eq:hong})?  To the best of our knowledge, there has been no
results in the literature for this problem. This paper shall address
these questions and substantially generalizes earlier results. We
shall prove that (\ref{eq:hong}) remains true even if all or some of
$r_k$ are not zero, but the variance of the limiting distribution,
being different, will depend on the values of $r_k$. Furthermore, we
derive the limiting distribution of $\sum_{k=1}^{s_n} \hat r_k^2$ when
the serial correlation is present. The latter result enables us to
calculate the asymptotic power of the Box-Pierce test with unbounded
lags.

\subsection{The ${\cal L}^\infty$ Theory}

Another natural omnibus choice is to use the maximum
autocorrelation as the test statistic. \cite{wu:2009} obtained a
stochastic upper bound for
\begin{equation}
  \label{eq:wu}
  \sqrt{n} \max_{1\leq k \leq s_n} |\hat\gamma_k-\gamma_k|,
\end{equation}
and argued that in certain situations the test based on
(\ref{eq:wu}) has a higher power over the Box-Pierce tests with
unbounded lags in detecting weak serial correlation. It turns out
that the uniform convergence of autocovariances is also closely
related to the estimation of orders of ARMA processes or linear
systems in general. The pioneer works in this direction were given
by E. J. Hannan and his collaborators, see for example
\cite{hannan:1974} and \cite{an:1982}. For a summary of these
works we recommend \cite[][Section~\S 5.3]{hannan:1988} and
references therein. In particular, \cite{an:1982} showed that if
$s_n = O[(\log n)^\alpha]$ for some $\alpha<\infty$, then with
probability one
\begin{align}
  \label{eq:hannan}
  \sqrt{n} \max_{1\leq k \leq s_n} |\hat\gamma_k-\gamma_k|
   = O\left(\log\log n\right).
\end{align}

% Assuming $s_n\to\infty$ and $s_n/n\to 0$, it will be
% useful to have uniform convergence rate of $\hat\gamma_k$ over $1\leq
% k\leq s_n$. For example, such a uniformity is necessary in estimating
% the autocovariance matrix, which is in turn crucial for almost all
% problems of statistical inferences. One of the main goals of this
% article is to study the asymptotic distribution of $\max_{1\leq k\leq
%   s_n}|\hat\gamma_k-\gamma_k|$.  Such a distributional theory can be
% used in constructing simultaneous confidence band of autocovariances,
% and in model selection. For example, in order to determine whether a
% moving average process of order $q$ is an appropriate model for a
% given set of observations, we need to test whether all the
% autocovariances $\gamma_k$ with $k\geq q$ are zero.

% To the best of our knowledge, there has been no result on the
% asymptotic distribution of $\max_{1\leq k\leq
%   s_n}|\hat\gamma_k-\gamma_k|$ in the literature. \cite{wu:2009}
% obtained a stochastic upper bound.

The question of deriving the asymptotic distribution of
(\ref{eq:wu}) is more challenging. Although \cite{wu:2009} was not
able to obtain the limiting distribution of (\ref{eq:wu}), his
work provided important insights into this problem. Assuming $k_n
\to \infty$, $k_n/n\to 0$ and $h\geq 0$, he showed that, for $T_k
= \sqrt n (\hat\gamma_{k}- \E \hat \gamma_{k})$,
\begin{equation}\label{eq:cltunbdd}
  \left(T_{k_n}\,,\, T_{k_n+h}\right)^\top
  \Rightarrow \mathcal{N}\left[0,
      \begin{pmatrix}
        \sigma_0 & \sigma_h \\
        \sigma_h & \sigma_0
      \end{pmatrix}\right], \quad
      \hbox{where } \sigma_h = \sum_{k\in\Z}
      \gamma_k \gamma_{k+h},
\end{equation}
and we use the superscript $\top$ to denote the transpose of a vector
or a matrix.  The asymptotic distribution in (\ref{eq:cltunbdd}) does
not depend on the speed of $k_n \to \infty$. It suggests that, at
large lags, the covariance structure of $(T_k)$ is asymptotically
equivalent to that of the Gaussian sequence
\begin{equation}
  \label{eq:gp}
  (G_k) := \left(\sum_{i\in\Z}\gamma_i\eta_{i-k}\right)
\end{equation}
where $\eta_i$'s are iid standard normal random variables. Define
the sequences $(a_n)$ and $(b_n)$ as
\begin{align}
  \label{eq:gumbel_constants}
  a_n=(2\log n)^{-1/2} \,\, \hbox{ and } \,\,
  b_n=(2\log n)^{1/2} - (8\log n)^{-1/2}(\log\log n + \log 4\pi).
\end{align}
According to \cite{berman:1964} (also see Remarks \ref{rk:berman}
and \ref{rk:absolute}), under the condition $\lim_{n \to \infty}
\E(G_0 G_n) \log n = 0$,
\begin{align*}
  \lim_{s\to\infty}P\left(\max_{1\leq i\leq s}|G_i| \leq
    \sqrt{\sigma_0}(a_{2s}\,x+b_{2s})\right) = \exp\{-\exp(-x)\}.
\end{align*}
Therefore, \cite{wu:2009} conjectured that under suitable
conditions, one has the Gumbel convergence
\begin{align}\label{eq:main1}
  \lim_{n\to\infty}
  P\left(\max_{1\leq k\leq s_n} |T_k|
    \leq \sqrt{\sigma_0}(a_{2s_n}\,x+b_{2s_n})\right)
    = \exp\{-\exp(-x)\}.
\end{align}
In a recent work, \cite{jirak:2011} proved this conjecture for linear
processes and for $s_n$ growing with at most logarithmic speed. We
shall prove (\ref{eq:main1}) in Section~\ref{sec:proof} for general
stationary processes; and our result allows $s_n$ to grow as $s_n =
O(n^\eta)$ for some $0 < \eta < 1$, and $\eta$ can be arbitrarily
close to $1$ under appropriate moment and dependence conditions. The
latter result substantially relaxes the severe restriction on the
growth speed in (\ref{eq:hannan}) and \cite{jirak:2011} and, moreover,
the obtained distributional convergence are more useful for
statistical inference. For example, other than testing for serial
correlation and estimating the order of a linear system,
(\ref{eq:main1}) can also be used to construct simultaneous confidence
intervals of autocovariances.

% Besides the requirement $s_n\to\infty$ and $s_n/n\to
% 0$, large value of $s_n$ is desired. It is not surprising that the
% maximum order that $s_n$ can take depends on the strength of
% dependence and the order of finite moments of $X_1$, as shown in
% Theorem~\ref{thm:single}. On the other hand, we would expect a higher
% order of $s_n$ if we observe $N$ iid replicates of the sequence
% $(X_1,\ldots,X_n)$. Since it is typical that the series is very long
% and the number of replicate is not very large, we shall pursue this
% problem under the condition $\limsup_{n\to\infty} N/n \leq c$,
% where $c$ is some positive constant.
%
% If the process is nonstationary, then it does not make sense to talk
% about autocovariances in general. However, it we observe iid
% replicates, we still can estimate $\gamma_{ij}:=\E X_iX_j$, and now the
% focus is on $\max_{1\leq i<h\leq n}|\hat\gamma_{ij}-\gamma_{ij}|$.
% The special case
% where $(X_1,\ldots,X_n)$ are iid was first studied by \cite{jiang:2004},
% followed by \cite{zhou:2007} and \cite{liu:2008}. Their method is
% Possion approximation \citep[see for example][]{arratia:1989}, which
% heavily depends on the fact that $X_i's$ are independent.

\subsection{Relations with the Random Matrix Theory}
In a companion paper, using the asymptotic theory of sample
autocovariances developed in this paper, \cite{wu:2010} studied
convergence properties of estimated covariance matrices which are
obtained by banding or thresholding. Their bounds are analogs under
the time series context to those of \cite{bickel:2008a,
  bickel:2008b}. There is an important difference between these two
settings: we assume that only one realization is available, while
\cite{bickel:2008a, bickel:2008b} require multiple iid copies of the
underlying random vector.

There has been some related works in the random matrix theory
literature that are similar to (\ref{eq:main1}). Suppose one has $n$
iid copies of a $p$-dimensional random vector, forming a $p \times n$
data matrix $\h{X}$. Let $\hat r_{ij}$, $1\le i,j \le p$, be the
sample correlations. \cite{jiang:2004} showed that the limiting
distribution of $\max_{1\leq i<j\leq p}|\hat r_{ij}|$, after suitable
normalization, is Gumbel provided that each column of $\h{X}$ consists
of $p$ iid entries and each entry has finite moment of some order
higher than 30, and $p/n$ converges to some constant. His work was
followed and improved by \cite{zhou:2007} and \cite{liu:2008}. In a
recent article, \cite{cai:2010} extended those results in two ways:
(i) the dimension $p$ could grow exponentially as the sample size $n$
provided exponential moment conditions; and (ii) they showed that the
test statistic $\max_{|i-j|>s_n} |\hat r_{ij}|$ also converges to the
Gumbel distribution if each column of $\h{X}$ is Gaussian and is
$s_n$-dependent. The latter generalization is important since it is
one of the very few results that allow dependent entries. Their method
is Poisson approximation \citep[see for example][]{arratia:1989},
which heavily depends on the fact that for each sample correlation to
be considered, the corresponding entries are
independent. \cite{schott:2005} proved that $\sum_{1\leq i<j\leq p}
\hat r_{ij}^2$ converges to normal distribution after suitable
normalization, under the conditions that each column of $\h{X}$
contains iid Gaussian entries and $p/n$ converges to some positive
constant. His proof heavily depends on the normality
assumption. Techniques developed in those papers are not applicable
here since we have {\it only one realization} and the dependence
structure among the entries can be quite complicated.

\subsection{A Summary of Results of the Paper}

We present the main results in Section~\ref{sec:main}, which
include a central limit theory of (\ref{eq:hong}) and the Gumbel
convergence (\ref{eq:main1}). The proofs are given in
Section~\ref{sec:proof}. In Section~\ref{sec:normalcomparison} we
prove a normal comparison principle, which is of independent
interest. Since summability conditions of joint cumulants are
commonly used in time series analysis (see for example
\cite{brillinger:2001} and \cite{rosenblatt:1985}) and is needed
in the proof of Theorem~\ref{thm:ljung}, we present a sufficient
condition in Section~\ref{sec:cum}. Some auxiliary lemmas are
collected in Section~\ref{sec:aux}. We also conduct a simulation
study in Section~\ref{sec:simulation}, where we design a
simulation-based block of blocks bootstrapping procedure that
improves the finite-sample performance.

\section{Main Results}
\label{sec:main}

To develop an asymptotic theory for time series, it is necessary
to impose suitable measures of dependence and structural
assumptions for the underlying process $(X_i)$. Here we shall
adopt the framework of \cite{wu:2005}. Assume that $(X_i)$ is a
stationary causal process of the form
\begin{eqnarray}
  \label{eq:wold}
  X_i = g(\cdots,\epsilon_{i-1},\epsilon_{i}),
\end{eqnarray}
where $\epsilon_i, i \in \Z$, are iid random variables, and $g$ is
a measurable function for which $X_i$ is a properly defined random
variable. For notational simplicity we define the operator
$\Omega_k$: suppose $X = h(\epsilon_j, \epsilon_{i-1}, \ldots)$ is
a random variable which is a function of the innovations
$\epsilon_l, l\le j$, then $\Omega_k(X) := h(\epsilon_j, \ldots,
\epsilon_{k+1}, \epsilon_k', \epsilon_{k-1}, \ldots)$, where
$(\epsilon_k')_{k \in \Z}$ is an iid copy of $(\epsilon_k)_{k \in
\Z}$. Namely $\epsilon_k$ in $X$ is replaced by $\epsilon_k'$.

For a random variable $X$ and $p > 0$, we write $X \in
\mathcal{L}^p$ if $\|X\|_p:=(\E |X|^p )^{1/p}<\infty$, and in
particular, use $\|X\|$ for the $\mathcal{L}^2$-norm $\|X\|_2$.
Assume $X_i \in \mathcal{L}^p$, $p > 1$. Define the {\it physical
dependence measure of order $p$} as
\begin{equation}
  \label{eq:physical}
  \delta_p(i)=\|X_i-\Omega_0(X_i)\|_p,
\end{equation}
which quantifies the dependence of $X_i$ on the innovation
$\epsilon_0$. Our main results depend on the decay rate of
$\delta_p(i)$ as $i\to\infty$. Let $p'=\min(2,p)$ and define
\begin{eqnarray}
  \label{eq:stable}
  \Theta_{p}(n)&=&\sum_{i=n}^{\infty}\delta_p(i),
  \quad
  \Psi_{p}(n)=\left(\sum_{i=n}^{\infty}\delta_{p}(i)^{p'}\right)^{1/p'},
  \quad\hbox{and}\quad \cr
  \Delta_p(n)&=&\sum_{i=0}^{\infty}\min\{\CC_p\Psi_{p}(n),\delta_p(i)\},
\end{eqnarray}
where $\CC_p$ is defined in (\ref{eq:burkholder}). It is easily
seen that $\Psi_p(\cdot) \le \Theta_p(\cdot) \le \Delta_p(\cdot)$.
We use $\Theta_p$, $\Psi_p$ and $\Delta_p$ as shorthands for
$\Theta_p(0)$, $\Psi_p(0)$ and $\Delta_p(0)$ respectively. We make
the convention that $\delta_p(k)=0$ for $k<0$.

There are several reasons that we use the framework (\ref{eq:wold})
and the dependence measure (\ref{eq:physical}).  First, the class of
processes that (\ref{eq:wold}) represents is huge and it includes
linear processes, bilinear processes, Volterra processes, and many
other time series models. See, for instance, \cite{tong:1990} and
\cite{wiener:1958}. Second, the physical dependence measure is easy to
work with and it is directly related to the underlying data-generating
mechanism.  Third, it enables us to develop an asymptotic theory for
complicated statistics of time series.

\subsection{Maximum deviations of sample autocovariances}

Note that $\hat\gamma_k$ is a biased estimate of $\gamma_k$ with
$\E \hat\gamma_k = (1-|k|/n) \gamma_k$. It is then more convenient
to consider the centered version $\max_{1\le k\le s_n} \sqrt{n}
|\hat\gamma_k - \E \hat\gamma_k|$ instead of $\max_{1\le k\le s_n}
\sqrt{n} |\hat\gamma_k - \gamma_k|$. Recall
(\ref{eq:gumbel_constants}) for $a_n$ and $b_n$.

\begin{theorem}
\label{thm:single} Assume $\E X_i=0$, $X_i \in \mathcal{L}^p$ for
some $p>4$, and $\Theta_p(m)=O(m^{-\alpha})$, $\Delta_p(m) =
O(m^{-\alpha'})$ for some $\alpha\geq\alpha'>0$.  If $s_n$
satisfies $s_n \to \infty$ and $s_n = O(n^\eta)$ with
  \begin{align}
    \label{eq:decay_rate}
    0<\eta<1, \quad \eta<\alpha p/2, \,\,\, \hbox{and} \,\,\,
    \eta \min\{2(p-2-\alpha p),\, (1-2\alpha')p\}<p-4,
  \end{align}
  % $s_n=O(n^\beta)$ for some
  % $\beta<\min\left\{1,{\alpha p}/{2}\right\}{(p-4)}/{p}$,
  % \begin{equation}
  %   \label{eq:size}
  %   \limsup \log(s_n)/\log(n) <\min\left\{1,
  %   \frac{\alpha p}{2}\right\}\frac{p-4}{p},
  % \end{equation}
then for all $x \in \R$,
\begin{align}\label{eq:main}
\lim_{n\to\infty}P\left(\max_{1\leq k\leq
s_n}|\sqrt{n}[\hat{\gamma}_k-(1-k/n)\gamma_k]|
 \leq \sqrt{\sigma_0}(a_{2s_n}\,x+b_{2s_n})\right)
 = \exp\{-\exp(-x)\}.
  \end{align}
\end{theorem}

In (\ref{eq:decay_rate}), if $p \le 2 + \alpha p$ or $1 \le 2
\alpha'$, then the second and third conditions are automatically
satisfied, and hence Theorem \ref{thm:single} allows a very wide range
of lags $s_n = O(n^\eta)$ with $0 < \eta < 1$. In this sense Theorem
\ref{thm:single} is nearly optimal.

% \begin{corollary}
%   Theorem~\ref{thm:single} and Theorem~\ref{thm:singlegmc} hold for
%   $\max_{1\leq k\leq s_n}|\sqrt{n}(\check\gamma_k-\gamma_k)|$. They
%   also hold for $\max_{1\leq k\leq
%     s_n}\sqrt{n}|\hat\gamma_k-\gamma_k|$ if the following condition is satisfied
%   \begin{equation}
%     \label{eq:30}
%     \max_{1\leq k\leq s_n}|\sqrt{n}(\E\hat\gamma_k-\gamma_k)|
%     \leq n^{-1/2}\max_{1\leq k\leq s_n}|k\gamma_k| = o\left((\log s_n)^{-1/2}\right).
% \end{equation}
% \end{corollary}
% \begin{proof}
%   The second assertion is obvious. To see the first one, observe that
%   \begin{align*}
%     \max_{1\leq k\leq s_n}|\sqrt{n}(\hat\gamma_k-\E\hat\gamma_k-\check\gamma_k+\gamma_k)|
%     = \frac{s_n}{n-s_n}\max_{1\leq k\leq s_n}|\sqrt{n}(\hat\gamma_k-\E\hat\gamma_k)|
%     = o_P\left(\frac{1}{\sqrt{\log s_n}}\right).
%   \end{align*}
% \end{proof}
% \begin{remark}
%   Condition (\ref{eq:30}) is very mild. For example, it holds if
%   $|\gamma_k|=O(1/k)$, which in turn is implied by
%   $\delta_{2}(k)=O(1/k)$, see (\ref{eq:fact4}).
% \end{remark}

For the maximum deviation $\max_{1\leq k<n} | \hat\gamma_k - \E
\hat\gamma_k|$ over the whole range $1\leq k < n$, it seems not
possible to derive a limiting distribution by using our method.
However, we can obtain a sharp bound $(n^{-1} \log n)^{1/2}$. The
upper bound is given in (\ref{eq:covorder}), while the lower
bounded can be obtained by applying Theorem \ref{thm:single} and
choosing a sufficiently small $\eta$ such that
(\ref{eq:decay_rate}) holds. Using Theorem \ref{thm:covorder},
\cite{wu:2010} derived convergence rates for the thresholded
autocovariance matrix estimates.

\begin{theorem}
  \label{thm:covorder}
  Assume $\E X_i=0$, $X_i \in \mathcal{L}^p$ for some $p>4$, and
  $\Theta_p(m)=O(m^{-\alpha})$, $\Delta_p(m)=O(m^{-\alpha'})$ for some
  $\alpha\geq\alpha'>0$. If
  \begin{align}
    \label{eq:decay_rate2}
    \alpha>1/2 \quad\hbox{or}\quad \alpha'p>2
  \end{align}
  then for $c_p=6(p+4)\,e^{p/4}\,\kappa_4\,\Theta_4$,
  \begin{eqnarray}
    \label{eq:covorder}
    \lim_{n\to\infty}
    P\left(\max_{1\leq k<n}|\hat\gamma_k-\E\hat\gamma_k|
      \leq c_p\sqrt{\frac{\log n}{n}}\right) =1.
  \end{eqnarray}
\end{theorem}

Since $\Theta_p(m) \ge \Psi_p(m)$, we can assume $\alpha \ge
\alpha'$. For a detailed discussion on their relationship, see
Remark~6 of \cite{wu:2010}. It turns out that for the special case
of linear processes the condition (\ref{eq:decay_rate}) can be
weakened to the following one:
\begin{align}
  \label{eq:decay_rate1}
  0 < \eta <1, \quad \eta<\alpha p/2,
  \quad \hbox{and} \quad (1-2\alpha)\eta<(p-4)/p.
\end{align}
See Remark~\ref{rk:decay_rate}. Furthermore, for linear processes
the condition (\ref{eq:decay_rate2}) can be relaxed to $\alpha
p>2$ as well.

In practice, the mean $\mu = \E X_0$ is often unknown and we can
estimate it by the sample mean $\bar X_n = (1/n)\sum_{i=1}^n X_i$.
The usual estimates of autocovariances and autocorrelations are
\begin{align}\label{eq:D26175}
  \breve\gamma_k=\frac{1}{n}
  \sum_{i=k+1}^n(X_{i-k} - \bar X_n)(X_i - \bar X_n)
  \quad \hbox{and} \quad
  \breve r_k = \breve \gamma_k / \breve \gamma_0.
\end{align}

\begin{corollary}
  \label{thm:single_corr}
Theorem~\ref{thm:single} and Theorem~\ref{thm:covorder} still hold
if we replace $\hat\gamma_k$ therein by $\breve\gamma_k$.
Furthermore,
  \begin{align*}%\label{eq:main_corr}
    \lim_{n\to\infty}
    P\left(\max_{1\leq k\leq s_n}
    \left|\sqrt{n}[\breve r_k-(1-k/n)r_k]\right|
      \leq (\sqrt{\sigma_0}/\gamma_0)(a_{2s_n}\,x+b_{2s_n})\right)
      = \exp\{-\exp(-x)\}.
  \end{align*}
\end{corollary}

\begin{proof}[Proof of Corollary~\ref{thm:single_corr}]
For the $\breve\gamma_k$ version of Theorem~\ref{thm:single}, it
suffices to show that
  \begin{align}
    \label{eq:38}
    \max_{1\leq k\leq s_n}
    \left|\sqrt{n}(\breve\gamma_k-\hat\gamma_k)\right|
    = o_P\left(\frac{1}{\sqrt{\log s_n}}\right).
  \end{align}
  Let $S_k=\sum_{i=1}^k X_i$. By Theorem~1 (iii) of \cite{wu:2007}, we
  have $\left\|\max_{1\leq k \leq n}\left|S_k\right|\right\| \leq
  2\sqrt{n}\Theta_2$. Since
  \begin{align*}
    \sum_{i=k+1}^n (X_{i-k}-\bar X_n)(X_i-\bar X_n) - \sum_{i=k+1}^nX_{i-k}X_i
    & = -\bar X_n\sum_{i=1}^{n-k}X_i + \bar X_n\sum_{i=1}^k X_i - k\bar
    X_n^2,
  \end{align*}
  we have (\ref{eq:38}). The proof of the $\breve\gamma_k$ version of
  Theorem~\ref{thm:covorder} is similar. The assertion on sample
  autocorrelations can be proved easily, and details are omitted.
\end{proof}

% Theorem~\ref{thm:single} and \ref{thm:singlegmc} cannot be extend

%\subsection{Stationary processes with iid repetitions}

\subsection{Box-Pierce tests}
Box-Pierce tests \citep{box:1970,ljung:1978} are commonly used in
detecting lack of fit of a particular time series model. After a
correct model has been fitted to a set of observations, one would
expect the residuals to be close to a sequence of iid random
variables, and therefore one should perform some tests for serial
correlations as model diagnostics. Suppose $(X_i)_{1\leq i\leq n}$
is an iid sequence, let $\hat r_k$ be its sample autocorrelations.
Then the distribution of $Q_n(K) := n \sum_{k=1}^K\hat r_K^2$ is
approximately $\chi^2_K$. Logically, it is not sufficient to
consider a fixed number of correlations as the number of
observations increases, because there may be some dependencies at
large lags. We present a normal theory about the Box-Pierce test
statistic, which allows the number of correlations included in
$Q_n$ to go to infinity.
% \cite{ljung:1978} suggested that the distribution of
% \begin{align*}
%   \label{eq:ljung}
%   n(n+2) \sum_{k=1}^m (n-k)^{-1}\hat r_k^2
% \end{align*}

\begin{theorem}
  \label{thm:ljung}
  Assume $X_i \in \LLL^8$, $\E X_i=0$ and $\sum_{k=0}^\infty
  k^6\delta_8(k)<\infty$. If $s_n \to \infty$ and $s_n
  =O(n^\beta)$ for some $\beta<1$, then
  \begin{align*}
    \frac{1}{\sqrt{s_n}}\sum_{k=1}^{s_n}
    \left[n(\hat\gamma_k - (1-k/n)\gamma_k)^2 - (1-k/n)\sigma_0\right]
    \Rightarrow \mathcal{N}\left(0,2\sum_{k\in\Z}\sigma_k^2\right).
  \end{align*}
\end{theorem}
To see the connection to the Box-Pierce test, we have the
following corollary on autocorrelations. Using the same argument,
we can show that the same asymptotic law holds for the similar
Ljung-Box test statistic $Q_{L B} = n(n+2) \sum_{k=1}^K \hat r_K^2
/ (n-k)$.

\begin{corollary}
 \label{thm:ljung_corr}
Under the conditions of Theorem~\ref{thm:ljung}, the same result
holds if $\hat \gamma_k$ is replaced by $\breve \gamma_k$.
Furthermore,
  \begin{align}
    \label{eq:ljung_corr}
    \frac{1}{\sqrt{s_n}}\sum_{k=1}^{s_n}
    \left[n(\hat r_k - (1-k/n)r_k)^2
        - (1-k/n)\sigma_0/\gamma_0^2\right]
    \Rightarrow \mathcal{N}\left(0,
    \frac{2}{\gamma_0^4}\sum_{k\in\Z}\sigma_k^2\right).
  \end{align}
\end{corollary}

\begin{remark}
  \label{rk:chisq}
  Theorem \ref{thm:ljung} clarifies an important historical issue in
  testing of correlations. If $\gamma_k = 0$ for all $k\geq 1$, which
  means $X_i$ are uncorrelated; then $\sigma_0 = \gamma_0^2$ and
  $\sigma_k = 0$ for all $|k| \ge 1$, and (\ref{eq:ljung_corr})
  becomes
  \begin{align}
    \label{eq:ljung_uncorr}
    \frac{1}{\sqrt{s_n}}\sum_{k=1}^{s_n}
    \left[n\hat r_k^2 - (1-k/n)\right]
    \Rightarrow \mathcal{N}\left(0,2\right).
  \end{align}
  In an influential paper, \cite{romano:1996} argued that, for fixed
  $K$, the chi-squared approximation for $Q_n(K)$ does not hold if
  $X_i$ are only uncorrelated but not independent. One of the main
  reasons is that for fixed $K$, $\hat r_1, \ldots, \hat r_K$ are not
  asymptotically independent if $X_i$ are not independent.  However,
  interestingly, the situation is different if the number of
  correlations included in $Q_n$ can increase to infinity.  According
  to (\ref{eq:cltunbdd}), $\sqrt n \hat \gamma_{k_n}$ and $\sqrt n
  \hat \gamma_{k_n+h}$ are asymptotically independent if $h > 0$ and
  $k_n \to \infty$, because the asymptotic covariance is $\sigma_h =
  0$. Therefore, the original Box-Pierce approximation of $Q_n(s_n)$
  by $\chi^2_{s_n}$, with unbounded $s_n$, is still asymptotically
  valid in the sense of (\ref{eq:ljung_uncorr}) since $(\chi^2_{s_n} -
  s_n) / \sqrt{s_n} \Rightarrow \mathcal{N}\left(0,2\right)$ as $s_n
  \to \infty$. This observation again suggests that the asymptotic
  behaviors for bounded and unbounded lags are different. A similar
  observation has been made in \cite{shao:2011}, whose result also
  suggests that (\ref{eq:ljung_uncorr}) is true under the assumption
  that $\delta_8(k) = O(\rho^k)$ for some $0<\rho<1$. Our condition
  $\sum_{k=1}^\infty k^6\delta_8(k)<\infty$ is much weaker.

%at least when $s_n=o(n^{2/3})$. They suggested to use resampling
%and subsampling scheme to obtain the distribution of $Q_n$.
\end{remark}

The next theorem consists of two separate but closely related
parts, one is on the estimation of $\sigma_0 = \sum_{k \in \Z}
\gamma_k^2$, and the other is related to the power of the
Box-Pierce test. Define the projection operator
\begin{eqnarray*}
\PP^j \cdot =  \E(\cdot|\FF_{-\infty}^j) - \E(\cdot |
\FF_{-\infty}^{j-1}), \mbox{ where } \FF_i^j = \langle \epsilon_i,
\epsilon_{i+1}, \ldots, \epsilon_j \rangle, \, i, j \in \Z.
\end{eqnarray*}

\begin{theorem}
  \label{thm:ljung_power}
  Assume $X_i\in\LLL^4$, $\E X_i=0$ and $\Theta_4<\infty$. If
  $s_n\to\infty$ and $s_n=o(\sqrt{n})$, then
  \begin{align}
    \label{eq:var_estimation}
    \sqrt{n}\left(\sum_{k=-s_n}^{s_n} \hat\gamma_k^2
    - \sum_{k=-s_n}^{s_n} \gamma_k^2\right)
    \Rightarrow \mathcal{N}(0,4\|D_0'\|^2),
  \end{align}
  where $D'_0=\sum_{i=0}^\infty \PP^0 (X_iY_i)$ with $Y_i=\gamma_0 X_i
  + 2\sum_{k=1}^\infty\gamma_k X_{i-k}$.
  Furthermore, if $\sum_{k=1}^\infty \gamma_k^2>0$, then
  \begin{align}
    \label{eq:ljung_power}
    \sqrt{n}\left(\sum_{k=1}^{s_n} \hat\gamma_k^2
    - \sum_{k=1}^{s_n} \gamma_k^2\right)
    \Rightarrow \mathcal{N}(0,4\|D_0\|^2),
  \end{align}
  where $D_0=\sum_{i=0}^\infty \PP^0 (X_iY_i)$ with
  $Y_i=\sum_{k=1}^\infty\gamma_kX_{i-k}$.
\end{theorem}

\begin{corollary}
  \label{thm:ljung_corr_power}
Under conditions of Theorem~\ref{thm:ljung_power}, the same
results hold if $\hat \gamma_k$ is replaced by $\breve \gamma_k$.
Furthermore, there exist positive numbers $\tau_1^2$ and
$\tau_2^2$ such that
  \begin{align*}
    % \label{eq:ljung_corr_power}
    \sqrt{n}\left(\sum_{k=1}^{s_n} \hat r_k^2
    - \sum_{k=1}^{s_n} r_k^2\right)
    \Rightarrow \mathcal{N}(0,\tau_1^2)
    \quad \hbox{and} \quad
    \sqrt{n}\left(\sum_{k=-s_n}^{s_n} \hat r_k^2
    - \sum_{k=-s_n}^{s_n} r_k^2\right)
    \Rightarrow \mathcal{N}(0,\tau_2^2).
  \end{align*}
\end{corollary}

As an immediate application, we consider testing whether $(X_i)$
is an uncorrelated sequence. According to (\ref{eq:ljung_uncorr}),
we can use the test statistic
\begin{align*}
  T_n:=\frac{1}{\sqrt{s_n}}
   \left[Q_n(s_n)-\frac{s_n(2n-s_n-1)}{2n}\right],
\end{align*}
whose asymptotic distribution under the null hypothesis is
$\mathcal{N}(0,2)$. The null is rejected when $T_n > \sqrt{2}
z_{1-\alpha}$, where $z_{1-\alpha}$ is the $(1-\alpha)$-th
quantile of a standard normal random variable $Z$. However, under
the alternative hypothesis $\sum_{k=1}^\infty r_k^2>0$, the
distribution of $T_n$ should be approximated according to
Corollary~\ref{thm:ljung_corr_power}, and the asymptotic power is
\begin{align*}
  P\left(T_n>\sqrt{2}z_{1-\alpha}\right)
  \approx P\left(\tau_1Z >
     \frac{\sqrt{2s_n}\cdot z_{1-\alpha}}{\sqrt{n}}
    + \frac{s_n(2n-s_n-1)}{2n^{3/2}}
    - \sqrt{n}\sum_{k=1}^{s_n}r_k^2\right),
\end{align*}
which increases to 1 as $n$ goes to infinity.

\section{A Simulation Study}
\label{sec:simulation}
%\subsection{Limiting distribution of }
%\label{sec:simulation-asymp}

Suppose $\left(r^{(0)}_k\right)$ is a sequence of
autocorrelations, one might be interested in the hypothesis test
that $r_k=r^{(0)}_k$ for all $k \geq 1$. This hypothesis is,
however, impossible to test in practice, except in some special
parametric cases. A more tractable hypothesis is
\begin{equation}
  \label{eq:hypothesis}
  \boldsymbol{\mathrm{H}}_0:\quad r_k=r^{(0)}_k
  \quad\hbox{for } 1 \leq k \leq s_n.
\end{equation}
In traditional asymptotic theory, one often assumes that $s_n$ is a
fixed constant, for example, the popular Box-Pierce test for serial
correlation. Our results in the previous section provide both ${\cal
  L}^{\infty}$ and ${\cal L}^2$ based tests, which allow $s_n$ to grow
as $n$ increases. Nonetheless, the asymptotic tests can perform poorly
when the sample size $n$ is not large enough, namely, there may exist
noticeable differences between the true and nominal probabilities of
rejecting $\boldsymbol{\mathrm{H}}_0$ (hereafter referred as error in
rejection probability or ERP). In a recent paper, \cite{horowitz:2006}
showed that the Box-Pierce test with bootstrap-based $p$-values can
significantly reduce the ERP. They used the blocks of blocks
bootstrapping with overlapping blocks (hereafter referred as BOB)
invented by \cite{kunsch:1989}.  For finite sample, our ${\cal L}^2$
based test is similar as the traditional Box-Pierce test considered in
their paper, so in this section our focus will be on the ${\cal
  L}^{\infty}$ based tests. We shall provide simulation evidence
showing that the BOB works reasonably well.

\begin{figure}[!htb]
\centering
\includegraphics[width=10cm]{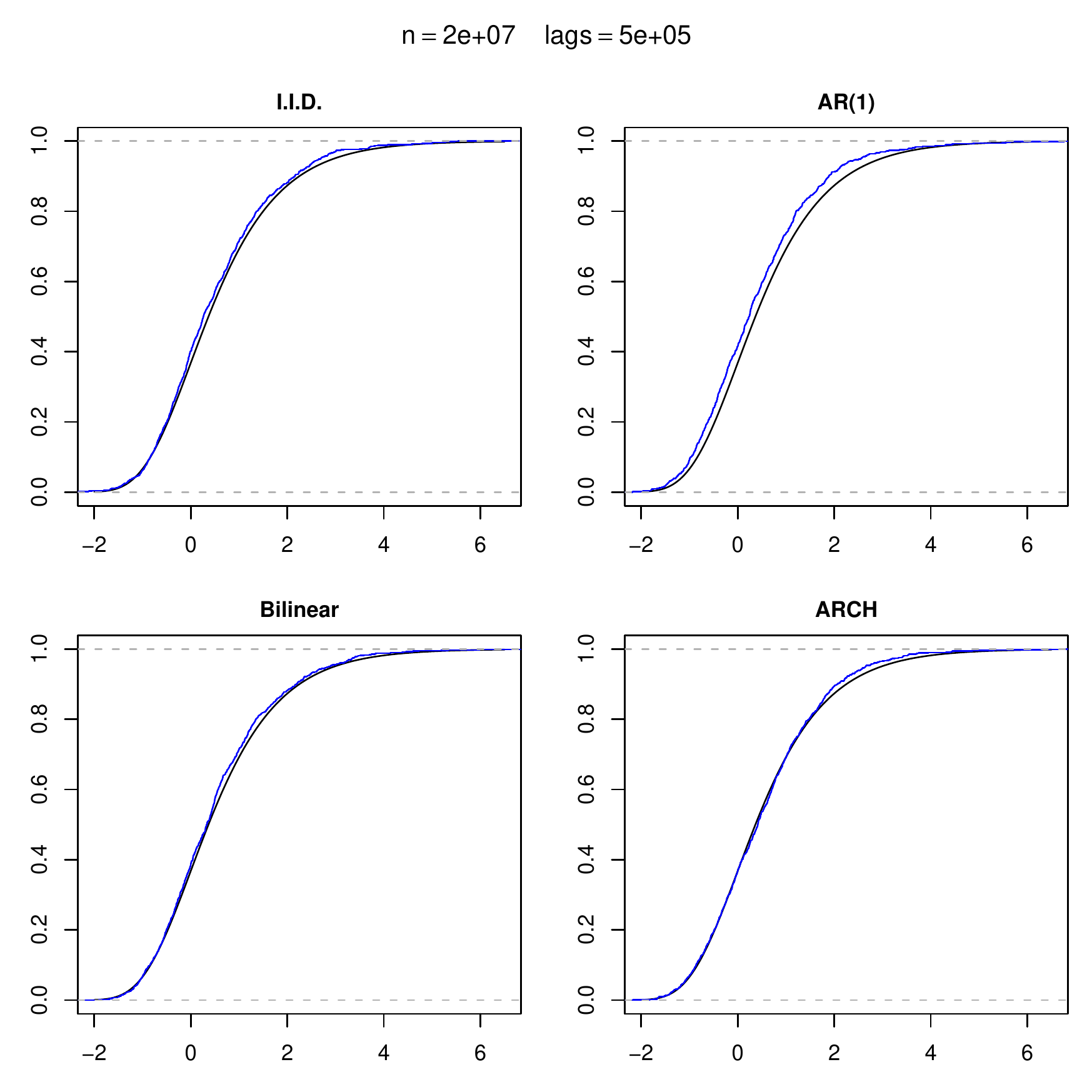}
\caption{\label{fig:all_in_one} Empirical distribution functions
for quantities in (\ref{eq:quantity_simulated}). {\footnotesize We
choose $b = 0.5$ for model (\ref{eq:ar}), $a = b = 0.4$ for model
(\ref{eq:bilin}), and $a = b = 0.25$ for model (\ref{eq:arch}).
The black line gives the true distribution function of the Gumbel
distribution. }}
\end{figure}

Throughout this section, we let the innovations $\epsilon_i$ be
iid standard normal random variables, and consider the following
four models.
\begin{align}
  & \hbox{I.I.D.: } & & X_i=\epsilon_i && \label{eq:iid}\\
  & \hbox{AR(1): } & & X_i=bX_{i-1}+\epsilon_i && \label{eq:ar}\\
  & \hbox{Bilinear: } & &
        X_i=(a+b\epsilon_i)X_{i-1}+\epsilon_i && \label{eq:bilin}\\
  & \hbox{ARCH: } & &
       X_i=\sqrt{a+bX_{i-1}^2} \cdot \epsilon_i. && \label{eq:arch}
\end{align}
We generate each process with length $n=2\times 10^7$, and compute
\begin{equation}
  \label{eq:quantity_simulated}
  a_{2s_n}^{-1}\left(\max_{1\leq k\leq s_n}\sqrt{n}
   \left|\hat r_k-(1-k/n)r_k\right|/\sqrt{\hat\sigma_0}-b_{2s_n}\right)
\end{equation}
with $s_n= 5\times 10^5$ and $\hat\sigma_0=\sum_{k=-t_n}^{t_n}\hat
r_k^2$, where $t_n$ is chosen as $t_n=\floor{n^{1/3}} = 271$. Based on
1000 repetitions, we plot the empirical distribution functions in
Figure~\ref{fig:all_in_one}. We see that these four empirical curves
are close to the one for the Gumbel distribution, which confirms our
theoretical results.

One the other hand, these empirical distributions are not very close
to the limiting one if the sample size is not large, because the
Gumbel type of convergence in (\ref{eq:main}) is slow. This is a
well-known phenomenon; see for example \cite{hall:1979}. It is
therefore not reasonable to use the limiting distribution to
approximate the finite sample distributions. To perform the test
(\ref{eq:hypothesis}), we repeat the BOB procedure as described in
\cite{horowitz:2006} (called SBOB in their paper). Since in the
bootstrapped tests, the test statistics are not to be compared with
the limiting distribution, we can ignore the norming constants in
(\ref{eq:quantity_simulated}) and simply use the following test
statistics
\begin{eqnarray*}
M_n = \max_{1\leq k\leq s_n} \left|r_k-(1-k/n)r_k^{(0)} \right|
 \mbox{ and }
 \mathcal M_n = {{M_n} \over { \sqrt{\hat\sigma_0}}},
\end{eqnarray*}
where $\mathcal M_n$ is the self-normalized version with
$\sigma_0$ estimated as $\hat \sigma_0 = \sum_{k=-t_n}^{t_n} \hat
r_k^2$, with $t_n = \min\{\floor{n^{1/3}}, s_n\}$. For simplicity,
we refer these two tests as $M$-test and $\mathcal M$-test,
respectively.
% More specifically, one application of our
% results is using (\ref{eq:quantity_simulated}) as a test statistic to
% test the hypothesis that $r_k=r_k^{(0)}$ for all $k \in \Z$. For a
% level $\alpha$ test, when $n$ is not very large, the $(1-\alpha)$-th
% quantile of the Gumbel distribution needs not be close to the one for
% $(\ref{eq:quantity_simulated})$. Therefore, finding the critical
% values is not a trivial problem. It is not the goal of this paper to
% explore this issue in great details. However, we shall use some
% simulation results to show that the blocks of blocks bootstrapping
% with overlapping blocks \citep{kunsch:1989} seems working well in some
% situations. This procedure was implemented by \cite{horowitz:2006} to
% bootstrap the Box-Pierce test.

From the series $X_1,\ldots,X_n$, for some specified number of
lags $s_n$ that will be included in the test and block size
$\mathfrak{b}_n$, form $Y_i=(X_i,X_{i+1},\ldots,X_{i+s_n})^\top$,
$1\le i \le n-s_n$ and blocks $\mathcal{B}_j = (Y_j, Y_{j+1},
\ldots, Y_{j+\mathfrak{b}_n-1})$, $1 \le j \le n - s_n -
\mathfrak{b}_n+1$. For simplicity assume $h_n = n/\mathfrak{b}_n$
is an integer. Suppose $Y_\sharp$ is obtained by sampling a block
$\mathcal{B}_\sharp$ from the set of blocks $\{\mathcal{B}_1,
\mathcal{B}_2, \ldots, \mathcal{B}_{n-s_n-\mathfrak{b}_n+1}\}$,
and then sampling a column from $\mathcal{B}_\sharp$, let
$\Cov_\sharp$ represent the covariance of the bootstrap
distribution of $Y_\sharp$, conditional on $(X_1,X_2,\ldots,X_n)$.
Denote by $Y_\sharp^j$ the $j$-th entry of $Y_\sharp$, set
\begin{equation*}
  r_k^{(e)}= {{\Cov_\sharp(Y_\sharp^1, Y_\sharp^{k+1})}
  \over
  \sqrt{\Cov_\sharp(Y_\sharp^1,Y_\sharp^1) \cdot
  \Cov_\sharp(Y_\sharp^{k+1},Y_\sharp^{k+1})}}.
\end{equation*}
The explicit formula of $r_k^{(e)}$ was also given in
\cite{horowitz:2006}. The BOB algorithm is as follows.
\begin{enumerate}
\item Sample $h_n$ times with replacement from $\{\mathcal{B}_1,
\mathcal{B}_2, \ldots, \mathcal{B}_{n-s_n-\mathfrak{b}_n+1}\}$ to
obtain blocks $\{\mathcal{B}^\ast_1, \mathcal{B}^\ast_2, \,
\ldots, \, \mathcal{B}^\ast_{h_n}\}$, which are laid end-to-end to
form a series of vectors $(Y^\ast_1, Y^\ast_2, \ldots, Y^\ast_n)$.
  % Denote by $Y_i^{\ast j}$ the
  % $j$-th entry of $Y^\ast_i$.

\item Pretend that $(Y^\ast_1,Y^\ast_2,\ldots,Y^\ast_n)$ is a
random sample of size $n$ from some $s_n$-dimensional population
distribution, let $r^\ast_k$ be the sample correlation of the
first entry and the $(k+1)$-th entry. Then calculate the test
statistic $M_n^\ast=\max_{1\leq k\leq s_n} \left|
r^\ast_k-r^{(e)}_k \right|$ and $\mathcal M_n^\ast =
M_n^\ast/\sqrt{\sigma_0^\ast}$, where $\sigma_0^\ast =
\sum_{k=-t_n}^{t_n} \left(r_k^\ast\right)^2$.

\item Repeat steps 1 and 2 for $N$ times. The bootstrap $p$-value
of the $M$-test is given by $\#(M_n^\ast>M_n)/N$. For a nominal
level $\alpha$, we reject $\boldsymbol{\mathrm{H}}_0$ if
$\#(M_n^\ast>M_n)/N<\alpha$. The $\mathcal M$-test is performed in
the same manner.
\end{enumerate}

We compare the BOB tests and the asymptotic tests for the four models
listed at the beginning of this section, with $a=.4$ for
(\ref{eq:ar}), $a=b=.4$ for (\ref{eq:bilin}) and $a=b=.25$ for
(\ref{eq:arch}). We set the series length as $n=1800$, and consider
four choices of $s_n$: $\floor{\log(n)}=7$, $\floor{n^{1/3}}=12$,
$\floor{\sqrt{n}}=42$ and 25. The BOB tests are performed with
$N=999$, and the asymptotic tests are carried out by comparing
$a_{2s_n}^{-1}(\sqrt{n}\mathcal M_n-b_{2s_n})$ with the corresponding
quantiles of the Gumbel distribution. The empirical rejection
probabilities based on 10,000 repetitions are reported in
Table~\ref{tab:sbob}. All probabilities are given in percentages. For
all cases, we see that the asymptotic tests are too conservative, and
the ERP are quite large. At the nominal level $1\%$, the rejection
probabilities are often less than or around $0.1\%$, and at most
$0.51\%$; while at nominal level $10\%$, they are often less than
$3\%$ and at most $6.4\%$. Except for the bilinear models with $s_n=7$
and $s_n=12$, the bootstrapped tests significantly reduce the ERP,
which are often less than $0.2\%$ at nominal level $1\%$, less than
$.5\%$ at level $5\%$, and less than $1\%$ at level $10\%$. The
performance of $M$-test and $\mathcal M$-test are similar, with the
former being slightly more conservative. The BOB tests are roughly
insensitive to the block size, which provides additional evidence of
the findings on BOB tests in \cite{davison:1997}.

{\small
\ctable[
  caption = {\small Empirical rejection probabilities (in percentages)},
  label = tab:sbob,
  pos = ht,
%  doinside = \small,
]{lrrrcrrrcrrrcrrr}{ \tnote[]{The values 1, 5, 10 in the 2nd row
indicate nominal levels in percentages. The numbers in the third
row starting with the model name ``I.I.D.'' are for the asymptotic
tests. The fourth row staring with $\mathfrak b_n=5$ is for BOB
$M$-tests with block size 5. The fifth row is for BOB $\mathcal
M$-tests with the same block size 5. Other rows should be read
similarly.} } { \hline
    Test & \multicolumn{3}{c}{$s_n=7$} && \multicolumn{3}{c}{$s_n=12$} && \multicolumn{3}{c}{$s_n=25$} && \multicolumn{3}{c}{$s_n=42$} \\ \cline{2-4} \cline{6-8} \cline{10-12} \cline{14-16}
& 1 & 5 & 10 &&  1 & 5 & 10 &&  1 & 5 & 10 &&  1 & 5 & 10  \\ \hline
I.I.D.               & .00 & .34 & 1.6 && .02 & .69 & 2.3 && .03 & .93 & 3.2 && .04 & 1.0 & 3.3 \\
$\mathfrak{b}_n=5$   & 1.3 & 5.1 & 10.0 && 1.1 & 5.2 & 9.8 && .95 & 4.7 & 9.3 && 1.0 & 4.7 & 9.6  \\
                     & 1.4 & 5.3 & 10.4 && 1.2 & 5.6 & 10.5 && 1.1 & 5.1 & 10.1 && 1.1 & 5.1 & 10.2  \\
$\mathfrak{b}_n=10$  & .83 & 4.8 & 10.0 && 1.1 & 4.9 & 9.6 && 1.1 & 4.9 & 10.1 && .65 & 4.3 & 8.9  \\
                     & .94 & 5.1 & 10.3 && 1.2 & 5.4 & 10.3 && 1.1 & 5.5 & 11.0 && .78 & 4.7 & 9.6  \\\hline
AR(1)                & .01 & .17 & 1.2 && .01 & .36 & 1.8 && .02 & .77 & 2.5 && .02 & .88 & 2.8  \\
$\mathfrak{b}_n=10$  & 1.3 & 5.7 & 10.9 && 1.3 & 5.5 & 11.4 && 1.3 & 5.5 & 10.9 && 1.1 & 5.7 & 11.2  \\
                     & 1.3 & 5.7 & 11.2 && 1.4 & 5.9 & 11.7 && 1.3 & 6.0 & 11.5 && 1.2 & 6.0 & 11.7  \\
$\mathfrak{b}_n=20$  & .98 & 5.5 & 10.9 && 1.0 & 5.8 & 11.3 && 1.1 & 5.3 & 10.6 && .86 & 4.9 & 10.5  \\
                     & 1.0 & 5.7 & 11.0 && 1.1 & 6.1 & 11.9 && 1.2 & 5.6 & 11.0 && .83 & 5.0 & 10.9  \\\hline
Bilinear             & .34 & 2.8 & 6.4 && .43 & 2.5 & 5.8 && .51 & 2.5 & 5.9 && .40 & 2.8 & 5.9  \\
$\mathfrak{b}_n=10$  & 2.8 & 8.7 & 14.4 && 1.8 & 7.1 & 12.7 && 1.2 & 6.1 & 12.0 && 1.2 & 5.4 & 10.9  \\
                     & 2.7 & 8.6 & 14.5 && 1.8 & 7.3 & 12.9 && 1.3 & 6.2 & 12.2 && 1.1 & 5.5 & 11.1  \\
$\mathfrak{b}_n=20$  & 2.7 & 8.4 & 14.6 && 2.1 & 7.2 & 13.5 && 1.5 & 6.3 & 12.0 && 1.3 & 5.2 & 10.8  \\
                     & 2.5 & 8.3 & 14.6 && 2.1 & 7.5 & 13.9 && 1.5 & 6.2 & 12.0 && 1.2 & 5.3 & 10.9  \\\hline
ARCH                 & .05 & .82 & 3.2 && .06 & 1.5 & 3.9 && .09 & 1.3 & 4.0 && .12 & 1.4 & 4.4  \\
$\mathfrak{b}_n=10$  & .99 & 5.0 & 10.5 && 1.2 & 4.9 & 9.7 && .80 & 4.6 & 9.9 && .82 & 4.7 & 9.3  \\
                     & 1.1 & 5.4 & 10.9 && 1.4 & 5.3 & 10.4 && .92 & 5.1 & 10.7 && .94 & 5.1 & 10.2  \\
$\mathfrak{b}_n=20$  & .86 & 5.1 & 10.5 && 1.0 & 5.0 & 10.3 && .69 & 4.8 & 9.7 && .63 & 4.3 & 8.9  \\
                     & .98 & 5.5 & 11.0 && 1.2 & 5.6 & 11.0 && .89 & 5.1 & 10.4 && .76 & 4.7 & 9.5  \\
\hline
}
}

% The bootstrapped test still performs poorly for bilinear model
% (\ref{eq:bilin}) with $a=.7$ and $b=.4$. This is possibly due to the
% fact that this bilinear process does not have high-order moments. A
% result in \cite{tong:1981} indicated that the 6-th order moment does
% not exist.
The bootstrapped tests still perform relatively poorly for
bilinear models when $s_n$ is small (7 and 12). This is possibly
due to the heavy-tailedness of the bilinear process.
\cite{tong:1981} gave necessary conditions for the existence of
even order moments. On the other hand, \cite{horowitz:2006} showed
that the iterated bootstrapping further reduce the ERP. It is of
interest to see whether the iterated procedure has the same effect
for the ${\cal L}^{\infty}$ based tests, in particular, whether it
makes the ERP reasonably small for the bilinear models when $s_n$
is small. The simulation for the iterated bootstrapping will be
computationally expensive and we do not pursue it here.

\section{Proofs}
\label{sec:proof} This section provides proofs for the results in
Section \ref{sec:main}. For readability we list the notation here.
For a random variable $X$, write that $X \in \mathcal{L}^p$,
$p>0$, if $\|X\|_p := (\E |X|^p )^{1/p} < \infty$. Write $\|X\| =
\|X\|_2$ if $p = 2$. To express centering of random variables
concisely, we define the operator $\E_0$ as $\E_0 X := X - \E X$.
For a vector $\h{x} = (x_1,\ldots,x_d)^\top \in \R^d$, let
$|\h{x}|$ be the usual Euclidean norm, $|\h{x}|_{\infty} :=
\max_{1 \le i \le d}|x_i|$, and $|\h{x}|_{\bullet} := \min_{1\leq
i \leq d}|x_i|$. For a square matrix $A$, $\rho(A)$ denotes the
operator norm defined by $\rho(A) := \max_{|\h{x}|=1} |A \h{x}|$.
Let us make some convention on the constants. We use $C$, $c$ and
$\mathcal{C}$ for constants. The notation $\CC_p$ is reserved for
the constant appearing in Burkholder's inequality, see
(\ref{eq:burkholder}). The values of $C$ may vary from place to
place, while the value of $c$ is fixed within the statement and
the proof of a theorem (or lemma). A constant with a symbolic
subscript is used to emphasize the dependence of the value on the
subscript.

The framework (\ref{eq:wold}) is particularly suited for two
classical tools for dealing with dependent sequences, martingale
approximation and $m$-dependence approximation.  For $i\leq j$,
define $\FF_i^j = \langle \epsilon_i, \epsilon_{i+1}, \ldots,
\epsilon_j \rangle$ be the $\sigma$-field generated by the
innovations $\epsilon_i, \epsilon_{i+1}, \ldots, \epsilon_j$, and
the projection operator $\HH_{i}^j (\cdot) = \E(\cdot|\FF_{i}^j)$.
Set $\FF_i := \FF_{i}^\infty$, $\FF^j := \FF_{-\infty}^j$, and
define $\HH_i$ and $\HH^j$ similarly. Define the projection
operator $\PP^j(\cdot) = \HH^j(\cdot) - \HH^{j-1}(\cdot)$, and
% $\PP_{j}(\cdot)= \E(\cdot|\FF_{j})-\E(\cdot|\FF_{j-1})$,
$\PP_i(\cdot) = \HH_i(\cdot) - \HH_{i+1}(\cdot)$, then
$(\PP^j(\cdot))_{j \in \Z}$ and $(\PP_{-i}(\cdot))_{i \in \Z}$
become martingale difference sequences with respect to the
filtrations $(\FF^j)$ and $(\FF_{-i})$, respectively. For $m\geq
0$, define $\tilde X_i=\HH_{i-m} X_i$, then $(\tilde X_i)_{i \in
\Z}$ is a $(m+1)$-dependent sequence.

\subsection{Some Useful Inequalities}
We collect in Proposition~\ref{thm:facts} some useful facts about
physical dependence measures and martingale and $m$-dependence
approximations. We expect that it will be useful in other
asymptotic problems that involve sample covariances. Hence for
convenience of other researchers, we provide explicit upper
bounds.

We now introduce a moment inequality (\ref{eq:mar_zyg}) which follows
from the Burkholder inequality \citep[see][]{burkholder:1988}. Let
$(D_i)$ be a martingale difference sequence and for every $i$, $D_i
\in \LLL^p$, $p > 1$, then
\begin{align}
  \label{eq:mar_zyg}
  \left\|D_1+D_2+\cdots+D_n\right\|_p^{p'} \leq \CC_p^{p'}
  \left(\|D_1\|_p^{p'} + \|D_2\|_p^{p'} + \cdots
   + \|D_n\|_p^{p'}\right),
\end{align}
where $p'=\min\{p,2\},$ and the constant
\begin{eqnarray}
  \label{eq:burkholder}
  \mathcal{C}_p = (p-1)^{-1} \hbox{ if } 1 < p < 2 \mbox{ and }
  = \sqrt{p-1} \hbox{ if } p \ge 2.
\end{eqnarray}
We note that when $p > 2$, the constant $\CC_p$ in (\ref{eq:mar_zyg})
equaled to $p-1$ in \cite{burkholder:1988}, and it was improved to
$\sqrt{p-1}$ by \cite{rio:2009}.

% Observe that for each $k \geq 0$, the sequence $(X_{i-k}X_i)_{i\in\Z}$
% also has a representation as (\ref{eq:wold}). Let
% $\delta_{k,p}(\cdot)$ be its physical dependence measures of order
% $p$. Define
% \begin{align*}

% \end{align*}
\begin{proposition}
  \label{thm:facts}
  \begin{enumerate}
  \item Assume $\E X_i=0$ and $p>1$. Recall that $p'=\min(p, 2)$.
    \begin{align}
      \label{eq:fact1}
      & \|\PP^0X_i\|_p \leq \delta_p(i) \quad\hbox{and}\quad
      \|\PP_0X_i\|_p \leq \delta_p(i)
      \\\label{eq:fact2}
      & \kappa_p:=\|X_0\|_p \leq \mathcal{C}_p\Psi_p\\\label{eq:fact3}
      & \left\|\sum_{i=1}^n c_iX_i\right\|_p
         \leq \mathcal{C}_p A_n\Theta_p, \mbox{ where }
         A_n = \left(\sum_{i=1}^n |c_i|^{p'}\right)^{1/p'} \\\label{eq:fact4}
      & |\gamma_k| \leq \zeta_2(k), \quad \hbox{where }
      \zeta_p(k):=\sum_{j=0}^\infty \delta_p(j)\delta_p(j+k) \\\label{eq:fact5}
      & \left\|\sum_{i=1}^n (X_{i-k}X_i -\gamma_k)\right\|_{p/2}
           \leq 2\CC_{p/2}\kappa_p\Theta_p\sqrt{n},
      \quad \hbox{when } p \geq 4 \\ \label{eq:fact5.5}
      & \left\|\sum_{i,j=1}^n c_{i,j} (X_iX_j-\gamma_{i-j})\right\|_{p/2}
      \leq 4\CC_{p/2}\CC_p \Theta_p^2 B_n \sqrt{n}, \quad \hbox{when } p \geq 4
      %& \delta_{k,p}(i) \leq \kappa_{2p}[\delta_{2p}(i)+\delta_{2p}(i-k)]
     \end{align}
     where
     $B_n^2 =\max \{\max_{1\leq i \leq n} \sum_{j=1}^n c_{i,j}^2,
        \, \max_{1\leq j \leq n} \sum_{i=1}^n c_{i,j}^2 \}.$
   \item For $m \geq 0$, define $\tilde X_i = \HH_{i-m} X_i$. For
     $p>1$, let $\tilde \delta_p(\cdot)$ be the physical dependence
     measures for the sequence $(\tilde X_i)$. Then
    \begin{align}
      \label{eq:fact6}
      & \tilde \delta_{p}(i) \leq \delta_p(i) \\ \label{eq:fact7}
      & \|X_0-\tilde X_0\|_p \leq \CC_p \Psi_p(m+1)\\ \label{eq:fact8}
      & \left\|\sum_{i=1}^n c_i(X_i-\tilde X_i)\right\|_p
         \leq \mathcal{C}_p A_n \Theta_{p}(m+1)\\
      \label{eq:fact9}
      & \left\|\sum_{i=k+1}^n \left(X_{i-k}X_i-\gamma_k
       -\tilde X_{i-k}\tilde X_i+\tilde\gamma_k\right)\right\|_p
     \leq 4\CC_p(n-k)^{1/p'} \kappa_{2p} \Delta_{2p}(m+1).
    \end{align}
  \end{enumerate}
\end{proposition}

\begin{proof}
The inequalities (\ref{eq:fact1}) and (\ref{eq:fact6}) are
obtained by the first principle. Since $X_{i-k} = \sum_{j\in\Z}
\PP^{j} X_{i-k}$ and $X_i = \sum_{j\in\Z}\PP^{j}X_i$, we have
\begin{eqnarray*}
  |\gamma_k|=\left|\sum_{j=-k}^\infty \E \left[(\PP^{-j}
      X_0) (\PP^{-j} X_{k})\right]\right| \leq
  \delta_2(j)\delta_2(j+k) \leq \zeta_{k},
\end{eqnarray*}
which proves (\ref{eq:fact4}). For (\ref{eq:fact5.5}), it can be
similarly proved as Proposition~1 of \cite{liu:2010}, and
(\ref{eq:fact8}) was given by Lemma~1 of the same paper.
(\ref{eq:fact3}) is a special case of (\ref{eq:fact8}). Define
$Y_i = X_{i-k} X_i$, then $(Y_i)$ is also a stationary process of
the form (\ref{eq:wold}). By H\"older's inequality, $\|Y_i -
\Omega_0(Y_i) \|_{p/2} \le 2 \kappa_p[\delta_p(i) +
\delta_p(i-k)]$. Applying (\ref{eq:fact3}) to $(Y_i)$, we obtain
(\ref{eq:fact5}). To see (\ref{eq:fact7}), we first write $X_m -
\tilde X_m = \sum_{j=1}^\infty \PP_{-j}X_m$.
  % $\tilde X_0 -X_0= \sum_{j=m}^\infty \left(\HH_{-j}^i X_0 -
  %   \HH_{-(j+1)}^i X_i\right)$. Since $\left\|\HH_{-j}^i X_0 -
  %   \HH_{-(j+1)}^i X_0\right\|_p \leq \delta_p(j+1)$, and
  % $\left(\HH_{-j}^i X_0 - \HH_{-(j+1)}^i X_0\right)_{j\geq m}$ is a
  % martingale difference sequence with respect to the filtration
  % $(\FF_{-(j+1),\infty})_{j\geq m}$, by Burkh\"older inequality, we
  % have
Since $\left\|\PP_{-j}X_m\right\|_p \leq \delta_p(m+j)$, and
$(\PP_{-j}X_m)_{j\geq 1}$ is a martingale difference sequence, by
(\ref{eq:mar_zyg}), we have
  \begin{equation*}
    \|X_0-\tilde X_0\|^{p'}_p \leq \mathcal{C}_p^{p'}
    \sum_{j=1}^{\infty}\left\|\PP_{-j}X_m\right\|_p^{p'}
    \leq \mathcal{C}_p^{p'} \sum_{j=1}^{\infty} [\delta_p(m+j)]^{p'}
    = \mathcal{C}_p^{p'} [\Psi_p(m+1)]^{p'}.
  \end{equation*}
The above argument also leads to (\ref{eq:fact2}). Using a similar
argument as in the proof of Theorem~2 of \cite{wu:2009}, we can
show (\ref{eq:fact9}). Details are omitted.
\end{proof}

\subsection{Proof of Theorem~\ref{thm:single}}
% We prove Theorem~\ref{thm:single} first.  It is a bit annoying that
% $\hat{\gamma}_k$ involves a sum of $X_{i-k}X_i$ over the range
% $k+1\leq i\leq n$, which is different for different $k$. So we define
% $R_{n,k}=\sum_{i=1}^n X_{i-k}X_i$. Based on the following lemma, we shall
% consider $R_{n,k}$ instead of $\hat{\gamma}_k$ in the sequel.
% \begin{lemma}
%   \label{thm:length}
%   Assume $X_i \in \mathcal{L}^p$ and $\Theta_{p}<\infty$ for some
%   $p\geq 4$. If $s_n=O(n^\beta)$ with $\beta<p/(p+4)$, then
%   \begin{equation*}
%     \label{eq:18}
%     \max_{1 \leq k \leq s_n} |R_{n,k}-n\hat{\gamma}_k-k\gamma_k|
%    =o_P\left(\sqrt{n/\log s_n}\right).
%   \end{equation*}
% \end{lemma}
% \begin{proof}
%   By (\ref{eq:fact3}) and (\ref{eq:fact5}), we have
%   $\|R_{n,k}-n\hat{\gamma}_k-k\gamma_k\|_{p/2} \leq
%   2\mathcal{C}_{p/2}\kappa_p\Theta_{p}\sqrt{k}$. Therefore, for any
%   $\epsilon>0$,
%   \begin{align*}
%     P & \left(\max_{1\leq k\leq s_n} |R_{n,k}-n\hat{\gamma}_k-k\gamma_k|
% > \epsilon\sqrt{\frac{n}{\log s_n}}\right)
%     \leq \sum_{k=1}^{s_n} \frac{(2\mathcal{C}_{p/2}\kappa_p
% \Theta_{p})^{p/2} k^{p/4} (\log s_n)^{p/4}}{\epsilon^{p/2} n^{p/4}}\cr
%     & \leq {C_p}{\epsilon^{-p/2}}{s_n^{p/4+1}}{n^{-p/4}}(\log s_n)^{p/4}
%     \leq {C_p}{\epsilon^{-p/2}} n^{\beta(p/4+1)-p/4}(\log s_n)^{p/4} = o(1).
%   \end{align*}
% \end{proof}

The proof is quite complicated and will be divided into several
steps. We first give the outline.

\subsubsection*{4.1.0. Outline}

\paragraph{\underline{\textnormal{\em Step 1: $m$-dependence approximation.}}}
Define $R_{n,k}=\sum_{i=k+1}^n (X_{i-k}X_i-\gamma_k)$. Set
$m_n=\floor{n^\beta}$, $0<\beta<1$. Define $\tilde X_i = \HH_{i-m_n}
X_i$, $\tilde \gamma_k = \E (\tilde X_0 \tilde X_k)$, and $\tilde
R_{n,k}=\sum_{i=k+1}^n (\tilde X_{i-k} \tilde X_i - \tilde \gamma_k
)$. We next show that it suffices to consider $\tilde R_{n,k}$.
\begin{lemma}
  \label{thm:mdep_app}
  Assume $\E X_i=0$, $X_i\in\LLL^{p}$, and $\Theta_p(m) =
  O(m^{-\alpha})$ for some $p>4$ and $\alpha>0$. If $s_n = O(n^\eta)$
  with $0<\eta<\alpha p/2$, then there exists a $\beta$ such that
  $\eta<\beta<1$ and
  \begin{align*}
    \max_{1\leq k\leq s_n} \left|  R_{n,k}-\tilde R_{n,k} \right|
    = o_P\left(\sqrt{n/\log s_n}\right).
  \end{align*}
\end{lemma}

%\noindent\underline{\em Step 2: Throw out small blocks.}\\[2mm]
\paragraph{\underline{\textnormal{\em Step 2: Throw out small blocks.}}}
Let $l_n = \floor{n^\gamma}$, where $\gamma \in (\beta, 1)$. For
each $t_n < k \leq s_n$, we apply the blocking technique and split
the integer interval $[k+1,n]$ into alternating large and small
blocks
\begin{equation}
  \label{eq:blk}
  \begin{aligned}
  & K_1=[k+1,s_n] \cr
  & H_j = [s_n+(j-1)(2m_n+l_n)+1,s_n+(j-1)(2m_n+l_n)+l_n];
       \quad 1 \leq j \leq w_n-1,\cr
  & K_{j+1}=[s_n+(j-1)(2m_n+l_n)+l_n+1,s_n+j(2m_n+l_n)];
   \quad 1\leq j \leq w_n-1; \quad \hbox{and} \cr
  & H_{w_n}=[s_n+(w_n-1)(2m_n+l_n)+1,n],
\end{aligned}
\end{equation}
where $w_n$ is the largest integer such that $s_n + (w_n-1)
(2m_n+l_n) + l_n \leq n$. Denote by $|H|$ the size of a block $H$.
By definition we know $l_n\leq |H_{w_n}| \leq 3l_n$ when $n$ is
large enough. For $1\leq j \leq w_n$ define
\begin{eqnarray*}
V_{k,j}=\sum_{i\in K_j,\,i>k} \left(\tilde X_{i-k} \tilde X_i -
\tilde \gamma_k\right)
 \mbox{ and }
U_{k,j}=\sum_{i\in H_j}
\left(\tilde X_{i-k} \tilde X_i - \tilde
  \gamma_k\right).
\end{eqnarray*}
Note that $w_n \sim n/(2 m_n+l_n) \sim n^{1-\gamma}$. We show that
the sums over small blocks are negligible.
\begin{lemma}
  \label{thm:small_blk}
  Assume the conditions of Theorem~\ref{thm:single}. Then
  \begin{align*}
    \max_{1 \leq k\leq s_n} \left|\sum_{j=1}^{w_n} V_{k,j} \right|
    = o_P\left(\sqrt{\frac{n}{\log s_n}}\right).
  \end{align*}
\end{lemma}

%\noindent\underline{\em Step 3: Truncate sums over large blocks.}\\[2mm]
\paragraph{\underline{\textnormal{\em Step 3: Truncate sums over large blocks.}}}
We show that it suffices to consider
\begin{eqnarray*}
\RR_{n,k} = \sum_{j=1}^{w_n} \bar U_{k,j}, \mbox{ where }
 \bar U_{k,j} = \E_0 \left(U_{k,j}I\{|U_{k,j}| \le \sqrt{n} /
(\log s_n)^3\}\right).
\end{eqnarray*}

\begin{lemma}
  \label{thm:large_blk}
  Assume the conditions of Theorem~\ref{thm:single}. Then
  \begin{align*}
    \max_{1\leq k\leq s_n} \left|\sum_{j=1}^{w_n}(U_{k,j}-\bar U_{k,j})\right|
    = o_P\left(\sqrt{\frac{n}{\log s_n}}\right).
  \end{align*}
\end{lemma}

%\noindent\underline{\em Step 4: Compare covariance structures.}\\[2mm]
\paragraph{\underline{\textnormal{\em Step 4: Compare covariance structures.}}}
In order to prove Lemma~\ref{thm:vec_mod_dev}, we need the
autocovariance structure of $\left(\RR_{n,k}/\sqrt{n}\right)$ to
be close to that of $(G_k)$. However, this only happens when $k$
is large. We show that there exists an $0 < \iota < 1$ such that
for $t_n = 3\floor{s_n^\iota}$, (i) $\max_{1\leq k\leq t_n}
|\RR_{n,k}/\sqrt{n}|$ does not contribute to the asymptotic
distribution; and (ii) the autocovariance structure of
$\left(\RR_{n,k}/\sqrt{n}\right)$ converges to that of $(G_k)$
uniformly on $t_n < k \leq s_n$.
\begin{lemma}
  \label{thm:small_lag}
  Under conditions of Theorem~\ref{thm:single}, there exists
  a constant $0<\iota<1$ such that for $t_n=3\floor{s_n^\iota}$,
  \begin{align}
    \label{eq:small_lag}
    \lim_{n\to\infty}P \left(\max_{1 \leq k \leq t_n}
    |\RR_{n,k}| > \sqrt{\sigma_0n\log s_n}\right)=0.
  \end{align}
\end{lemma}

\begin{lemma}
  \label{thm:cov_struc}
  Let conditions of Theorem~\ref{thm:single} be satisfied. Recall that
  $t_n=3\floor{s_n^\iota}$ from Lemma~\ref{thm:small_lag}. There exist
  constants $C_p>0$ and $0<\ell<1$ such that for any $t_n<k\leq
  k+h\leq s_n$,
    \begin{align*}%\label{eq:33}
      & |\Cov(\RR_{n,k},\RR_{n,k+h})/n-\sigma_h| \leq C_p\,s_n^{-\ell}.
    \end{align*}
\end{lemma}

\paragraph{\underline{\textnormal{\em Step 5: Moderate deviations.}}}
Let $t_n=3\floor{s_n^\iota}$ be as in Lemma~\ref{thm:small_lag}.
For $t_n<k_1<k_2<\ldots<k_d\leq s_n$, define $\h{\RR}_n =
(\RR_{n,k_1},\RR_{n,k_2},\ldots,\RR_{n,k_d})^\top$ and
$\h{V}=(G_{k_1},G_{k_2},\ldots,G_{k_d})^\top$, where $(G_k)$ is
defined in (\ref{eq:gp}). Let $\Sigma_n=\Cov(\h{\RR}_n)$ and
$\Sigma=\Cov(\h{V})$. For fixed $x \in \R$, set $z_n = a_{2s_n} x
+ b_{2s_n}$, where the constants $a_n$ and $b_n$ are defined in
(\ref{eq:gumbel_constants}). In the following lemma we provide a
moderate deviation result for $\h{\RR}_n$.

\begin{lemma}
  \label{thm:vec_mod_dev}
  Assume conditions of Theorem~\ref{thm:single}. Then there exists a
  constant $C_{p,d}>1$ such that for all $t_n<k_1<k_2<\ldots<k_d\leq
  s_n$,
  \begin{align*}
   \left|P\left(\left|\h{\RR}_n/\sqrt{n}\right|_\bullet \geq z_n \right)
    - P\left(\left|\h{V}\right|_\bullet \geq z_n \right)\right|
  \le C_{p,d} {{ P\left(\left|\h{V}\right|_\bullet \geq z_n \right)}
   \over{(\log s_n)^{1/2}}}
   + C_{p,d} \,\exp\left\{-{{(\log s_n)^2} \over{C_{p,d}}}\right\}.
  \end{align*}
\end{lemma}

\subsubsection{Step 1: $m$-dependence approximation}

\begin{proof}[Proof of Lemma~\ref{thm:mdep_app}]
  Recall that $m_n=\floor{n^\beta}$ with $\eta<\beta<1$. We claim
  \begin{align}
    \label{eq:mdep_app}
    \left\| R_{n,k} - \tilde R_{n,k}\right\|_{p/2}
    \leq 6\,\mathcal{C}_{p/2}\Theta_p\Theta_p(m_n-k+1)\cdot \sqrt{n}.
  \end{align}
  It follows that for any $\lambda>0$
  \begin{align*}
    P & \left( \max_{1\leq k\leq s_n} \left|  R_{n,k}-\tilde R_{n,k} \right|
      > \lambda \sqrt{n/\log s_n}\right)
    \leq {{ (\log s_n)^{p/4}}\over{n^{p/4} \lambda^{p/2}}}
    \sum_{k=1}^{s_n} \| R_{n,k}-\tilde R_{n,k}\|^{p/2}_{p/2} \\
    & \qquad \qquad \qquad \leq C_p \lambda^{-p/2}
           s_n (\log s_n)^{p/4} n^{-\alpha\beta p/2}
    \leq C_p \lambda^{-p/2} n^{\eta - \alpha\beta p/2} (\log n)^{p/4}.
  \end{align*}
  Therefore, if $\alpha p/2>\eta$, then there exists a $\beta$ such
  that $\eta<\beta<1$ and $\eta - \alpha\beta p/2 <0$, and hence the
  preceding probability goes to zero as $n\to\infty$. The
  proof of Lemma~\ref{thm:mdep_app} is complete.

We now prove claim (\ref{eq:mdep_app}). For each $1 \le k \le
s_n$, we have
  \begin{align*}
    \| R_{n,k} - \tilde R_{n,k} \|_{p/2}
    \leq & \left\|\sum_{i=k+1}^n ( X_{i-k} - \tilde X_{i-k})
               \tilde X_i\right\|_{p/2}
    + \left\|\sum_{i=k+1}^n (\HH_{i-m_n} X_{i-k})
                  ( X_i-\tilde X_i)\right\|_{p/2} \\
    & + \left\|\sum_{i=k+1}^n
     \E_0\left[( X_{i-k} - \HH_{i-m_n} X_{i-k})( X_i-\tilde X_i)
         \right]\right\|_{p/2}
  \end{align*}
Observe that $ (\tilde X_i\PP_{i-k-j} X_{i-k} )_{1\leq i \leq n}$
is a backward martingale difference sequence with respect to
$\FF_{i-k-j}$ if $j>m_n$, so by the inequality (\ref{eq:mar_zyg}),
  \begin{align*}
    \left\|\sum_{i=k+1}^n ( X_{i-k}
                    - \tilde X_{i-k})\tilde X_i\right\|_{p/2}
    & \leq \sum_{j=m+1}^\infty
     \left\|\sum_{i=k+1}^n \tilde X_i\PP_{i-k-j} X_{i-k}
        \right\|_{p/2} \cr
    & \leq \sum_{j=m+1}^\infty \sqrt{n}\mathcal{C}_{p/2}
      \|\tilde X_{j+k}\PP_0 X_j\|_{p/2} \cr
    & \leq \mathcal{C}_{p/2} \Theta_p\Theta_p(m_n+1)\cdot \sqrt{n}.
  \end{align*}
Similarly we have $\|\sum_{i=k+1}^n (\HH_{i-m_n} X_{i-k})( X_i -
\tilde X_i) \|_{p/2} \leq \sqrt{n} \mathcal{C}_{p/2} \Theta_p
\Theta_p(m_n+1)$. Similarly as (\ref{eq:fact8}), we get $\|\tilde
X_{i-k} - \HH_{i-m_n} X_{i-k}\|_p \leq \Theta_p(m_n-k+1)$. Let
$Y_{n,i}:=( X_{i-k} - \HH_{i-m_n} X_{i-k})( X_i-\tilde X_i)$. Then
  \begin{align*}
    \left\|Y_{n,i}-\Omega_0(Y_{n,i})\right\|_{p/2} \leq
    2\left[\delta_p(i)\Theta_p(m_n-k+1) + \delta_p(i-k)\Theta_p(m_n+1)\right].
  \end{align*}
  Therefore, by (\ref{eq:fact3}), it follows that
  \begin{align*}
    \left\|\sum_{i=k+1}^n \E_0\left[( X_{i-k} - \HH_{i-m_n} X_{i-k})
      ( X_i-\tilde X_i)\right]\right\|_{p/2}
    \leq 4\,\mathcal{C}_{p/2}\Theta_p\Theta_p(m_n-k+1)\cdot\sqrt{n},
  \end{align*}
  and the proof of (\ref{eq:mdep_app}) is complete.
\end{proof}

\subsubsection{Step 2: Throw out small blocks}
In this section, as well as many other places in this article, we
often need to split an integer interval $[s,t] = \{s, s+1, \ldots,
t\} \subset \N$ into consecutive blocks $\mathcal{B}_1, \ldots,
\mathcal{B}_w$ with the size $m$. Since $s-t+1$ may not be a
multiple of $m$, we make the convention that unless the size of
the last block is specified clearly, it has the size $m \leq
|\mathcal{B}_w| <2m$, and all the other ones have the same size
$m$.

\begin{proof}[Proof of Lemma~\ref{thm:small_blk}]
  It suffices to show that for any $\lambda>0$,
  \begin{align*}
    %\label{eq:43}
    \lim_{n\to\infty}\sum_{k=1}^{s_n}
     P\left(\left|\sum_{j=1}^{w_n} V_{k,j}\right|
      \geq \lambda \sqrt{\frac{n}{\log s_n}}\right) =0.
  \end{align*}
  Observe that $V_{k,j}, 1\leq j\leq w_n$, are independent. By
  (\ref{eq:fact5}), $\|V_{k,j}\| \leq
  2|K_j|^{1/2} \kappa_4 \Theta_4$. By Corollary~1.6 of
  \cite{nagaev:1979}, for any $M>1$, there exists a constant $C_{M}>1$
  such that
  \begin{equation}
    \label{eq:48}
      \begin{aligned}
        P\left(\left|\sum_{j=1}^{w_n} V_{k,j}\right|
             \geq \lambda\sqrt{\frac{n}{\log s_n}}\right)
        & \leq \sum_{j=1}^{w_n} P\left(|V_{k,j}|
             \geq C_{M}^{-1}\lambda\sqrt{{n}/{\log s_n}}\right)
        +\left(\frac{4e^2\kappa_4^2\Theta_4^2\sum_{j=1}^{w_n}|K_j|}
          {C_{M}^{-1}\lambda^2 n/\log s_n}\right)^{C_{M}/2} \\
        & \leq \sum_{j=1}^{w_n} P\left(|V_{k,j}|
             \geq C_{M}^{-1}\lambda\sqrt{{n}/{\log n}}\right)
        + C_{M} \left(n^{\beta-\gamma}\log n\right)^{C_{M}/2} \\
        & \leq \sum_{j=1}^{w_n} P\left(|V_{k,j}|
           \geq C_{M}^{-1} \sqrt{{n}/{\log n}}\right) +
           n^{-M}.
  \end{aligned}
  \end{equation}
  where we resolve the constant $\lambda$ into the constant $C_M$ in
  the last inequality. It remains to show that
  \begin{align}
    \label{eq:44}
    \lim_{n\to\infty}\sum_{k=1}^{s_n} \sum_{j=1}^{w_n}
    P\left(|V_{k,j}| \geq q_1\delta \phi_n\right)=0,
    \mbox{ where } \phi_n = \sqrt{\frac{n}{\log n}},
  \end{align}
  holds
  for any $\delta>0$, where $q_1$ is the smallest integer such that
  $\beta^{q_1}<\min\{(p-4)/p, \, (p-2-2\eta)/(p-2)\}$. This choice of
  $q_1$ will be explained later. We adopt the technique of
  successive $m$-dependence approximations from \cite{liu:2010} to
  prove (\ref{eq:44}).

  For $q\geq 1$, set $m_{n,q}=\floor{n^{\beta^q}}$. Define $X_{i,q} =
  \HH_{i-m_{n,q}} X_i$, $\gamma_{k,q} = \E (X_{0,q} X_{k,q})$, and
\begin{eqnarray*}
V_{k,j,q} = \sum_{i \in K_j, i>k} (X_{i-k,q}X_{i,q}-\gamma_{k,q}).
\end{eqnarray*}
In particular, $m_{n,1}$ is same as $m_n$ defined in Step 2, and
$V_{k,j,1}=V_{k,j}$. Without loss of generality assume $s_n \leq
\floor{n^\eta}$. Let $q_0$ be such that $\beta^{q_0+1} \leq \eta <
\beta^{q_0}$. We first consider the difference between $V_{k,j,q}$
and $V_{k,j,q+1}$ for $1\leq q <q_0$. Split the block $K_j$ into
consecutive small blocks $\mathcal{B}_1, \ldots,
\mathcal{B}_{w_{n,q}}$ with size $2m_{n,q}$. Define
  \begin{equation}
    \label{eq:49}
    \begin{aligned}
      V_{k,j,q,t}^{(0)} = \sum_{i \in \mathcal{B}_t}
        (X_{i-k,q}X_{i,q}-\gamma_{k,q}) \quad \hbox{and} \quad
      V_{k,j,q,t}^{(1)} = \sum_{i \in \mathcal{B}_t}
        (X_{i-k,q+1}X_{i,q+1}-\gamma_{k,q+1}).
    \end{aligned}
  \end{equation}
  Observe that $V_{k,j,q,t_1}^{(0)}$ and $V_{k,j,q,t_2}^{(0)}$ are
  independent if $|t_1-t_2|>1$. Similar as (\ref{eq:48}), for any
  $M>1$, there exists a constant $C_M>1$ such that, for
  sufficiently large $n$,
  \begin{equation}
    \label{eq:53}
  \begin{aligned}
    P\left(\left|V_{k,j,q}-V_{k,j,q+1}\right| \geq \delta\phi_n\right)
    & =P\left[\left|\sum_{t=1}^{w_{n,q}} \left(V^{(0)}_{k,j,q,t}-V^{(1)}_{k,j,q,t}\right)\right|
      \geq \delta \phi_n\right] \\
    & \leq \sum_{t=1}^{w_{n,q}} P\left(\left|V^{(0)}_{k,j,q,t}-V^{(1)}_{k,j,q,t}\right|
      \geq C_M^{-1} \phi_n \right) + n^{-M}.
  \end{aligned}
  \end{equation}
  Similarly as (\ref{eq:mdep_app}), we have
  $\left\|V^{(0)}_{k,j,q,t}-V^{(1)}_{k,j,q,t}\right\|_{p/2} \leq C_p
  |\mathcal{B}_t|^{1/2} m_{n,q+1}^{-\alpha}$. It follows that
  \begin{align*}
    \sum_{k=1}^{s_n} \sum_{j=1}^{w_n}
    P\left(\left|V_{k,j,q}-V_{k,j,q+1}\right|
       \geq \delta \phi_n\right)
    & \leq C_{p,M} n^{\eta} n^{1-\gamma} \left(n^{-M} +
      \frac{n^\gamma m_{n,q}^{p/4}
        m_{n,q+1}^{-\alpha p/2}}{m_{n,q}
        (n/\log n)^{p/4}} \right) \\
    & \leq C_{p,M} \left(n^{\eta+1-\gamma-M} + n^{\eta} n^{1-p/4}
      m_{n,q}^{p/4-1-\alpha\beta p/2}\right).
  \end{align*}
  Under the condition (\ref{eq:decay_rate1}), there exists a
  $0<\beta<1$, such that
  \begin{align*}
    %\label{eq:59}
    \sum_{k=1}^{s_n} \sum_{j=1}^{w_n}
    P\left(\left|V_{k,j,q}-V_{k,j,q+1}\right|
    \geq \delta \phi_n \right)
    \leq C_{p,M} \left(n^{\eta+1-\gamma-M} + n^{\eta+ 1-p/4
      +\beta^q (p/4-1-\alpha\beta p/2)} \right) \to 0.
  \end{align*}

  Recall that $q_1$ is the smallest integer such that
  $\beta^{q_1}<\min\{(p-4)/p, (p-2-2\eta)/(p-2)\}$. We now consider
  the difference between $V_{k,j,q}$ and $V_{k,j,q+1}$ for $q_0 \leq q
  <q_1$. The problem is more complicated than the preceding
  case $1\leq q<q_0$, since now it is possible that $m_{n,q}<k$ for
  some $1\leq k \leq s_n$. We consider three cases.

  {\em Case 1: $k \geq 2m_{n,q}$.}
  Partition the block $K_j$ into consecutive smaller blocks
  $\mathcal{B}_1, \ldots, \mathcal{B}_{w_{n,q}}$ with same size
  $m_{n,q}$. Define $V_{k,j,q,t}^{(0)}$ and $V_{k,j,q,t}^{(1)}$ as in
  (\ref{eq:49}). Observe that
  $\left(V_{k,j,q,t}^{(0)}-V_{k,j,q,t}^{(1)}\right)_{t \hbox{
      \scriptsize{is odd}}}$ is a martingale difference sequence with
  respective to the filtration $\left(\xi_t:=\langle
    \epsilon_{l}:\,l\leq \max\left\{\mathcal{B}_t\right\}
    \rangle\right)_{t \hbox{ \scriptsize{is odd}}}$, and so is the
  sequence and filtration labelled by even $t$. Set $\xi_0=\langle
  \epsilon_l:\,l<\min\{\mathcal{B}_1\} \rangle$ and $\xi_{-1} =
  \langle \epsilon_l:\,l<\min\{\mathcal{B}_1\}-m_{n,q} \rangle$. For
  each $1 \leq t \leq w_{n,q}$, define
  \begin{align*}
    \mathcal{V}^{(l)}_{t}
     =\E\left[\left(V^{(l)}_{k,j,q,t}\right)^2 | \xi_{t-2}\right]
    = \sum_{i_1,i_2 \in \mathcal{B}_t} X_{i_1-k,q+l}X_{i_2-k,q+l}
     \gamma_{i_1-i_2,q+l}
  \end{align*}
  for $l=0,1$. By Lemma~1 of \cite{haeusler:1984}, for any $M>1$,
  there exists a constant $C_M>1$ such that
  \begin{equation}
    \label{eq:51}
  \begin{aligned}
    P & \left(\left|V_{k,j,q}-V_{k,j,q+1}\right|
     \geq \delta \phi_n\right)
    \leq \sum_{t=1}^{w_{n,q}}
     P\left(\left|V^{(0)}_{k,j,q,t}-V^{(1)}_{k,j,q,t}\right|
      \geq \sqrt{\frac{n}{(\log n)^3}}\right) + n^{-M} \\
    & \quad + \sum_{l=0,1} 2
     \left\{P\left[\sum_{t \hbox{ {\scriptsize is odd}}} \mathcal{V}^{(l)}_t
      \geq \frac{C_M^{-1}n}{(\log n)^2}\right]
    + P\left[\sum_{t \hbox{ {\scriptsize is even}}} \mathcal{V}^{(l)}_t
      \geq \frac{C_M^{-1}n}{(\log n)^2}\right] \right\}.
  \end{aligned}
  \end{equation}
  By (\ref{eq:fact4}), $\sum_{k \in \Z} |\gamma_{k,q+l}|^2 \leq
  \Theta_2^2$, and hence by (\ref{eq:fact5.5}),
  $\|\mathcal{V}^{(l)}_{t}\|_{p/2} \leq C_p m_{n,q}^{1/2}$. Observe
  that $\mathcal{V}^{(0)}_{t_1}$ and $\mathcal{V}^{(0)}_{t_1}$ are
  independent if $|t_1-t_2|>1$, so similarly as (\ref{eq:48}), we have
  \begin{align*}
    P\left[\sum_{t \hbox{ {\scriptsize is odd}}} \mathcal{V}^{(l)}_t
      \geq \frac{C_M^{-1}n}{(\log n)^2}\right]
    & \leq n^{-M} + \sum_{t\hbox{ {\scriptsize is odd}}} P \left[\mathcal{V}^{(l)}_t
      \geq \frac{C_M^{-2}n}{(\log n)^2}\right] \\
    & \leq n^{-M} + C_{p,M} \cdot w_{n,q} \cdot n^{-p/2} (\log n)^p \cdot m_{n,q}^{p/4}.
  \end{align*}
  The same inequality holds for the sum over even $t$. For the first
  term in (\ref{eq:51}), we claim that
  \begin{align}
    \label{eq:50}
    \left\|V^{(0)}_{k,j,q,t}-V^{(1)}_{k,j,q,t}\right\|_{p}
     \leq C_p \cdot m_{n,q}^{1/2} \cdot m_{n,q+1}^{-\alpha},
  \end{align}
  which together with the preceding two inequalities implies that
  \begin{align*}
    P & \left(\left|V_{k,j,q}-V_{k,j,q+1}\right|
      \geq \delta \phi_n\right)
    \leq C_{p,M} \,w_{n,q} \cdot n^{-p/2} (\log n)^{3p/2}
    \left( m_{n,q}^{p/2} \cdot m_{n,q+1}^{-\alpha p}
     + m_{n,q}^{p/4}\right) + n^{-M}.
  \end{align*}
  It follows that under condition (\ref{eq:decay_rate1}), there
  exists a $0<\beta<1$ such that
  \begin{equation}
    \label{eq:52}
    \begin{aligned}
    & \sum_{k=2m_{n,q}}^{s_n} \sum_{j=1}^{w_n}
    P \left(\left|V_{k,j,q}-V_{k,j,q+1}\right|
       \geq \delta \phi_n\right) \\
    & \quad \leq n^{1+\eta-M}
     + C_{p,M}\cdot n^{1+\eta-p/2} (\log n)^{3p/2}
    \left[ n^{\beta^q (p/2-1-\alpha\beta p)}
     + n^{\beta^q (p/4-1)}\right]
    =o(1).
    \end{aligned}
  \end{equation}

  {\em Case 2: $k \leq m_{n,q+1}/2$.}
  Partition the block $K_j$ into consecutive smaller blocks
  $\mathcal{B}_1, \ldots, \mathcal{B}_{w_{n,q}}$ with size
  $3m_{n,q}$. Define $V_{k,j,q,t}^{(0)}$ and $V_{k,j,q,t}^{(1)}$ as in
  (\ref{eq:49}). Similarly as (\ref{eq:mdep_app}), we have
  \begin{align*}
    \left\|V^{(0)}_{k,j,q,t}-V^{(1)}_{k,j,q,t}\right\|_{p/2}
     \leq C_p \cdot m_{n,q}^{1/2} \cdot m_{n,q+1}^{-\alpha}.
  \end{align*}
  Similar as (\ref{eq:53}), for any
  $M > 1$, there exist a constant $C_M > 1$ such that
  \begin{equation*}
  \begin{aligned}
    P\left(\left|V_{k,j,q}-V_{k,j,q+1}\right| \geq \delta \phi_n\right)
    & \leq \sum_{t=1}^{w_{n,q}}
     P\left(\left|V^{(0)}_{k,j,q,t}-V^{(1)}_{k,j,q,t}\right|
      \geq C_M^{-1} \phi_n \right) + n^{-M} \\
    & \leq n^{-M} + C_{p,M} \cdot w_{n,q} \cdot n^{-p/4} (\log n)^{p/4}
    \cdot m_{n,q}^{p/4}\cdot m_{n,q+1}^{-\alpha\beta p/2}.
  \end{aligned}
  \end{equation*}
  It follows that that under condition (\ref{eq:decay_rate1}), there
  exists a $0<\beta<1$ such that
  \begin{equation}
    \label{eq:56}
    \begin{aligned}
    & \sum_{k=1}^{m_{n,q+1}/2} \sum_{j=1}^{w_n}
    P \left(\left|V_{k,j,q}-V_{k,j,q+1}\right|
     \geq \delta \phi_n\right) \\
    & \quad \leq n^{1+\eta-M} + C_{p,M}\cdot n^{1-p/4} (\log n)^{p/4}
    \cdot \left(n^{\beta^q}\right)^{p/4-\alpha\beta p/2}
    =o(1).
    \end{aligned}
  \end{equation}

  {\em Case 3: $m_{n,q+1}/2 < k < 2m_{n,q}$.} We use the same argument
  as in Case 2. But this time we claim that
  \begin{align}
    \label{eq:57}
    \left\|V^{(0)}_{k,j,q,t}-V^{(1)}_{k,j,q,t}\right\|_{p/2}
     \leq C_p \left[ m_{n,q}^{1/2} \cdot m_{n,q+1}^{-\alpha}
    + m_{n,q}\zeta_p(k)\right],
  \end{align}
  where $\zeta_p(k)$ is defined in (\ref{eq:fact4}). Since
  $\sum_{k=m}^\infty [\zeta_p(k)]^{p/2} \leq
  \left[\sum_{k=m}^\infty\zeta_p(k)\right]^{p/2} = O(m^{-\alpha
    p/2})$, under condition (\ref{eq:decay_rate}), there exist
  constants $C_{p,M}>1$ and $0<\beta<1$ such that for $M$ large enough
  \begin{equation}
    \label{eq:58}
    \begin{aligned}
    & \sum_{k>m_{n,q+1}/2}^{2m_{n,q}-1} \sum_{j=1}^{w_n}
    P \left(\left|V_{k,j,q}-V_{k,j,q+1}\right| \geq \delta \phi_n\right)
    \leq C_{p,M}\cdot n^{1-p/4} (\log n)^{p/4} m_{n,q}^{p/4-\alpha\beta p/2} \\
    & \qquad + n^{1+\eta-M}
    + C_{p,M}\cdot n^{1-p/4} (\log n)^{p/4}
     \cdot m_{n,q}^{p/2-1} \sum_{k>m_{n,q+1}/2}^{2m_{n,q}-1}
      [\zeta_p(k)]^{p/2} \\
    & \quad \leq n^{1+\eta-M} + C_{p,M}\cdot n^{1-p/4} (\log n)^{p/4}
    \cdot m_{n,q}^{p/2-1-\alpha\beta p/2}
    =o(1).
    \end{aligned}
  \end{equation}
  Alternatively, if we use the bound from (\ref{eq:fact9}),
  $\left\|V^{(0)}_{k,j,q,t}-V^{(1)}_{k,j,q,t}\right\|_{p/2} \leq C_p
  m_{n,q}^{1/2} \cdot m_{n,q+1}^{-\alpha'}$, it is still true that
  under condition (\ref{eq:decay_rate}), there exist
  constants $C_{p,M}>1$ and $0<\beta<1$ such that for $M$ large enough
  \begin{equation}
    \label{eq:10}
    \begin{aligned}
      \sum_{k>m_{n,q+1}/2}^{2m_{n,q}-1} & \sum_{j=1}^{w_n}
      P \left(\left|V_{k,j,q}-V_{k,j,q+1}\right|
       \geq \delta \phi_n\right) \\
      & \leq n^{1+\eta-M} + C_{p,M}\cdot n^{1-p/4} (\log n)^{p/4}
       \cdot m_{n,q}^{p/2-1-\alpha'\beta p/2}
      =o(1).
    \end{aligned}
\end{equation}
  Combine (\ref{eq:52}), (\ref{eq:56}), (\ref{eq:58}) and (\ref{eq:10}),
  we have shown that
  \begin{align}
    \label{eq:61}
    \lim_{n\to\infty}\sum_{k=1}^{s_n} \sum_{j=1}^{w_n}
    P\left(\left|V_{k,j,q}-V_{k,j,q+1}\right| \geq \delta \phi_n\right) = 0.
  \end{align}
  for $1\leq q <q_1$. Therefore, to prove (\ref{eq:44}), it suffices to show
  \begin{align}
    \label{eq:60}
    \lim_{n\to\infty}\sum_{k=1}^{s_n} \sum_{j=1}^{w_n}
    P\left(|V_{k,j,q_1}| \geq \delta \phi_n \right)=0
  \end{align}
  By considering two cases (i) $2m_{n,q_1} \leq k \leq s_n$ and (ii)
  $1\leq k <2m_{n,q_1}$ under the condition
  $\beta^{q_1}<\min\{(p-4)/p, (p-2-2\eta)/(p-2)\}$, and using similar
  arguments as those in proving (\ref{eq:61}), we can obtain
  (\ref{eq:60}). The proof of Lemma~\ref{thm:small_blk} is complete.

  We now turn to the proof of the two claims (\ref{eq:50}) and
  (\ref{eq:57}). For (\ref{eq:57}), we have
  \begin{align*}
    \left\|V^{(0)}_{k,j,q,t}-V^{(1)}_{k,j,q,t}\right\|_{p/2}
    \leq & \left\|\sum_{i\in\mathcal{B}_t}
       ( X_{i-k,q} - X_{i-k,q+1})X_{i,q+1}\right\|_{p/2}
    + \left\|\sum_{i\in\mathcal{B}_t}
     \E_0\left[X_{i-k,q+1}( X_{i,q}-X_{i,q+1})\right]\right\|_{p/2} \\
    & + \left\|\sum_{i\in\mathcal{B}_t}
     \E_0\left[( X_{i-k,q}-X_{i-k,q+1})(X_{i,q}-X_{i,q+1})\right]\right\|_{p/2}
    =: I + \II + \III.
  \end{align*}
  Similarly as in the proof of (\ref{eq:mdep_app}), we have
  \begin{align*}
    I \leq \CC_{p/2}\Theta_p\Theta_p(m_{n,q+1}+1)\cdot\sqrt{3m_{n,q}}
     \quad\hbox{and}\quad
    \III \leq 4\,\CC_{p/2}\Theta_p\Theta_p(m_{n,q+1}+1)\cdot\sqrt{3m_{n,q}}.
  \end{align*}
  For the second term $\II$, write
  \begin{align*}
    \E_0\left[X_{i-k,q+1}( X_{i,q}-X_{i,q+1})\right]
    = \sum_{l_1=0}^{m_{n,q+1}} \sum_{l_2=m_{n,q+1}+1}^{m_{n,q}}
     \E_0\left[(\PP_{i-k-l_1}X_{i-k})(\PP_{i-l_2}X_{i})\right].
  \end{align*}
  For a pair $(l_1,l_2)$ such that $i-k-l_1 \neq i-l_2$, by the
  inequality (\ref{eq:mar_zyg}), we have
  \begin{align*}
    \left\|\sum_{i \in \mathcal{B}_t}
     (\PP_{i-k-l_1}X_{i-k})(\PP_{i-l_2}X_{i})\right\|_{p/2}
    \leq \CC_{p/2} \delta_p(l_1)\delta_p(l_2)\cdot \sqrt{3m_{n,q}}.
  \end{align*}
  For the pairs $(l_1,l_2)$ such that $i-k-l_1 = i-l_2$, by the
  triangle inequality
  \begin{align*}
    \left\|\sum_{i\in\mathcal{B}_t}\sum_{l=0}^{m_{n,q+1}}
      \E_0\left[(\PP_{i-k-l}X_{i-k})(\PP_{i-k-l}X_{i})\right]\right\|_{p/2}
    \leq 3m_{n,q}\cdot2\sum_{l=0}^{m_{n,q+1}} \delta_p(l)
     \delta_p(k+l) \leq 6m_{n,q} \zeta_p(k).
  \end{align*}
  Putting these pieces together, the proof of (\ref{eq:57}) is
  complete.  The key observation in proving (\ref{eq:50}) is that
  since $k \geq 2m_{n,q}$, $X_{i-k,q}$ and $X_{i,q}$ are independent,
  hence the product $X_{i-k,q}X_{i,q}$ has finite $p$-th moment.
  The rest of the proof is similar to that of
  (\ref{eq:57}). Details are omitted.
\end{proof}

\begin{remark}
  \label{rk:decay_rate}
  Condition (\ref{eq:decay_rate}) is only used to deal with Case~3,
  while (\ref{eq:decay_rate1}) suffices for the rest of the proof. In
  fact, for linear processes, one can show that the term
  $m_{n,q}\zeta_p(k)$ in (\ref{eq:57}) can be removed, so we have
  (\ref{eq:58}) under condition (\ref{eq:decay_rate1}) and do not
  need (\ref{eq:10}). So (\ref{eq:decay_rate1}) suffices for
  Theorem~\ref{thm:single}. Furthermore, for nonlinear
  processes with $\delta_p(k)=O\left[k^{-(1/2+\alpha)}\right]$,
  the term $m_{n,q}\zeta_p(k)$ can also be removed from
  (\ref{eq:57}). Details are omitted.
\end{remark}

\subsubsection{Step 3: Truncate sums over large blocks}

\begin{proof}[Proof of Lemma~\ref{thm:large_blk}]
  We need to show for any $\lambda>0$
  \begin{align*}
    % \label{eq:43}
    \lim_{n\to\infty}\sum_{k=1}^{s_n}
     P\left(\left|\sum_{j=1}^{w_n} (U_{k,j}-\bar U_{k,j})\right|
      \geq \lambda \sqrt{\frac{n}{\log s_n}}\right) =0.
  \end{align*}
  Using (\ref{eq:fact5}), elementary calculation gives
  \begin{align}
    \label{eq:31}
    \left\|{\tilde U}_{k,j}-\bar{U}_{k,j}\right\|^2
    \leq \frac{\E |{\tilde U}_{k,j}|^{p/2}}{(\sqrt{n}/\log s_n)^{p/2-2}}
    \leq \frac{(2\CC_{p/2}\kappa_p\Theta_p)^{p/2}|H_j|^{p/4}
     (\log s_n)^{3(p-4)/2}}{n^{(p-4)/4}}.
  \end{align}
  Similarly as (\ref{eq:48}), for any $M>1$, there exists a constant
  $C_{M}>1$ such that
  \begin{equation*}
    %\label{eq:48}
      \begin{aligned}
        P \left(\left|\sum_{j=1}^{w_n} (U_{k,j}-\bar U_{k,j})\right|
          \geq \lambda \sqrt{\frac{n}{\log s_n}}\right)
        \leq & \sum_{j=1}^{w_n} P\left(|U_{k,j}-\bar U_{k,j}|
          \geq C_{M}^{-1} \lambda \sqrt{\frac{n}{\log s_n}}\right) \\
        & +\left(\frac{C_p\sum_{j=1}^{w_n}|H_j|^{p/4}(\log n)^{3p/2}}
          {C_{M}^{-1} \lambda^2 n^{p/4}}\right)^{C_{M}/2} \\
        \leq & \sum_{j=1}^{w_n} P\left(|U_{k,j}-\bar U_{k,j}|
          \geq C_{M}^{-1}\sqrt{\frac{n}{\log s_n}}\right) + n^{-M}.
  \end{aligned}
  \end{equation*}
  Therefore, it suffices to show that for any $\delta>0$,
  \begin{align*}
    \lim_{n\to\infty}\sum_{k=1}^{s_n} \sum_{j=1}^{w_n}
    P\left(|U_{k,j}-\bar U_{k,j}| \geq \delta\sqrt{\frac{n}{\log n}}\right)=0.
  \end{align*}
Since we can use the same arguments as those for (\ref{eq:44}),
Lemma~\ref{thm:large_blk} follows.
\end{proof}

\subsubsection{Step 4: Compare covariance structures}

\begin{proof}[Proof of Lemma~\ref{thm:small_lag}]
  Since $|\bar U_{k,j}| \leq 2\sqrt{n}/(\log s_n)$ and $\E \bar
  U_{k,j}^2 \leq \E U_{k,j}^2 \leq 4(\kappa_4\Theta_4)^2|H_j|$, by
  Bernstein's inequality \citep[cf. Fact~2.3,][]{einmahl:1997},
  we have
  \begin{align*}
    P \left(|\RR_{n,k}| > \sqrt{\sigma_0 n\log s_n}\right)
    \leq \exp\left\{-\frac{(\sigma_0 n\log s_n)/2}
      {4(\kappa_4\Theta_4)^2 n
          + n \sqrt{\sigma_0/(\log s_n)}}\right\}.
  \end{align*}
Therefore, for any $0< \iota < \sigma_0/[8(\kappa_4 \Theta_4)^2]$,
(\ref{eq:small_lag}) holds.
\end{proof}

\begin{proof}[Proof of Lemma~\ref{thm:cov_struc}]
  For $1\leq j\leq w_n$, by (\ref{eq:31}), we have
  \begin{align*}
    \left|\E (\bar U_{k,j}\bar U_{k+h,j})
     - \E ({\tilde U}_{k,j}{\tilde U}_{k+h,j})\right|
    & \leq \|\bar U_{k,j}-{\tilde U}_{k,j}\|\|\bar U_{k+h,j}\|
    + \|{\tilde U}_{k,j}\|\|\bar U_{k+h,j} - {\tilde U}_{k+h,j}\| \cr
    & \leq 4\kappa_4\Theta_4|H_j|^{1/2}
    \frac{(2\CC_{p/2}\kappa_p\Theta_p)^{p/4}|H_j|^{p/8}(\log s_n)^{3(p-4)/4}}{n^{(p-4)/8}}\cr
    & \leq C_{p} |H_j| n^{-(1-\gamma)(p-4)/8}(\log n)^{3(p-4)/4}.
  \end{align*}
  Let $S_{k,j}=\sum_{i\in H_j}(X_{i-k}X_i- \gamma_k)$, by
  (\ref{eq:fact5}) and (\ref{eq:mdep_app}), we have
  \begin{align*}
    \left|\E (S_{k,j}S_{k+h,j})
     - \E (\tilde U_{k,j}\tilde U_{k+h,j})\right|
    & \leq \|S_{k,j}-\tilde U_{k,j}\|\|S_{k+h,j}\|
     + \|\tilde U_{k,j}\|\|S_{k+h,j} - \tilde U_{k+h,j}\| \cr
    & \leq 4\kappa_4\Theta_4|H_j|^{1/2}\cdot
      6\Theta_4\Theta_4(m_n-k+1) |H_j|^{1/2}
    \leq C|H_j|n^{-\alpha\beta}.
  \end{align*}
  Since $\Theta_4(m) = O(m^{-\alpha})$, elementary calculation shows
  that $\Delta_4(m) = O(n^{-\alpha^2/(1+\alpha)})$, which together
  with Lemma~\ref{thm:covconvergence} implies that if $k>t_n$,
  \begin{align*}
    \left|\E (\tilde U_{k,j}\tilde U_{k+h,j})/|H_j| - \sigma_h\right|
    & \leq \Theta_{4}^3\left(16\Delta_4(t_n/3+1) +
      6\Theta_4\sqrt{t_n/l_n} + 4\Psi_4(t_n/3+1)\right)\cr
    & \leq C \left(s_n^{-\alpha^2\iota/(1+\alpha)}+n^{-(1-\iota)\gamma/2}\right).
  \end{align*}
  % Since $|\E \bar U_{k,w_n}\bar U_{k+h,w_n}| \leq \|\tilde
  % U_{k,w_n}\|\cdot\|\tilde U_{k+h,w_n}\| \leq
  % 4\kappa_4^2\Theta_4^2|H_{w_n}| \leq 8\Theta_4^4l_n$, we have
Choose $\ell$ such that $0<\ell<\min\{(1-\eta)(p-4)/8,\,
\alpha\beta, \, \alpha^2\iota/(1+\alpha),\, (1-\iota)\gamma/2,\,
\gamma-\beta\}$. Then
  \begin{align*}
    |\Cov(\RR_{n,k}, \RR_{n,k+h})/n-\sigma_h|
     & \leq C_p \Big(n^{-(1-\eta)(p-4)/8}(\log n)^{(p-4)/4}
      + n^{-\alpha\beta} \cr
     & \qquad  + s_n^{-\alpha^2\iota/(1+\alpha)}
      +n^{-(1-\iota)\gamma/2}\Big)
     + {{2w_n m_n\sigma_0} \over n}
      \leq C_p\,s_n^{-\ell}
  \end{align*}
and the lemma follows.
\end{proof}

\subsubsection{Step 5: Moderate deviations.}

\begin{proof}[Proof of Lemma~\ref{thm:vec_mod_dev}]
  Note that for $\h{x}, \h{y} \in \R^{d}$, $|x+y|_\bullet\leq
  |x|_\bullet+|y|$. Let $\h{Z}\sim\mathcal{N}(0,I_d)$ and
  $\theta_n=(\log s_n)^{-1}$. Since $|\bar U_{k,j}| \leq
  2\sqrt{n}/(\log s_n)^3$, by Fact~2.2 of \cite{einmahl:1997},
  \begin{align*}
    P(|\h{\RR}_n/\sqrt{n}|_\bullet\geq z_n)
      & \leq P(|\Sigma_n^{1/2}\h{Z}|_{\bullet} \geq z_n - \theta_n)
    + P (|\h{\RR}_n/\sqrt{n} - \Sigma_n^{1/2}\h{Z}| \geq \theta_n) \\
    & \leq P(|\Sigma_n^{1/2}\h{Z}|_{\bullet} \geq z_n - \theta_n)
    + C_{p,d}\, \exp\left\{-C_{p,d}^{-1}(\log s_n)^2\right\}.
  \end{align*}
  By Lemma~\ref{thm:eigenbdd}, the smallest eigenvalue of $\Sigma$ is
  bounded from below by some $c_d>0$ uniformly on $1 \leq
  k_1<k_2<\cdots<k_d$. By Lemma~\ref{thm:cov_struc} we have
  $\rho(\Sigma^{1/2}_n-\Sigma^{1/2}) \leq
  c_d^{-1/2}\cdot\rho(\Sigma_n-\Sigma)\leq C_{p,d} \,s_n^{-\ell}$, where
  the first inequality is taken from Problem~7.2.17 of
  \cite{horn:1990}. It follows that by (\ref{eq:normbdd}) and
  elementary calculations that
  \begin{align*}
    P(|\Sigma_n^{1/2}\h{Z}|_{\bullet} \geq z_n - \theta_n)
    & \leq P(|\Sigma^{1/2}\h{Z}|_{\bullet} \geq z_n - 2\theta_n)
    + P\left[\left|\left(\Sigma^{1/2}_n-\Sigma^{1/2}\right)\h{Z}\right|
      \geq \theta_n \right] \\
    & \leq P(|\Sigma^{1/2}\h{Z}|_{\bullet} \geq z_n - 2\theta_n)
     + C_{p,d}\, \exp\left\{s_n^{-\ell}\right\}.
  \end{align*}
  By Lemma~\ref{thm:lemma1}, we have
  \begin{align*}
    P(|\Sigma^{1/2}\h{Z}|_{\bullet} \geq z_n - 2\theta_n)
    \leq \left[1+C_{p,d}(\log s_n)^{-1/2}\right]
     P(|\Sigma^{1/2}\h{Z}|_{\bullet} \geq z_n).
  \end{align*}
  Putting these pieces together and observing that $\h{V}$ and
  $\Sigma^{1/2}\h{Z}$ have the same distribution, we have
  \begin{align*}
    P(|\h{\RR}_n/\sqrt{n}|_\bullet\geq z_n)
    \leq \left[1+C_{p,d}(\log s_n)^{-1/2}\right] P(|\h{V}|_\bullet\geq z_n)
    + C_{p,d}\, \exp\left\{-C_{p,d}^{-1}(\log s_n)^2\right\},
  \end{align*}
which together with a similar lower bound completes the proof of
Lemma~\ref{thm:vec_mod_dev}.
\end{proof}

\subsubsection{Proof of Theorem~\ref{thm:single}}

After these preparation steps, we are now ready to prove
Theorem~\ref{thm:single}.
\begin{proof}[Proof of Theorem~\ref{thm:single}]
  Set $z_n = a_{2s_n}\,x+b_{2s_n}$. It suffices to show
  \begin{align}\label{eq:gumbel}
    \lim_{n\to\infty}P\left(\max_{t_n< k\leq s_n}|\RR_k/\sqrt{n}|
      \leq \sqrt{\sigma_0} z_n\right) = \exp\{-\exp(-x)\}.
  \end{align}
  Without loss of generality assume $\sigma_0=1$.
  Define the events $A_k=\{G_k\geq z_n\}$ and
  $B_k=\{\RR_k/\sqrt{n}\geq z_n\}$. Let
  \begin{align*}
    Q_{n,d} = \sum_{t_n< k_1<\ldots<k_d\leq s_n}
     P(A_{k_1}\cap \cdots \cap A_{k_d})
     \quad\hbox{and}\quad
    \tilde Q_{n,d} = \sum_{t_n< k_1<\ldots<k_d\leq s_n}
    P(B_{k_1}\cap \cdots \cap B_{k_d}).
  \end{align*}
  By the inclusion-exclusion formula, we know for any $q\geq 1$
  \begin{align}
    \label{eq:inex}
    \sum_{d=1}^{2q} (-1)^{d-1}\tilde Q_{n,d}
     \leq P\left(\max_{t_n< k\leq s_n}|\RR_k/\sqrt{n}|
      \geq a_{2s_n}\,x+b_{2s_n}\right)
      \leq \sum_{d=1}^{2q-1} (-1)^{d-1}\tilde Q_{n,d}.
  \end{align}
By Lemma~\ref{thm:vec_mod_dev}, $|\tilde Q_{n,d}-Q_{n,d}| \leq
C_{p,d}(\log s_n)^{-1/2}Q_{n,d} + s_n^{-1}. $ By
Lemma~\ref{thm:normalcomparison} with elementary calculations, we
know $\lim_{n\to \infty}Q_{n,d} = e^{-dx}/d!$, and hence $\lim_{n
\to \infty}\tilde Q_{n,d} = e^{-dx}/d!$. By letting $n$ go to
infinity first and then $d$ go to infinity in (\ref{eq:inex}), we
obtain (\ref{eq:gumbel}), and the proof is complete.
\end{proof}

\subsection{Proof of Theorem~\ref{thm:covorder}}

\begin{proof}[Proof of Theorem~\ref{thm:covorder}]
We start with an $m$-dependence approximation that is similar to
the proof of Theorem~\ref{thm:single}. Set $m_n=\floor{n^\beta}$
for some $0<\beta<1$. Define $\tilde X_i = \HH_{i-m_n} X_i$,
$\tilde \gamma_k = \E (\tilde X_0 \tilde X_k)$, and $\tilde
R_{n,k}=\sum_{i=k+1}^n (\tilde X_{i-k} \tilde X_i - \tilde
\gamma_k )$. Similarly as the proof of Lemma~\ref{thm:small_blk},
we have under the condition (\ref{eq:decay_rate2})
  \begin{align*}
    \max_{1 \leq k <n} |R_{n,k}-\tilde R_{n,k}|
     = o_P\left(\sqrt{{n}/{\log n}}\right).
  \end{align*}
  For $\tilde R_{n,k}$, we consider two cases according to whether $k
  \geq 3m_n$ or not.

  {\em Case 1: $k \geq 3m_n$.}  We first split the interval $[k+1,n]$
  into the following big blocks of size $(k-m_n)$
  \begin{equation*}
    \begin{aligned}
      & H_j=[k+{j-1}(k-m_n)+1,k+j(k-m_n)]
          \quad \hbox{for } 1 \leq j \leq w_n-1 \\
      & H_{w_n}=[k+(w_n-1)(k-m_n)+1,n],
    \end{aligned}
  \end{equation*}
  where $w_n$ is the smallest integer such that $k+w_n(k-m_n) \geq
  n$. For each block $H_j$, we further split it into small blocks of
  size $2m_n$
  \begin{equation*}
    \begin{aligned}
      & K_{j,l}=[k+(j-1)(k-m_n)+(l-1)2m_n+1,k+(j-1)(k-m_n)+2lm_n] \quad \hbox{for } 1 \leq l < v_j \\
      & K_{j,v_j}=[k+(v_j-1)(k-m_n)+(l-1)2m_n+1,k+(j-1)(k-m_n)+|H_j|]
    \end{aligned}
  \end{equation*}
  where $v_j$ is the smallest integer such that
  $2m_nv_j\geq|H_j|$. Now define $U_{k,j,l}=\sum_{i \in K_{j,l}}
  \tilde X_{i-k}\tilde X_i$ and
  \begin{align}
    \label{eq:11}
    \tilde R_{n,k}^{u,1}=\sum_{j \equiv u \!\!\!\!\pmod 3}\sum_{l \hbox{ \scriptsize{odd}}} U_{k,j,l} \quad \hbox{and} \quad
    \tilde R_{n,k}^{u,2}=\sum_{j \equiv u \!\!\!\!\pmod 3}\sum_{l \hbox{ \scriptsize{even}}} U_{k,j,l}
  \end{align}
  for $u=0,1,2$. Observe that each $\tilde R_{n,k}^{u,o}$
  ($u=0,1,2; \;o=1,2$) is a sum of independent random variables. % Since
  % $\tilde X_i$ and $\tilde X_{i-k}$ are independent, $U_{k,j,l}$ have
  % finite $p$-th moment, and similarly as
  By (\ref{eq:fact5}), $\|U_{k,j,l}\| \leq
  2\kappa_4\Theta_4|U_{k,j,l}|^{1/2}$. By Corollary~1.7 of
  \cite{nagaev:1979} where we take $y_i=\sqrt{n}$ in their result, we
  have for any $\lambda>0$
  \begin{equation}
    \label{eq:40}
    \begin{aligned}
   P & \left(|\tilde R_{n,k}| \geq 6 \lambda \sqrt{n\log n}\right)
    \leq \sum_{u=0}^2 \sum_{o=1,2}
    P\left(\left|\tilde R_{n,k}^{u,o}\right| \geq \lambda \sqrt{n\log n}\right) \\
   & \leq \sum_{u=0}^2 \sum_{o=1,2} \sum_{j,l}^{\ast} P \left(|U_{k,j,l}|
     \geq \lambda \sqrt{n\log n}\right)
      + 12 \left(\frac{C_p \,n^{1-\beta} \cdot
       n^{\beta p/4}}{n^{p/4}}\right)^{p\sqrt{\log n}/(p+4)} \\
      & \quad + 12\exp\left\{-\frac{2 \lambda^2}
                  {(p+4)^2\cdot e^{p/2}\cdot\kappa_4^2\cdot\Theta_4^2}
        \cdot\log n\right\} =: I_{n,k} + \II_{n,k} + \III_{n,k},
    \end{aligned}
  \end{equation}
where the range of $j,l$ in the sum $\sum_{j,l}^{\ast}$ is as in
(\ref{eq:11}). Clearly, $\sum_{k=3m_n}^{n-1}
  \II_{n,k}=o(1)$. Similarly as the proof of
  Lemma~\ref{thm:small_lag}, we can show that $\sum_{k=3m_n}^{n-1}
  I_{n,k}=o(1)$. Therefore, if $\epsilon = c_p/6$, then
  $\sum_{k=3m_n}^{n-1} \III_{n,k}=O(n^{-1})$.

  {\em Case 2: $1\leq k < 3m_n$.} This case is easier. By splitting
  the interval $[k+1,n]$ into blocks with size $4m_n$ and using a
  similar argument as (\ref{eq:40}), we have
  \begin{align*}
    \lim_{n\to\infty}\sum_{k=1}^{3m_n-1} P \left(|\tilde R_{n,k}| \geq c_p \sqrt{n\log n}\right) =0.
  \end{align*}
  The proof is complete.

  % Set $m_n=\floor{n^\gamma}$ and $\tilde X_i=\HH_{i-m_n}^i X_i$. For
  % every $3m_n \leq k < n$, define $\tilde R_{n,k}=\sum_{i=k+1}^n\tilde
  % X_{i-k}\tilde X_i$. By (\ref{eq:fact9}), we have
  % \begin{equation}
  %   \label{eq:39}
  %   \begin{aligned}
  %     P \left(\left|R_{n,k} - \tilde R_{n,k}\right| \geq \epsilon \sqrt{n \log n}\right)
  %     \leq \frac{(4\CC_{p/2}\kappa_p)^{p/2}n^{-\alpha\gamma p/2}}{\epsilon^{p/2}(\log n)^{p/4}}
  %     =C_p\epsilon^{-p/2}n^{-\alpha\gamma p/2}(\log n)^{-p/4}.
  %   \end{aligned}
  % \end{equation}

  % For $1 \leq k \leq 3\floor{n^\gamma}$, pick $m_n =
  % 6\floor{n^\gamma}$, by a similar argument as
  % Lemma~\ref{thm:small_lag}, we have
  % \begin{equation}
  %   \label{eq:42}
  %   \begin{aligned}
  %     P \left( |R_{n,k}| \geq  \epsilon \sqrt{n\log n}\right)
  %     \leq C_p\epsilon^{-p/2}\left(n^{-\alpha\gamma p/2} + n^{-(p/4-1)(1-\gamma)}\right)(\log n)^{-p/4}
  %     + 4\exp\left\{{\epsilon^2\log n}/{C_p}\right\}.
  %   \end{aligned}
  % \end{equation}
  % If $\alpha \geq 2/(p-4)$, then there exists a $0<\gamma<1$ such that
  % \begin{align*}
  %   \alpha\gamma p/2 \geq 1, \quad (p/2-1)(1-\gamma) \geq 1 \quad \hbox{and} \quad (p/4-1)(1-\gamma) \geq \gamma.
  % \end{align*}
  % Therefore, the proof is complete in view of (\ref{eq:41}) and (\ref{eq:42}).
\end{proof}

%\subsection{Stationary processes with iid repetitions}

\subsection{Box-Pierce tests}
Similarly as the proof of Theorem~\ref{thm:single}, we use
$m$-dependence approximations and blocking arguments to prove
Theorem~\ref{thm:ljung}. We first outline the intermediate steps and
give the main proof in Section~\ref{sec:proof_ljung}, and then provide
proofs of the intermediate lemmas in Section~\ref{sec:proof_blk} and
Section~\ref{sec:proof_cleared}. We prove
Theorem~\ref{thm:ljung_power} in Section~\ref{sec:proof_power}, and
prove Corollary~\ref{thm:ljung_corr} and \ref{thm:ljung_corr_power} in
Section~\ref{sec:proof_corollary}.

\subsubsection{Proof of Theorem~\ref{thm:ljung}}
\label{sec:proof_ljung}

\paragraph{\underline{\textnormal{\em Step 1: $m$-dependence approximation.}}}
Recall that $R_{n,k}=\sum_{i=k+1}^n (X_{i-k}X_i-\gamma_k)$.  Without
loss of generality, assume $s_n \leq \floor{n^\beta}$.  Set
$m_n=2\floor{n^\beta}$. Let $\tilde X_i=\HH_{i-m_n}^i X_i$ and $\tilde
R_{n,k}=\sum_{i=k+1}^n (\tilde X_{i-k}\tilde X_i - \tilde\gamma_k
)$. By (\ref{eq:fact5}) and (\ref{eq:mdep_app}), we know if
$\Theta_4(m)=o(m^{-\alpha})$ for some $\alpha > 0$, then for all
$1\leq k\leq s_n$
\begin{align*}
  \E  |R_{n,k}^2 - \tilde R_{n,k}^2 | \leq \|R_{n,k}+\tilde R_{n,k}\| \cdot \|R_{n,k}-\tilde R_{n,k}\|
  \leq C\, \Theta_4^3 \cdot n \cdot \Theta_4\left(m_n/2\right) = o\left(n^{1-\alpha\beta}\right).
\end{align*}
The condition $\sum_{k=0}^\infty k^6\delta_8(k)<\infty$ implies that
$\Theta_4(m)=O(m^{-6})$. Therefore, under the conditions of
Theorem~\ref{thm:ljung}, we have
\begin{align*}
  %\label{eq:ljung_mdep}
  \frac{1}{n\sqrt{s_n}} \sum_{k=1}^{s_n}
  \E_0\left(R_{n,k}^2 - \tilde R_{n,k}^2\right)
  =o_P(1).
\end{align*}

\paragraph{\underline{\textnormal{\em Step 2: Throw out small blocks.}}}
Let $l_n=\floor{n^\eta}$, where $\eta \in (\beta,1)$. Split the
interval $[1,n]$ into alternating small and large blocks similarly as
(\ref{eq:blk}):
\begin{equation*}
  \begin{aligned}
  & K_0=[1,s_n] \cr
  & H_j = [s_n+(j-1)(2m_n+l_n)+1,s_n+(j-1)(2m_n+l_n)+l_n]
  \quad 1 \leq j \leq w_n\cr
  & K_{j}=[s_n+(j-1)(2m_n+l_n)+l_n+1,s_n+j(2m_n+l_n)]; \quad 1\leq j \leq w_n-1;
  \quad \hbox{and} \cr
  & K_{w_n}=[s_n+(w_n-1)(2m_n+l_n)+l_n+1,n],
\end{aligned}
\end{equation*}
where $w_n$ is the largest integer such that
$s_n+(w_n-1)(2m_n+l_n)+l_n\leq n$.  Define $U_{k,0}=0$,
$V_{k,0}=\sum_{i\in K_0,i>k}(\tilde X_{i-k}\tilde
X_{i}-\tilde\gamma_k)$, and $U_{k,j}=\sum_{i\in H_j}(\tilde
X_{i-k}\tilde X_{i}-\tilde\gamma_k)$, $V_{k,j}=\sum_{i\in K_j}(\tilde
X_{i-k}\tilde X_{i}-\tilde\gamma_k)$ for $1\leq j \leq w_n$. Set
$\RR_{n,k}=\sum_{j=1}^{w_n}U_{k,j}$. Observe that by construction,
$U_{k,j},1\leq j \leq w_n$ are iid random variables. In the following
lemma we show that it suffices to consider $\RR_{n,k}$.
\begin{lemma}
  \label{thm:blk}
Assume $X_i \in \LLL^8$, $\E X_i=0$, and $\sum_{k=0}^\infty
k^6\delta_8(k) < \infty$, then
  \begin{align*}
    \frac{1}{n\sqrt{s_n}}\sum_{k=1}^{s_n}
    \E_0\left(\tilde R_{n,k}^2 - \RR_{n,k}^2\right) = o_P(1).
  \end{align*}
\end{lemma}

\paragraph{\underline{\textnormal{\em Step 3: Central limit theorem concerning $\RR_{n,k}$'s.}}}
\begin{lemma}
  \label{thm:ljung_cleared}
  Assume $X_i \in \LLL^8$, $\E X_i=0$, and $\sum_{k=0}^\infty k^6\delta_8(k) < \infty$, then
  \begin{align*}
    \frac{1}{n\sqrt{s_n}}\sum_{k=1}^{s_n} \left(\RR_{n,k}^2 - \E\RR_{n,k}^2\right)
    \Rightarrow \mathcal{N}\left(0,2\sum_{k\in\Z}\sigma_k^2\right).
  \end{align*}
\end{lemma}

We are now ready to prove Theorem~\ref{thm:ljung}.
\begin{proof}[Proof of Theorem~\ref{thm:ljung}]
  By Lemma~\ref{thm:blk} and Lemma~\ref{thm:ljung_cleared}, we know
  \begin{align*}
    \frac{1}{n\sqrt{s_n}}\sum_{k=1}^{s_n}
     \left(R_{n,k}^2-\E R_{n,k}^2\right)
    \Rightarrow
    \mathcal{N}\left(0, 2\sum_{k\in\Z}\sigma_k^2\right).
  \end{align*}
  It remains to show that
  \begin{align}
    \label{eq:ljung_bias}
    \lim_{n\to\infty} \frac{1}{n\sqrt{s_n}}\sum_{k=1}^{s_n}
     \left[\E R_{n,k}^2-(n-k)\sigma_0\right] = 0.
  \end{align}
  We need Lemma~\ref{thm:covconvergence} with a slight
  modification. Observe that in equation (\ref{eq:21}), we now have
  $\sum_{j=1}^{m_n} \Theta_2(j)^2<\infty$, and hence
  \begin{align*}
    \left|\E R_{n,k}^2 - (n-k)\sigma_0\right|
     \leq C \left[(n-k)\Delta_4(\floor{k/3}+1) + \sqrt{n-k}\right]
  \end{align*}
  With the condition $\Theta_8(m)=o(m^{-6})$, elementary calculations
  show that $\Delta_4(m)=o(m^{-5})$. Then (\ref{eq:ljung_bias})
  follows, and the proof is complete.
\end{proof}

\subsubsection{Step 2: Throw out small blocks.}
\label{sec:proof_blk}

Let $\mathcal{A}_2$ be the collection of all double arrays
$A=(a_{ij})_{i,j\geq 1}$ such that
$$\|A\|_{\infty}:=\max\left\{\sup_{i\geq
    1}\sum_{j=1}^{\infty}|a_{ij}|,\, \sup_{j\geq
    1}\sum_{i=1}^{\infty}|a_{ij}|\right\}<\infty.$$ For
$A,B\in\mathcal{A}_2$, define $AB=(\sum_{k=1}^{\infty} a_{ik}
b_{kj})$. It is easily seen that $AB \in \mathcal{A}_2$ and
$\|AB\|_{\infty} \le \|A\|_{\infty} \|B\|_{\infty}$. Furthermore,
this fact implies the following proposition, which will be useful
in computing sums of products of cumulants. For $d \geq 0$, let
$\mathcal{A}_d$ be the collection of all $d$-dimensional array
$A=A(i_1,i_2,\ldots,i_d)$ such that
\begin{equation*}
  \|A\|_\infty:=\max_{1\leq j\leq d}
  \left\{\sup_{i_j\geq 1}\sum_{\{i_k:\,k\neq j\}}
   |A(i_1,i_2,\ldots,i_d)|\right\} < \infty.
\end{equation*}
Note that $\mathcal{A}_0=\R$, and $\|A\|_\infty=|A|$ if $A \in
\mathcal{A}_0$.

\begin{proposition}
\label{thm:array} For $k\ge 0$, $l\ge 0$ and $d \ge 1$, if $A \in
\mathcal{A}_{k+d}$ and $B \in \mathcal{A}_{l+d}$, define an array
$C$ by
  \begin{align*}
    C(i_1,\ldots,i_k,i_{k+1},\ldots,i_{k+l})
     =\sum_{j_1,\ldots,j_d\geq 1}
     A(i_1,\ldots,i_k,j_1,\ldots,j_d)
     B(j_1,\ldots,j_d,i_{k+1},\ldots,i_{k+l})
  \end{align*}
  then $C \in \mathcal{A}_{k+l}$, and $\|C\|_{\infty}\leq
  \|A\|_{\infty} \|B\|_{\infty}$.
\end{proposition}

In Lemma~\ref{thm:mag} we present an upper bound for
$\Cov(R_{n,k},R_{n,h})$. We formulate the lemma in a more general way
for later uses in the proofs of Lemma~\ref{thm:blk} and
Lemma~\ref{thm:ljung_cleared}.  For a $k$-dimensional random vector
$(Y_1,\ldots,Y_k)$ such that $\|Y_i\|_k<\infty$ for $1\leq i \leq k$,
denote by $\Cum(Y_1,\ldots,Y_k)$ its $k$-th order joint cumulant. For
the stationary process $(X_i)_{i\in\Z}$, we write
\begin{equation*}
  \gamma(k_1,k_2,\ldots,k_{d})
   := \Cum(X_0,X_{k_1},X_{k_2},\ldots,X_{k_{d}}).
\end{equation*}
We need the assumption of summability of joint cumulants in
Lemma~\ref{thm:mag}, Lemma~\ref{thm:blk} and
Lemma~\ref{thm:ljung_cleared}. For this reason, we provide a
sufficient condition in Section~\ref{sec:cum}.
\begin{lemma}
  \label{thm:mag}
  Assume $X_i \in \LLL^4$, $\E X_i=0$, $\Theta_2<\infty$ and
  $\sum_{k_1,k_2,k_3\in\Z} |\gamma(k_1,k_2,k_3)|<\infty$.
  For $k,h \geq 1$, $l_n\geq t_n > 0$ and $s_n \in \Z$, set
  $U_k=\sum_{i=1}^{l_n}(X_{i-k}X_i-\gamma_k)$ and
  $V_h=\sum_{j=s_n+1}^{s_n+t_n}(X_{j-h}X_{j}-\gamma_j)$,
  then we have
  \begin{align*}
    |\E (U_k V_h)| \leq t_n \Xi(k,h)
  \end{align*}
  where $\left[\Xi(k,h)_{k,h\geq 1}\right]$ is a symmetric double
  array of non-negative numbers such that $\Xi\in\mathcal{A}_2$, and
  \begin{align*}
    \|\Xi\|_\infty \leq 2\Theta_2^4
    + \sum_{k_1,k_2,k_3\in\Z} |\gamma(k_1,k_2,k_3)|.
  \end{align*}
\end{lemma}

\begin{proof}
  Write
  \begin{align*}
    \E(U_k V_h)
    = & \sum_{i=1}^{l_n}\sum_{j=1}^{t_n}
    \E[(X_{i-k}X_i-\gamma_k)(X_{s_n+j-h}X_{s_n+j}-\gamma_h)]\cr
    = & \sum_{i=1}^{l_n}\sum_{j=1}^{t_n}
     [\gamma(-k,j+s_n-i-h,j+s_n-i) \cr
    &\quad\quad  + \gamma_{j+s_n-i+k-h}\gamma_{j+s_n-i}
     + \gamma_{j+s_n-i+k}\gamma_{j+s_n-i-h} ].
  \end{align*}
  For the sum of the second term, we have
  \begin{align*}
    \left|\sum_{i=1}^{l_n}\sum_{j=1}^{t_n}
     \gamma_{j+s_n-i+k-h}\gamma_{j+s_n-i}\right|
    = & \bigg| \sum_{d=1}^{t_n-1}
     (\gamma_{s_n+d+k-h}\gamma_{s_n+d})(t_n-d) \cr
     &\quad + t_n\sum_{d=t_n-l_n}^0
        \gamma_{s_n+d+k-h}\gamma_{s_n+d}  \cr
    & \quad  + \sum_{d=1-l_n}^{t_n-l_n-1}
         (\gamma_{s_n+d+k-h}\gamma_{s_n+d})(l_n+d)\bigg|\cr
    \leq & t_n\sum_{d \in \Z} |\gamma_{s_n+d+k-h}\gamma_{s_n+d}|
    \cr
    \leq & t_n\sum_{d \in \Z}\zeta_{d+k-h}\zeta_{d}.
  \end{align*}
Similarly, for the sum of the last term
  \begin{align*}
    \left|\sum_{i=1}^{l_n}\sum_{j=1}^{t_n}
       \gamma_{j+s_n-i+k}\gamma_{j+s_n-i-h}\right|
    \leq & t_n\sum_{d\in\Z} \zeta_{d+k+h}\zeta_{d}.
  \end{align*}
  Observe that $\sum_{h=1}^{\infty}\sum_{d \in \Z} \zeta_{d+k-h}
  \zeta_{d} \leq \left(\sum_{d\in\Z}\zeta_d\right)^2 \leq \Theta_2^4$
  and similarly $\sum_{h=1}^{\infty}\sum_{d \in
    \Z}\zeta_{d+k+h}\zeta_{d} \leq \Theta_2^4$. For the sum of the
  first term, it holds that
  \begin{align*}
    \left|\sum_{i=1}^{l_n}\sum_{j=1}^{t_n}
     \gamma(-k,j+s_n-i-h,j+s_n-i) \right|
     \leq t_n\sum_{d\in\Z}|\gamma(-k,d-h,d)|.
  \end{align*}
Utilizing the summability of cumulants, the proof is complete.
\end{proof}

In the proof of Lemma~\ref{thm:blk}, we need the concept of {\em
  indecomposable partitions}. Consider the table
\begin{center}
  \begin{tabular}{lcl}
    $(1,1)$ & \ldots & $(1,J_1)$ \\
    \vdots  &        & \vdots \\
    $(I,1)$ & \ldots & $(I,J_I)$
  \end{tabular}
\end{center}
Denote the $j$-th row of the table by $\vartheta_j$. A partition
$\h{\nu} = \{\nu_1,\ldots,\nu_q\}$ of the table is said to be {\it
  indecomposable} if there are no sets $\nu_{i_1},\ldots,\nu_{i_k}$
($k<q$) and rows $\vartheta_{j_1},\ldots,\vartheta_{j_l}$ ($l<I$) such that $\nu_{i_1}
\cup \cdots \cup \nu_{i_k} = \vartheta_{j_1} \cup \cdots \cup \vartheta_{j_l}$.

\begin{proof}[Proof of Lemma~\ref{thm:blk}]
  Write
  \begin{align*}
    \sum_{k=1}^{s_n} \E_0(\tilde R_{n,k}^2 - \RR_{n,k}^2)
    & = 2\sum_{k=1}^{s_n} \E_0\left[\RR_{n,k}(\tilde
    R_{n,k}-\RR_{n,k})\right]
    + \sum_{k=1}^{s_n} \E_0(\tilde R_{n,k}-\RR_{n,k})^2 \cr
    & =: 2 I_n +  \II_n.
  \end{align*}
  Using Lemma~\ref{thm:ljung_cleared}, we know $\II_n / (n\sqrt{s_n})
  = o_P(1)$. We can express $I_n$ as
  \begin{align}
    I_n = \sum_{a=0}^1\sum_{b=0}^1 I_{n, a b}
    = I_{n, 00} + I_{n, 01} + I_{n, 10} + I_{n, 11}.
    \label{eq:35}
  \end{align}
  where for $a,b=0,1$ (assume without loss of generality that $w_n$ is
  even),
  \begin{align*}
    I_{n, a b} = \sum_{k=1}^{s_n} \E_0\left(\sum_{j=0}^{w_n/2}U_{k,2j-a} \sum_{j=0}^{w_n/2}V_{k,2j-b}\right).
  \end{align*}
  Consider the first term in (\ref{eq:35}), write
  \begin{align*}
    \E (I_{n, 00}^2) &
    = \sum_{k,h=1}^{s_n} \E\left[ \sum_{j=1}^{w_n/2}
      \E_0(U_{k,2j}V_{k,2j})\cdot\E_0(U_{h,2j}V_{h,2j})\right]\cr
    &\quad + \sum_{k,h=1}^{s_n} \sum_{j_1 \neq j_2}
       \E(U_{k,2j_1}U_{h,2j_1})\E(V_{k,2j_2}V_{h,2j_2}) \cr
    &\quad + \sum_{k,h=1}^{s_n} \sum_{j_1 \neq j_2}
      \E(U_{k,2j_1}V_{h,2j_1})\E(V_{k,2j_2}U_{h,2j_2}) \cr
    & :=A_n + B_n + C_n.
  \end{align*}
By Lemma~\ref{thm:mag}, it holds that
  \begin{eqnarray*}
    |B_n|
     &\leq& \sum_{k,h=1}^{s_n} \sum_{j_1,j_2=0}^{w_n/2}
      l_n|K_{2j_2}| \cdot \left[\tilde\Xi(k,h)\right]^2 \cr
    &\leq& w_nl_n \cdot (w_nm_n+2l_n) \sum_{k,h=1}^{s_n}
    \left[\tilde\Xi_n(k,h)\right]^2 = o(n^2s_n),
  \end{eqnarray*}
where $\tilde\Xi_n(k,h)$ is the $\Xi(k,h)$ (defined in
Lemma~\ref{thm:mag}) for the sequence $(\tilde X_i)$. Similarly,
  \begin{eqnarray*}
    |C_n| &\leq& \sum_{k,h=1}^{s_n} \sum_{j_1,j_2=1}^{w_n/2}
    |K_{2j_1}|\cdot |K_{2j_2}|
    \cdot\left[\tilde\Xi_n(k,h)\right]^2\cr
    &\leq& (w_nm_n+l_n)^2 \sum_{k,h=1}^{s_n}
    \left[\tilde\Xi_n(k,h)\right]^2 = o(n^2s_n).
  \end{eqnarray*}
To deal with $A_n$, we express it in terms of cumulants
  \begin{eqnarray*}
    A_n & = &\sum_{k,h=1}^{s_n} \sum_{j=1}^{w_n/2}
     [ \Cum(U_{k,2j},V_{k,2j},U_{h,2j},V_{h,2j}) \cr
    &&\quad\quad\quad\quad\quad
              + \E (U_{k,2j}U_{h,2j})\E (V_{k,2j}V_{h,2j}) \cr
    &&\quad\quad\quad\quad\quad
              + \E (U_{k,2j}V_{h,2j})\E (V_{k,2j}U_{h,2j})]\cr
    & =:& D_n + E_n + F_n.
  \end{eqnarray*}
Apparently $|E_n| = o(n^2 s_n)$ and $|F_n| = o(n^2 s_n)$. Using
the multilinearity of cumulants, we have
\begin{align*}
\Cum(U_{k,2j},V_{k,2j},U_{h,2j},V_{h,2j})
 = \sum_{i_1,i_2\in H_{2j}} \sum_{j_1,j_2\in K_{2j}}
    \Cum(\tilde X_{i_1-k} \tilde X_{i_1},
        \tilde X_{j_1-k} \tilde X_{j_1},
        \tilde X_{i_2-h} \tilde X_{i_2},
        \tilde X_{j_2-h} \tilde X_{j_2})
\end{align*}
for $1 \le k, h \le s_n$. By Theorem~II.2 of
\cite{rosenblatt:1985}, we know
  \begin{align}
    \label{eq:cum_decompose}
    \Cum\left(\tilde X_{i_1-k}\tilde X_{i_1},\tilde X_{j_1-k}\tilde X_{j_1},
      \tilde X_{i_2-h}\tilde X_{i_2},\tilde X_{j_2-h}\tilde X_{j_2}\right)
    = \sum_{\h{\nu}}\prod_{q=1}^b\Cum(\tilde X_i,\;i\in\nu_q)
  \end{align}
where the sum is over all indecomposable partitions
$\h{\nu}=\{\nu_1,\ldots,\nu_q\}$ of the table
  \begin{center}
    \begin{tabular}{ll}
      $i_1-k$ & $i_1$ \\
      $j_1-k$ & $j_1$ \\
      $i_2-h$ & $i_2$ \\
      $j_2-h$ & $j_2$
    \end{tabular}
  \end{center}
  By Theorem~\ref{thm:cum}, the condition $\sum_{k=0}^\infty
  k^6\delta_8(k)<\infty$ implies that all the joint cumulants up to
  order eight are absolutely summable. Therefore, using
  Proposition~\ref{thm:array}, we know
  \begin{align*}
    \sum_{k,h=1}^{s_n} \left|\Cum(U_{k,2j},V_{k,2j},U_{h,2j},V_{h,2j})\right| = O(|K_{2j}|s_n^2),
  \end{align*}
and it follows that $|D_n| = O\left((w_nm_n+l_n)s_n^2\right) =
o(n^2s_n). $ We have shown that $\E (I_{n, 00}^2) = o(n^2s_n), $
which, in conjunction with similar results for the other three
terms in (\ref{eq:35}), implies that $\E (I_n^2) = o(n^2 s_n)$ and
hence $I_n/(n\sqrt{s_n}) = o_P(1)$. The proof is now complete.
\end{proof}

\subsubsection{Step 3: Central limit theorem concerning $\RR_{n,k}$'s.}
\label{sec:proof_cleared}

\begin{proof}[Proof of Lemma~\ref{thm:ljung_cleared}]
  Let $\Upsilon_n(k,h) := \E (U_{k,1} U_{h,1})$ and $\upsilon_n(k,h)
  := \Upsilon_n(k,h)/l_n$. By Lemma~\ref{thm:mag} we know
  $|\upsilon_n(k,h)|\leq \tilde \Xi_n(k,h)$.  Write
  \begin{align*}
    \sum_{k=1}^{s_n} \E_0 \RR_{n,k}^2
    = & \sum_{k=1}^{s_n} \left[ \sum_{j=1}^{w_n}\left(U_{k,j}^2-\Upsilon_n(k,k)\right)
      + 2\sum_{j=1}^{w_n}\left(U_{k,j}\sum_{l=1}^{j-1}U_{k,l}\right) \right] \cr
    = & \sum_{j=1}^{w_n} \left[\sum_{k=1}^{s_n} \left(U_{k,j}^2 - \Upsilon_n(k,k)\right)\right]
    + 2
    \sum_{j=1}^{w_n}\left(\sum_{k=1}^{s_n}U_{k,j}\sum_{l=1}^{j-1}U_{k,l}\right).
  \end{align*}
Using similar a argument as the one for dealing with the term
$A_n$ in Lemma~\ref{thm:blk}, we know
  \begin{align*}
    \sum_{j=1}^{w_n}\left\|\sum_{k=1}^{s_n}\left(U_{k,j}^2
       -\Upsilon_n(k,k)\right)\right\|^2 = o(n^2s_n),
  \end{align*}
and it follows that
  \begin{align*}
    \frac{1}{n\sqrt{s_n}}\sum_{j=1}^{w_n}
    \left[\sum_{k=1}^{s_n} \left(U_{k,j}^2 - \Upsilon_n(k,k)\right)\right]
    = o_P(1).
  \end{align*}
  Therefore, it suffices to consider
  \begin{align*}
     \sum_{j=1}^{w_n}\left(\sum_{k=1}^{s_n}U_{k,j}
     \sum_{l=1}^{j-1}U_{k,l}\right)
     =:\sum_{j=1}^{w_n} D_{n,j}.
  \end{align*}
Let $\mathcal{G}_{n,j}=\langle D_{n,1},\ldots,D_{n,j}\rangle$. Observe
that $(D_{n,j})$ is a martingale difference sequence with respect
to $(\mathcal{G}_{n,j})$. We shall apply the martingale central limit
theorem. Write
  \begin{align*}
    \E\left(D_{n,j}^2|\mathcal{G}_{n,j-1}\right) - \E D_{n,j}^2
    &= \sum_{k,h=1}^{s_n} \Upsilon_n(k,h)
     \left(\sum_{l=1}^{j-1}U_{k,l}\sum_{l=1}^{j-1}U_{h,l}
        - (j-1)\Upsilon_n(k,h)\right)\cr
    & =  \sum_{k,h=1}^{s_n} \Upsilon_n(k,h)
     \left(\sum_{l=1}^{j-1}U_{k,l}U_{h,l}
          -(j-1)\Upsilon_n(k,h)\right) \cr
    & + \sum_{k,h=1}^{s_n} \Upsilon_n(k,h)
      \left(\sum_{l=1}^{j-1}U_{k,l}\sum_{q=1}^{l-1}U_{h,q}
      + \sum_{l=1}^{j-1}U_{h,l}\sum_{q=1}^{l-1}U_{k,q}\right)\cr
    & =: I_{n,j} + \II_{n,j}
  \end{align*}
  For the first term, by  Lemma~\ref{thm:mag}, we have
  % \begin{align*}
  %   \left\|\sum_{j=1}^{w_n} I_{n,j}\right\|_{p/4}
  %   & = \left\|\sum_{j=1}^{w_n-1} (w_n-j) \sum_{k,h=1}^{s_n} \Upsilon_n(k,h)
  %     \left[U_{k,j}U_{h,j}-\Upsilon_n(k,h)\right]\right\|_{p/4}\cr
  %   & \leq \CC_{p/4}\left(\sum_{j=1}^{w_n-1}(w_n-j)^{p/4}\right)^{4/p}
  %   \sum_{1\leq k,h\leq s_n}|\Upsilon_n(k,h)|\left\|(U_{k,j}U_{h,j}-\Upsilon_n(k,h))\right\|_{p/4}\cr
  %   & \leq \CC_{p/4} w_n^2 l_n^2 w_n^{4/p-1} 8\Theta_p^2 \cdot \sum_{1\leq k,h\leq s_n}\upsilon_n(k,h)
  %   = o(n^2s_n).
  % \end{align*}
  \begin{align*}
    \left\|\sum_{j=1}^{w_n} I_{n,j}\right\|^2 = & \left\|\sum_{j=1}^{w_n-1} (w_n-j) \sum_{k,h=1}^{s_n} \Upsilon_n(k,h)
      \left[U_{k,j}U_{h,j}-\Upsilon_n(k,h)\right]\right\|^2\cr
    = & \sum_{j=1}^{w_n-1}(w_n-j)^2\left[\sum_{k,h} |\Upsilon_n(k,h)|\left\|(U_{k,j}U_{h,j}-\Upsilon_n(k,h))\right\|\right]^2\cr
    \leq & w_n^3l_n^4 \left[\sum_{k,h} |\upsilon_n(k,h)| \cdot 4\Theta_8^2 \right]^2 = o(n^4s_n^2).
  \end{align*}
  Using Lemma~\ref{thm:mag} and Proposition~\ref{thm:array}, we obtain
  \begin{align*}
    & \left\|\sum_{j=1}^{w_n} \II_{n,j}\right\|^2 = \left\|\sum_{j=1}^{w_n-1} (w_n-j)
      \sum_{k,h}\Upsilon_n(k,h)\left(U_{k,j}\sum_{l=1}^{j-1}U_{h,l} + U_{h,j}\sum_{l=1}^{j-1}U_{k,l}\right)\right\|^2 \cr
    = & 2\sum_{j=1}^{w_n-1}(w_n-j)^2(j-1)\sum_{1\leq k_1,h_1,k_2,h_2\leq s_n}\Upsilon_n(k_1,h_1)\Upsilon_n(k_2,h_2)
    \left[\Upsilon_n(k_1,k_2)\Upsilon_n(h_1,h_2) + \Upsilon_n(k_1,h_2)\Upsilon_n(h_1,k_2)\right]\cr
    \leq & 4 n^4 \sum_{1\leq k_1,h_1,k_2,h_2\leq s_n}
    \left|\upsilon_n(k_1,h_1)\upsilon_n(h_1,h_2)\upsilon_n(h_2,k_2)\upsilon_n(k_2,k_1)\right| = O(n^4s_n) = o(n^4s_n^2).
  \end{align*}
  Therefore, we have
  \begin{align*}
    \frac{1}{n^2s_n}\left[\sum_{j=1}^{w_n}\E\left(D_{n,j}^2|\mathcal{G}_{n,j-1}\right)
      -\sum_{j=1}^{w_n}\E D_{n,j}^2 \right] \stackrel{p}{\to} 0.
  \end{align*}
  Using Lemma~\ref{thm:mag} and Lemma~\ref{thm:covconvergence}, we know
  \begin{align*}
    \frac{1}{n^2s_n}\sum_{j=1}^{w_n}\E D_{n,j}^2 = \frac{1}{2n^2s_n} w_n(w_n-1)l_n^2\sum_{k,h=1}^{s_n}[\upsilon_n(k,h)]^2
    \to \frac{1}{2}\sum_{k\in\Z}\sigma_k^2,
  \end{align*}
  and it follows that
  \begin{align}
    \label{eq:varcon}
    \frac{1}{n^2s_n}\sum_{j=1}^{w_n}\E\left(D_{n,j}^2|\mathcal{G}_{n,j-1}\right)
    \stackrel{p}{\to} \frac{1}{2}\sum_{k\in\Z}\sigma_k^2.
  \end{align}
  To verify the Lindeberg condition, we compute
\begin{align*}
    \E D_{n,j}^4 = & \sum_{k_1,k_2,k_3,k_4=1}^{s_n}
    \E\left(U_{k_1,j}U_{k_2,j}U_{k_3,j}U_{k_4,j}\right)\cr
    & \times \E\left[\left(\sum_{l=1}^{j-1}U_{k_1,l}\right) \left(\sum_{l=1}^{j-1}U_{k_2,l}\right)
      \left(\sum_{l=1}^{j-1}U_{k_3,l}\right)\left(\sum_{l=1}^{j-1}U_{k_4,l}\right)\right]\cr
    \leq & \sum_{k_1,k_2,k_3,k_4=1}^{s_n} \left|\E(U_{k_1,j}U_{k_2,j}U_{k_3,j}U_{k_4,j})\right| \cdot
    2\mathcal{C}_4^4(j-1)^2l_n^2\Theta_8^8
\end{align*}
  We express $\E(U_{k_1,1}U_{k_2,1}U_{k_3,1}U_{k_4,1})$ in terms of cumulants
  \begin{align*}
    \E(U_{k_1,1}U_{k_2,1}U_{k_3,1}U_{k_4,1}) & = \Cum(U_{k_1,1},U_{k_2,1},U_{k_3,1},U_{k_4,1})
    + \E(U_{k_1,1}U_{k_2,1})\E(U_{k_3,1}U_{k_4,1}) \cr
    & \quad + \E(U_{k_1,1}U_{k_3,1})\E(U_{k_2,1}U_{k_4,1}) + \E(U_{k_1,1}U_{k_4,1})\E(U_{k_2,1}U_{k_3,1})\cr
    & =:A_n+B_n+E_n+F_n
  \end{align*}
  From Lemma~\ref{thm:mag}, it is easily seen that
  \begin{align*}
    \sum_{k_1,k_2,k_3,k_4=1}^{s_n} |B_n| \leq l_n^2 \sum_{k_1,k_2,k_3,k_4=1}^{s_n} \tilde\Xi_n(k_1,k_2)\cdot\tilde\Xi_n(k_3,k_4)
    = O(l_n^2s_n^2),
  \end{align*}
  and similarly $\sum_{k_1,k_2,k_3,k_4=1}^{s_n} |E_n| = O(l_n^2s_n^2)$
  and $\sum_{k_1,k_2,k_3,k_4=1}^{s_n} |F_n| = O(l_n^2s_n^2)$. By multilinearity of cumulants,
  \begin{align*}
    A_n = \sum_{i_1,i_2,i_3,i_4=1}^{l_n}
    \Cum(\tilde X_{i_1-k_1}\tilde X_{i_1},\tilde X_{i_2-k_2}\tilde X_{i_2},
    \tilde X_{i_3-k_3}\tilde X_{i_3},\tilde X_{i_4-k_4}\tilde
    X_{i_4}).
  \end{align*}
Each cumulant in the preceding equation is to be further
simplified similarly as (\ref{eq:cum_decompose}). Using
summability of joint cumulants up to order eight and
Proposition~\ref{thm:array}, we have
  \begin{align*}
    \sum_{k_1,k_2,k_3,k_4=1}^{s_n} |A_n|
     = O(l_ns_n^3) =o(l_n^2s_n^2).
  \end{align*}
Using orders obtained for $|A_n|$, $|B_n|$, $|E_n|$ and $|F_n|$,
we obtain $\sum_{j=1}^{w_n} \E D_{n,j}^4 = o(n^4s_n^2)$. Then, by
(\ref{eq:varcon}), we can apply Corollary 3.1. of \cite{hall:1980}
to obtain
  \begin{align*}
    \frac{1}{n\sqrt{s_n}} \sum_{j=1}^{w_n} D_{n,j}
     \Rightarrow \mathcal{N}\left(0,
      \frac{1}{2}\sum_{k\in\Z}\sigma_k^2\right),
  \end{align*}
and the lemma follows.
\end{proof}

\subsubsection{Proof of Theorem~\ref{thm:ljung_power}}
\label{sec:proof_power}

\begin{proof}[Proof of Theorem~\ref{thm:ljung_power}]
We shall only prove (\ref{eq:ljung_power}), since
(\ref{eq:var_estimation}) can be obtained by very similar
arguments. Write $\hat\gamma_k = \E_0 \hat\gamma_k + \gamma_k -
(\gamma_k-\E\hat\gamma_k)$, and hence
  \begin{align*}
    \sum_{k=1}^{s_n} (\hat\gamma_k^2 - \gamma_k^2)
    & = 2\sum_{k=1}^{s_n} \gamma_k \E_0 \hat\gamma_k
    + \sum_{k=1}^{s_n} ( \E_0\hat\gamma_k)^2
    - 2\sum_{k=1}^{s_n} \frac{k}{n}\gamma_k \E_0\hat\gamma_k
    -  2\sum_{k=1}^{s_n} \frac{k}{n}\gamma_k^2
    + \sum_{k=1}^{s_n}\frac{k^2}{n^2}\gamma_k^2 \cr
    & = : 2I_n + \II_n + \III_n + \IV_n + V_n.
  \end{align*}
  Using the conditions $\Theta_4<\infty$ and $s_n=o(\sqrt{n})$, it is
  easily seen that $\sqrt{n}\IV_n\to 0$ and
  $\sqrt{n}V_n\to 0$. Furthermore
  \begin{align*}
    \sqrt{n}\|\III_n\| \leq 2\sqrt{n}\sum_{k=1}^{s_n} \frac{k}{n}|\gamma_k|\cdot \frac{2\Theta_4^2}{\sqrt{n}} \to 0
    \quad \hbox{and} \quad \sqrt{n}\E\II_n\leq \sqrt{n}\sum_{k=1}^{s_n}\frac{4\Theta_4^4}{n} \to 0.
  \end{align*}
  Define $Y_i=\sum_{k=1}^\infty \gamma_kX_{i-k}$. For the term $I_n$, write
  \begin{align*}
    nI_n & = \sum_{i=1}^n \E_0(X_iY_i) - \sum_{i=1}^n
     \E_0\left(X_i\sum_{k=s_n+1}^\infty\gamma_kX_{i-k}\right)
    + \sum_{k=1}^{s_n}
     \gamma_k\left(\sum_{i=1}^{k}(X_{i-k}X_i-\gamma_k)\right) \cr
    & =: A_n + B_n + E_n
  \end{align*}
  Clearly $\|E_n\|/\sqrt{n} \leq \sum_{k=1}^{s_n} |\gamma_k|
  2\Theta_4^2\sqrt{k}/\sqrt{n} \to 0$. Define
  $W_{n,i}=X_i\sum_{k=s_n+1}^\infty\gamma_kX_{i-k}$, then
  \begin{align*}
    \|\PP^0 W_{n,i}\| \leq \left\{
    \begin{array}{ll}
      \delta_{4}(i)\cdot\Theta_4\sum_{k=s_n+1}^\infty |\gamma_k|
       & \hbox{if}\quad 0 \leq i \leq s_n \\
      \Theta_4\delta_{4}(i)\sum_{k=s_n+1}^\infty |\gamma_k|
       + \Theta_4\sum_{k=s_n+1}^i |\gamma_k|\delta_4(i-k)
      & \hbox{if}\quad i>s_n.
    \end{array}\right.
  \end{align*}
  It follows that
  \begin{align*}
    \|B_n/\sqrt{n}\| \leq 2\Theta_4^2\sum_{k=s_n+1}^\infty |\gamma_k| \to 0.
  \end{align*}
  Set $Z_i=X_iY_i$, then $(Z_i)$ is a stationary process of the form
  (\ref{eq:wold}). Furthermore
  \begin{align*}
    \|\PP^0 Z_i\| \leq  \delta_{4}(i)\cdot\Theta_4\sum_{k=1}^\infty |\gamma_k|
      + \Theta_4\sum_{k=1}^i |\gamma_k|\delta_4(i-k).
  \end{align*}
  Since $\sum_{i=0}^\infty \|\PP^0 Z_i\|<\infty$, utilizing Theorem 1
  in \cite{hannan:1973} we have $A_n/\sqrt{n} \Rightarrow
  \mathcal{N}(0,\|D_0\|^2)$, and then (\ref{eq:ljung_power}) follows.
\end{proof}

\subsubsection{Proof of Corollary~\ref{thm:ljung_corr} and \ref{thm:ljung_corr_power}}
\label{sec:proof_corollary}

\begin{proof}[Proof of Corollary~\ref{thm:ljung_corr} and
  \ref{thm:ljung_corr_power}] By (\ref{eq:fact3}), we
  know $\|n\bar X_n\|_4 \le \sqrt{3n} \Theta_4$, and it follows that
  \begin{align*}
    \left\|\sum_{i=k+1}^n(X_{i-k}-\bar X_n)(X_i - \bar X_n)
     - \sum_{i=k+1}^n X_{i-k}X_i\right\|
    \leq 9 \Theta_4^2.
  \end{align*}
  % and $n\|\hat\gamma_0-\gamma_0\| = \left\|\sum_{i=1}^n(X_i^2-\gamma_0) 
  % - n\bar X_n^2\right\| \leq (2\sqrt{n}+1)\Theta_4^2$. 
  Theorem~\ref{thm:ljung} holds for $\breve\gamma_k$ because
  \begin{align*}
    \frac{n}{\sqrt{s_n}} \sum_{k=1}^{s_n}
    \E\left|(\hat\gamma_k-\E\hat\gamma_k)^2
      - (\breve\gamma_k-\E\hat\gamma_k)^2\right|
    & \leq \frac{n}{\sqrt{s_n}} \sum_{k=1}^{s_n}
      \left\|\hat\gamma_k+\breve\gamma_k-2\E\hat\gamma_k\right\|
    \cdot \left\|\hat\gamma_k-\breve\gamma_k\right\| \cr
    & \leq \frac{n}{\sqrt{s_n}} \sum_{k=1}^{s_n}
     \left(\frac{4\Theta_4^2}{\sqrt{n}} + \frac{9\Theta_4^2}{n}\right)
    \cdot \frac{9\Theta_4^2}{n} \to 0.
  \end{align*}
  In Theorem~\ref{thm:ljung_power}, (\ref{eq:ljung_power}) holds with
  $\hat \gamma_k$ replaced by $\breve \gamma_k$ because
  \begin{align*}
    \sqrt{n}\sum_{k=1}^{s_n}
       \E\left|\hat\gamma_k^2-\breve\gamma_k^2\right|
    \leq \sqrt{n}\sum_{k=1}^{s_n}
     \|\hat\gamma_k+\breve\gamma_k\|\cdot\|\hat\gamma_k-\breve\gamma_k\|
    \leq \sqrt{n}\sum_{k=1}^{s_n} \left(2|\gamma_k|
    + \frac{4\Theta_4^2}{\sqrt{n}} + \frac{9\Theta_4^2}{n}\right)
     \frac{9\Theta_4^2}{n} \to 0,
  \end{align*}
  and (\ref{eq:var_estimation}) can be proved similarly.  Now we turn
  to the sample autocorrelations. Write
  \begin{align*}
    \sum_{k=1}^{s_n} \left\{[\hat r_k - (1-k/n)r_k]^2
     - [\hat\gamma_k/\gamma_0 - (1-k/n)r_k]^2\right\}
    & = \sum_{k=1}^{s_n}
     {{2 (\E_0\hat\gamma_k)[\hat\gamma_k(\gamma_0-\hat\gamma_0)]}
     \over{\gamma_0^2\hat\gamma_0}}
    + {{\hat\gamma_k^2(\gamma_0-\hat\gamma_0)^2}
     \over {\gamma_0^2\hat\gamma_0^2}}.
  \end{align*}
  Since
  \begin{align*}
    \sum_{k=1}^{s_n} \E\left| (\E_0\hat\gamma_k)
     \hat\gamma_k(\gamma_0-\hat\gamma_0)\right|
    \leq \sum_{k=1}^{s_n} 2\CC_3\Theta_6^2\frac{1}{\sqrt{n}} \cdot
     \left(|\gamma_k|+2\CC_3\Theta_6^2\frac{1}{\sqrt{n}}\right)
    \cdot 2\CC_3\Theta_6^2\frac{1}{\sqrt{n}}
    = o\left(\frac{\sqrt{s_n}}{n}\right)
  \end{align*}
  and similarly $\sum_{k=1}^{s_n}
  \E\left|\hat\gamma_k^2(\gamma_0-\hat\gamma_0)^2\right| =
  o(\sqrt{s_n}/n)$, (\ref{eq:ljung_corr}) follows by applying the
  Slutsky theorem. To show the limit theorems in
  Corollary~\ref{thm:ljung_corr_power}, note that using the
  Cramer-Wold device, we have
  \begin{eqnarray*}
  \left[\sqrt{n}(\hat\gamma_0^2-\gamma_0^2),
  \sqrt{n}\left(\sum_{k=1}^{s_n}\hat\gamma_k^2
      - \sum_{k=1}^{s_n}\gamma_k^2\right)\right]
  \end{eqnarray*}
  converges to a bivariate normal distribution. Then
  Corollary~\ref{thm:ljung_corr_power} follows by applying the delta
  method.
\end{proof}

\section{A Normal Comparison Principle}
\label{sec:normalcomparison} In this section we shall control tail
probabilities of Gaussian vectors by using their covariance
matrices. Denote by $\varphi_d ((r_{ij}); x_1,\ldots, x_d)$ the
density of a $d$-dimensional multivariate normal random vector
$\h{X}=(X_1,\ldots,X_d)^\top$ with mean zero and covariance matrix
$(r_{ij})$, where we always assume $r_{ii}=1$ for $1 \leq i \leq
d$ and $(r_{ij})$ is nonsingular.
% For a index subset $A \subset [d]$ with size $|A|$, we
% use $\varphi_{|A|}((r_{ij}),X_k=x_k, k\in A)$ to denote the marginal
% density of the sub-vector $X_A:=(X_k)_{k \in A}$.
For $1\leq h < l \leq d$, we use
$\varphi_{2}((r_{ij});X_h=x_h,X_l=x_l)$ to denote the marginal
density of the sub-vector $(X_h,X_l)^\top$. Let
$$Q_d\left((r_{ij});z_1,\ldots,z_d\right) = \int_{z_1}^\infty
\cdots \int_{z_d}^{\infty}
\varphi_d\left((r_{ij}),x_1,\ldots,x_d\right) \dd x_{d} \cdots \dd
x_{1}.$$ The partial derivative with respect to $r_{hl}$ is
obtained similarly as equation (3.6) of \cite{berman:1964} by
using equation (3) of \cite{plackett:1954}
\begin{align}
  \label{eq:deriv}
   & \frac{\partial Q_d\left((r_{ij});z_1,\ldots,z_d\right)}{\partial
    r_{hl}} \cr
  & = \left(\prod_{k \neq h,l}\int_{z_k}^\infty\right)
  \varphi_d\left((r_{ij});x_1,\ldots,x_{h-1},
  z_h,x_{h+1},\ldots,x_{l-1},z_l,x_{l+1},\ldots,x_d\right)
  \prod_{k \neq h,l} \dd x_k.
\end{align}
where $\left(\prod_{k \neq h,l}\int_{z_k}^\infty\right)$ stands
for $\int_{z_1}^\infty \cdots \int_{z_{h-1}}^\infty
\int_{z_{h+1}}^\infty \cdots \int_{z_{l-1}}^\infty
\int_{z_{l+1}}^\infty \cdots \int_{z_d}^\infty$. If all the $z_k$
have the same value $z$, we use the simplified notation
$Q_d\left((r_{ij});z\right)$ and $\partial Q_d
((r_{ij});z)/\partial r_{hl}$. The following simple facts about
conditional distribution will be useful. For four different
indicies $1 \leq h,l,k,m \leq d$, we have
\begin{align}
\label{eq:condmean}
  \E(X_k|X_h=X_l=z) & = \frac{r_{kh}+r_{kl}}{1+r_{hl}}z,
 \\ \label{eq:condvar}
\Var(X_k|X_h=X_l=z) & =
\frac{1-r_{hl}^2-r_{kh}^2-r_{kl}^2+2r_{hl}r_{kh}r_{kl}}
     {1-r_{hl}^2},
\\ \label{eq:condcov}
\Cov(X_k,X_m|X_h=X_l=z) & = r_{km} -
\frac{r_{hk}r_{hm}+r_{lk}r_{lm}-r_{hl}r_{hk}r_{lm}-r_{hl}r_{hm}r_{lk}}
     {1-r_{hl}^2}.
\end{align}

\begin{lemma}
  \label{thm:normbdd}
  For every $z>0$, $0<s<1$, $d \geq 1$ and $\epsilon>0$, there exists positive
  constants $C_d$ and $\epsilon_d$ such that for $0 < \epsilon < \epsilon_d$
  \begin{enumerate}
  \item if $|r_{ij}| < \epsilon$ for all $1 \leq i < j \leq d$, then
    \begin{align}
      \label{eq:normbdd1}
      Q_d\left((r_{ij});z\right) & \leq C_d
      \exp\left\{-\left(\frac{d}{2} - C_d \epsilon \right)z^2\right\} \\ \label{eq:normbdd3}
      Q_{d}\left((r_{ij});z,\ldots,z\right) & \leq C_d \,f_d(\epsilon,1/z)\,
      \exp\left\{-\left(\frac{d}{2} - C_d \epsilon
        \right)z^2\right\} \\ \label{eq:normbdd2}
      Q_d\left((r_{ij});sz,z,\ldots,z\right) & \leq C_d
      \exp\left\{-\left(\frac{s^2+d-1}{2} - C_d \epsilon
        \right)z^2\right\}
    \end{align}
    where $f_{2k}(x,y)=\sum_{l=0}^k x^ly^{2(k-l)}$ and
    $f_{2k-1}(x,y)=\sum_{l=0}^{k-1}x^ly^{2(k-l)-1}$ for $k \geq 1$;
  \item if for all $1 \leq i < j \leq d+1$ such that $(i,j)\neq(1,2)$,
    $|r_{ij}| \leq \epsilon$, then
    \begin{equation}
      \label{eq:normbdd4}
      Q_{d+1}\left((r_{ij});z\right) \leq C_d
      \exp\left\{-\left(\frac{(1-|r_{12}|)^2+d}{2}- C_d\epsilon\right)z^2\right\}.
    \end{equation}
  \end{enumerate}
% In these inequalities, $C$ is a constant which depends on $\epsilon$
% and $d$, but is independent of $z$.
\end{lemma}

\begin{proof}
The following facts about normal tail probabilities are
well-known:
  \begin{equation}
    \label{eq:normbdd}
    P(X_1\geq x) \leq \frac{1}{\sqrt{2\pi}x} e^{-x^2/2} \hbox{ for } x>0
    \quad \hbox{and} \quad
    \lim_{x \to \infty}\frac{P(X_1 \geq x)}{(1/x)(2\pi)^{-1/2}
    \exp\left\{-x^2/2\right\}}=1,
  \end{equation}
  By (\ref{eq:normbdd}), the inequalities
  (\ref{eq:normbdd1}) -- (\ref{eq:normbdd2}) with $\epsilon=0$ are
  true for the random vector with iid standard normal entries. The
  idea is to compare the desired probability with the corresponding
  one for such a vector. We first prove (\ref{eq:normbdd1}) by
  induction. When $d=1$, the inequality is trivially true. When $d=2$,
  by (\ref{eq:deriv}), there exists a number $r_{12}'$ between $0$ and
  $r_{12}$ such that
\begin{eqnarray*}
 |Q_2((r_{ij});z) - Q_2(I_2;z)| &\leq& \varphi((r_{ij}'),z,z)
  |r_{12}| \cr
  &\leq& C \exp\left\{-\frac{z^2}{1+|r_{12}'|}\right\} \leq C
  \exp\left\{-(1-\epsilon)z^2\right\},
\end{eqnarray*}
which, together with $Q_2(I_2;z) \leq C \exp\{-z^2\}$, implies
(\ref{eq:normbdd1}) for $d=2$ with $\epsilon_2=1/2$ and some $C_2
> 1$. Now for $d \geq 3$, assume (\ref{eq:normbdd1}) holds for all
dimensions less than $d$. There exists a matrix $(r_{ij}') =
\theta(r_{ij}) + (1-\theta) I_d$ for some $0<\theta<1$ such that
  \begin{equation}
    \label{eq:2}
    Q_d\left((r_{ij});z\right) -
    Q_d\left((I_d;z\right) = \sum_{1 \leq h,l \leq d}
    \frac{\partial Q_d}{\partial r_{hl}}((r_{ij}');z,\ldots,z) r_{hl}.
  \end{equation}
By (\ref{eq:condmean}), $\E(X_k|X_h=X_l=z) \leq
2\epsilon'z/(1-\epsilon')$ for $k \neq h,l$. Therefore, by
writing the density in (\ref{eq:deriv}) as the product of the
density of $(X_h,X_l)$ and the conditional density of
$\h{X}_{-\{h,l\}}$ given $X_h=X_l=z$, where $\h{X}_{-\{h,l\}}$
denotes the sub-vector
$(X_1,\ldots,X_{h-1},X_{h+1},\ldots,X_{l-1},X_{l+1},\ldots,X_d)^\top$;
we have
\begin{equation}
  \label{eq:3}
  \left|\frac{\partial Q_d}{\partial
      r_{hl}}((r_{ij}');z,\ldots,z)\right| \leq
  \varphi_2((r_{ij}');X_h=X_l=z)
  Q_{d-2}((r_{ij|hl}');(1-3\epsilon)z),
\end{equation}
where $(r_{ij|hl}')$ is the correlation matrix of the conditional
distribution of $\h{X}_{-\{h,l\}}$ given $X_h$ and $X_l$. By
(\ref{eq:condvar}) and (\ref{eq:condcov}), we know for $k, m \in [d]\setminus\{h,l\}$
and $k \neq m$,
$$\Var(X_k|X_h=X_l=z) \geq 1-3\epsilon^2-2\epsilon^3 \quad
\hbox{and} \quad
\Cov(X_k,X_m|X_h=X_l=z) \leq \frac{\epsilon(1+\epsilon)}{1-\epsilon}.$$
Therefore, all the off-diagonal entries of  $(r_{ij|hl}')$ are less
than $2\epsilon$ if we let $\epsilon<1/5$. Applying the induction hypothesis,
if $2\epsilon<\epsilon_{d-2}$, then
$$Q_{d-2}((r_{ij|hl}');(1-3\epsilon)z) \leq C_{d-2}
\exp\left\{-\left(\frac{d-2}{2}-2C_{d-2}\epsilon\right)(1-3\epsilon)^2z^2\right\},$$
and equation (\ref{eq:3}) becomes
$$\left|\frac{\partial Q_d}{\partial
      r_{hl}}((r_{ij}');z,\ldots,z)\right| \leq CC_{d-2}
  \exp\left\{-(1-\epsilon)z^2\right\} \cdot
  \exp\left\{-\left(\frac{d-2}{2}-\left(2C_{d-2}+3(d-2)\right)\epsilon\right)z^2\right\}.$$
Therefore, (\ref{eq:normbdd1}) holds for
$\epsilon_d<\min\{1/5,\epsilon_{d-2}/2\}$ and some $C_d > 2C_{d-2}+3(d-2)+1$.

% Now let us consider (\ref{eq:normbdd2}). Again it is trivial when
% $d=1$. When $d=2$, we have for some $r_{12}'$ between $0$ and $r_{12}$,
% $$|Q_2((r_{ij});sz,z) - Q_2(I_2;sz,z)| \leq \varphi((r_{ij}'),z,sz)
%   |r_{12}| \leq C \exp\left\{-\frac{1+2r_{12}'s+s^2}{2(1-r_{12}'^2)}z^2\right\} \leq C
%   \exp\left\{-\frac{(1+s^2)z^2}{1+\epsilon'}\right\},$$
Using very similar arguments, inequality (\ref{eq:normbdd2}) can be
proved by applying (\ref{eq:normbdd1}); and inequality
(\ref{eq:normbdd4}) can be obtained by employing both
(\ref{eq:normbdd1}) and (\ref{eq:normbdd2}). To prove inequality
(\ref{eq:normbdd3}), which is a refinement of (\ref{eq:normbdd1}), it
suffices to observe that, by (\ref{eq:normbdd}), (\ref{eq:2}) and (\ref{eq:3})
\begin{align*}
  Q_d\left((r_{ij});z\right) & \leq Q_d(I_d;z) + \sum_{1\leq h,l\leq
    d} C\,\epsilon\, \exp\{-(1-\epsilon)z^2\}
  Q_{d-2}((r_{ij|hl}');(1-3\epsilon)z) \\
  & \leq C_d \frac{1}{z^d} \exp\left\{\frac{dz^2}{2}\right\} + C_d \,\epsilon\,
  \exp\{-(1-\epsilon)z^2\} \sum_{1\leq h,l\leq
    d}
  Q_{d-2}((r_{ij|hl}');(1-3\epsilon)z);
\end{align*}
and apply the induction argument.
\end{proof}

\begin{lemma}
  \label{thm:normalcomparison}
  Let $(X_n)$ be a stationary mean zero Gaussian process. Let
  $r_k=Cov(X_0,X_k)$. Assume $r_0=1$, and $\lim_{n\to\infty} r_n(\log
  n) =0$. Let $a_n=(2\log n)^{-1/2}$, $b_n=(2\log n)^{1/2} - (8\log
  n)^{-1/2}(\log\log n + \log 4\pi)$, and $z_n=a_nz+b_n$ for $z \in
  \R$. Define the event $A_i=\{X_i \geq z_n\}$, and
  $$Q_{n,d} = \sum_{1\leq i_1<\ldots<i_d\leq n}P(A_{i_1}\cap \cdots
  \cap A_{i_d}).$$
Then $\lim_{n \to \infty} Q_{n,d} = e^{-dz}/d\,!$ for all $d \geq
1$.
\end{lemma}

\begin{proof}
  Note that $z_n^2 = 2\log n - \log\log n -\log (4\pi) + 2z +
  o(1)$. If $(X_n)$ consists of iid random variables, by the equality in
  (\ref{eq:normbdd}),
  \begin{eqnarray*}
    \lim_{n \to \infty} Q_{n,d}
     &=& \lim_{n\to\infty} {n
    \choose d} Q_d(I_d,z_n) \cr
    &=& \lim_{n \to \infty} {n \choose d}
    \frac{1}{(2\pi)^{d/2}z_n^d}
          \exp\left\{-\frac{dz_n^2}{2}\right\}
     = \frac{e^{-dz}}{d!}.
  \end{eqnarray*}
When the $X_n$'s are dependent, the result is still trivially true
when $d=1$. Now we deal with the $d\geq 2$ case. Let $\gamma_k =
\sup_{j \geq k}|r_j|$, then $\gamma_1<1$ by stationarity, and
$\lim_{n\to\infty}\gamma_n\log  n=0$. Consider an ordered subset
  $J=\left\{t,{t+l_1},{t+l_1+l_2},\ldots,{t+l_1+\cdots+l_{d-1}}\right\}
  \subset [n]$, where $l_1,\ldots,l_{d-1} \geq 1$.
  % Suppose $l_{\ell_1},\ldots,l_{\ell_a} \leq L$, and all the other $l_j$ are
  % larger than $L$.
  We define an equivalence relation $\sim$ on $J$ by saying $k \sim j$
  if there exists $k_1,\ldots,k_p \in J$ such that
  $k=k_1<k_2<\cdots<k_p=j$, and $k_h-k_{h-1} \leq L$ for $2 \leq h
  \leq p$. For any $L\geq 2$, denote by $s(J,L)$ the number of $l_j$
  which are less than or equal to $L$. To similify the notation, we
  sometimes use $s$ instead of $s(J,L)$. $J$ is divided into $d-s$
  equivalence classes
  $\mathcal{B}_1,\ldots,\mathcal{B}_{d-s}$. Suppose $s \geq 1$, assume
  w.l.o.g. that $|\mathcal{B}_1|\geq 2$. Pick $k_0,k_1 \in \mathcal{B}_1$,
  and $k_p \in \mathcal{B}_p$ for $2\leq p \leq d-s$, and set
  $K=\{k_0,k_1,k_2,\ldots,k_{d-s}\}$. Define $Q_J=P(\cap_{k \in
    J}A_k)$ and $Q_K$ similarly, then $Q_J \leq Q_K$. By
  (\ref{eq:normbdd4}) of Lemma~\ref{thm:normbdd}, there exists a
  number $M>1$ depending on $d$ and the sequence $(\gamma_k)$, such
  that when $L>M$,
\begin{eqnarray*}
Q_K &\leq&  C_{d-s}
\exp\left\{-\left(\frac{(1-\gamma_1)^2+d-s}{2}-
  C_{d-s}\gamma_L\right)z_n^2\right\} \cr
  &\le &
  C_{d-s}\exp\left\{-\left(\frac{d-s}{2}
   +\frac{(1-\gamma_1)^2}{3} \right) z_n^2\right\}.
\end{eqnarray*}
  Note that $z_n^2 = 2\log n - \log\log n + O(1)$. Pick
  $L_n=\max\{\floor{n^{\alpha}},M\}$ for some $\alpha
  <2(1-\gamma_1^2)/3d$. For any $1 \leq a \leq d-1$, since there are
  at most $L_n^an^{d-a}$ ordered subset $J \subset [n]$ such that
  $s(J,L_n)=a$, we know the sum of $Q_J$ over these $J$ is dominated
  by
  % $$n^{d-a}\cdot n^{\frac{2(d-1)(1-\gamma_1)^2}{3d}} \cdot n^{-(d-a)}
  % \cdot n^{-\frac{2(1-\gamma_1)^2}{3}}$$
  $$C_{d-a}\exp\left\{\log n \left( (d-a) + \frac{2(d-1)(1-\gamma_1)^2}{3d}
      -(d-a) - \frac{2(1-\gamma_1)^2}{3} \right)\right\}$$ when $n$ is
  large enough, which converges to zero. Therefore, it suffices to
  consider all the ordered subsets $J$ such that $l_j>L_n$ for all
  $1\leq j \leq d-1$.

  Let $J=\{t_1,\ldots,t_d\}\subset [n]$ be an ordered subset such that
  $t_i-t_{i-1} > L_n$ for $2 \leq i \leq d$, and $\mathcal{J}(d,L_n)$
  be the collection of all such subsets. Let $(r_{ij})$ be the
  $d$-dimensional covariance matrix of $\h{X}_J$.  There exists a
  matrix $R_J=\theta(r_{ij})_{i,j\in J}+(1-\theta)I_d$ for some $0 <
  \theta < 1$ such that
  \begin{equation*}
    %\label{eq:4}
    Q_J-Q_d(I_d,z_n) = \sum_{h,l \in J, h<l} \frac{\partial
      Q_d}{\partial r_{hl}}[R_J;z_n]r_{ij}.
  \end{equation*}
Let $R_H$, $H = J \setminus \{h,l\}$, be the correlation matrix of
the conditional distribution of $\h{X}_H$ given $X_h$ and $X_l$.
By (\ref{eq:normbdd3}) of Lemma~\ref{thm:normbdd}, for $n$ large
enough
  \begin{align*}
    \frac{\partial Q_{d}}{\partial r_{hl}}[R_J;z_n] & \leq C
    \exp\left\{-\frac{z_n^2}{1+\gamma_{l-h}}\right\} \cdot
    Q_{d-2}\left(R_K;(1-3\gamma_{L_n})z_n\right) \\
    & \leq C C_{d-2}f_{d-2}(\gamma_{L_n},1/z_n)
    \exp\left\{-\frac{z_n^2}{1+\gamma_{l-h}}\right\} \cr
    & \quad \times
      \exp\left\{-\left(\frac{d-2}{2} - 2C_{d-2} \gamma_{L_n}
        \right)(1-3\gamma_{L_n}))^2z_n^2\right\} \\
    & \leq C_df_{d-2}(\gamma_{L_n},1/z_n)
     \exp\left\{-\left(\frac{d}{2} - (2C_{d-2}+3(d-2))
        \gamma_{L_n} -\gamma_{h-l} \right)z_n^2\right\}\\
    & \leq C_df_{d-2}(\gamma_{L_n},1/z_n)
         \exp\left\{-\left(\frac{d}{2} - C_d
        \gamma_{L_n} -\gamma_{h-l} \right)z_n^2\right\}.
  \end{align*}
  It follows that
  \begin{eqnarray}
  \label{eq:5}
  && \sum_{J \in \mathcal{J}(d,L_n)} |Q_J-Q_d(I_d;z_n)| \cr
  && \leq C_d f_{d-2}(\gamma_{L_n},1/z_n)\sum_{J \in
      \mathcal{J}(d,L_n)} \sum_{1\leq i<j\leq d}
    \exp\left\{-\left(\frac{d}{2} - C_d \gamma_{L_n} -\gamma_{t_j-t_i}
      \right)z_n^2\right\} \gamma_{t_j-t_i} \cr
  && = C_d f_{d-2}(\gamma_{L_n},1/z_n)
     \sum_{1\leq i<j\leq d} \sum_{J \in \mathcal{J}(d,L_n)}
      \exp\left\{-\left(\frac{d}{2} - C_d
        \gamma_{L_n} -\gamma_{t_j-t_i} \right)z_n^2\right\}
    \gamma_{t_j-t_i}.
  \end{eqnarray}
  For each fixed pair $1\leq i<j\leq d$, the inner sum in (\ref{eq:5})
  is bounded by
  \begin{align}\nonumber
    & C_d f_{d-2}(\gamma_{L_n},1/z_n)\sum_{l=L_n+1}^{n-1} (n-l)^{d-1} %\gamma_l
    \exp\left\{-\left(\frac{d}{2} - C_d \gamma_{L_n} -\gamma_{l}
      \right)z_n^2\right\} \gamma_l \\\label{eq:7}
    \leq & C_d f_{d-2}(\gamma_{L_n},1/z_n)(\log n)^{d/2} n^{-d}
    \sum_{l=L_n+1}^{n-1}(n-l)^{d-1} \exp\left\{\left( C_d \gamma_{L_n}
        + \gamma_{l} \right)2\log n\right\}\gamma_l \\ \label{eq:6}
    \leq & C_df_{d-2}(\gamma_{\floor{n^\alpha}},1/z_n)\,
     \gamma_{\floor{n^\alpha}}(\log n)^{d/2}
    \exp\left\{2\left( C_d +1 \right)
     \gamma_{\floor{n^\alpha}}\log n \right\}.
  \end{align}
  Since $\lim_{n \to \infty} \gamma_n \log n=0$, it
  also holds that $\lim_{n \to \infty}
  \gamma_{\floor{n^\alpha}} \log n
  =0$. Note that $\lim_{n\to\infty} (\log n)^{1/2}/z_n = 2^{-1/2}$,
  it follows that $\lim_{n\to\infty}
  f_{d-2}(\gamma_{\floor{n^\alpha}},1/z_n) (\log n)^{d/2-1} =
  2^{-d/2+1}$. Therefore, the term in (\ref{eq:6}) converges to zero,
  and the proof is complete.
\end{proof}

\begin{remark}\label{rk:berman}
This lemma provides another proof of Theorem 3.1 in
\cite{berman:1964}, which gives the asymptotic distribution of the
maximum term of a stationary Gaussian process. They also showed
that the theorem is true if the condition $\lim_{n\to\infty}r_n
\log n = 0$ is replaced by $\sum_{n=1}^\infty r_n^2<\infty$. Under
the later condition, if we replace $\gamma_{t_j-t_j}$ by
$|r_{t_j-t_i}|$ in (\ref{eq:5}), $\gamma_l$ by $|r_l|$ in
(\ref{eq:7}), then the term in (\ref{eq:7}) converges to zero, and
hence our result remains true.
\end{remark}

\begin{remark}\label{rk:absolute}
  In the proof, the upper bounds on $Q_J$ and $|Q_J-Q(I_d;z_n)|$ are
  expressed through the absolute values of the correlations, so we can
  obtain the same bounds for probabilities of the form $P(\cap_{1\leq
    i \leq d} \{(-1)^{a_i}X_{t_i} \geq z_n\})$ for any
  $(a_1,\ldots,a_d) \in \{0,1\}^d$. Therefore, our result can be used
  to show the asymptotic distribution of the maximum absolute term of
  a stationary Gaussian process. Specifically, we have
  \begin{align*}
    \lim_{n\to\infty}P\left(\max_{1\leq i\leq n}|X_i| \leq a_{2n}\,x+b_{2n}\right) = \exp\{-\exp(-x)\}.
  \end{align*}
  \cite{deo:1972} obtained this result under the condition
  $\lim_{n\to\infty} r_n (\log n)^{2+\alpha}=0$ for some
  $\alpha>0$, whereas we only need $\lim_{n\to\infty} r_n \log
  n =0$.
\end{remark}

\section{Summability of Cumulants}
\label{sec:cum}
For a $k$-dimensional random vector $(Y_1,\ldots,Y_k)$ such that
$\|Y_i\|_k<\infty$ for $1\leq i \leq k$, the $k$-th order joint
cumulant is defined as
\begin{equation}
  \label{eq:cumulant}
  \Cum(Y_1,\ldots,Y_k) =  \sum (-1)^{p-1}(p-1)!\prod_{j=1}^p\left(\E \prod_{i\in\nu_j}Y_i\right),
\end{equation}
where the summation extends over all partitions
$\{\nu_1,\ldots,\nu_p\}$ of the set $\{1,2,\ldots,k\}$ into $p$
non-empty blocks. For a stationary process $(X_i)_{i\in\Z}$, we
abbreviate
\begin{equation*}
  \gamma(k_1,k_2,\ldots,k_{d}) := \Cum(X_0,X_{k_1},X_{k_2},\ldots,X_{k_{d}}),
\end{equation*}
Summability conditions of cumulants are often assumed in the
spectral analysis of time series, see for example
\cite{brillinger:2001} and \cite{rosenblatt:1985}. Recently, such
conditions were used by \cite{anderson:2008} in studying the
spectral properties of banded sample covariance matrices. While
such conditions are true for some Gaussian processes, functions of
Gaussian processes \citep{rosenblatt:1985}, and linear processes
with iid innovations \citep{anderson:1971}, they are not easy to
verify in general. \cite{wu:2004} showed that the summability of
joint cumulants of order $d$ holds under the condition that
$\delta_d(k)=O(\rho^k)$ for some $0<\rho<1$. We present in
Theorem~\ref{thm:cum} a generalization of their result. To
simplify the proof, we introduce the {composition} of an integer.
A {\it composition} of a positive integer $n$ is an ordered
sequence of strictly positive integers
$\{\upsilon_1,\upsilon_2,\ldots,\upsilon_q\}$ such that
$\upsilon_1+\cdots+\upsilon_q=n$. Two sequences that differ in the
order of their terms define different compositions. There are in
total $2^{n-1}$ different compositions of the integer $n$. For
example, we are giving in the following all of the eight
compositions of the integer 4.
\begin{align*}
  \{1,1,1,1\} \quad \{1,1,2\} \quad \{1,2,1\}
  \quad \{1,3\}  \quad \{2,1,1\}  \quad \{2,2\}
  \quad \{3,1\}  \quad \{4\}.
\end{align*}

\begin{theorem}\label{thm:cum}
  Assume $d\geq 2$, $X_i \in \LLL^{d+1}$ and $\E X_i=0$. If
\begin{eqnarray}
\label{eq:kddcm}
  \sum_{k=0}^{\infty}k^{d-1}\delta_{d+1}(k)<\infty,
\end{eqnarray}
  then
  \begin{equation}\label{eq:smcmD25}
    \sum_{k_1,\ldots,k_{d}\in \Z} |\gamma(k_1,k_2,\ldots,k_{d})| < \infty.
  \end{equation}
\end{theorem}
\begin{proof}
By symmetry of the cumulant in its arguments and stationarity of
the process, it suffices to show
  \begin{equation*}
    \sum_{0 \leq k_1\leq k_2\leq\cdots\leq k_{d}} |\gamma({k_1}, {k_2},\ldots,{k_{d}})| < \infty.
  \end{equation*}
  Set $X(k,j):=\HH_jX_k$, we claim
  % \begin{align*}
  %   \gamma({k_1}, {k_2},\ldots,{k_{d}}) = \sum
  %   \Cum & \left[X_0,\HH_{1}X_{k_1},\ldots,\HH_1X_{k_{\upsilon_1-1}}, X_{k_{\upsilon_1}}-\HH_1X_{k_{\upsilon_1}},
  %     \phantom{\HH_{k_{\upsilon_{q}}+1}}\right.\cr
  %   & \left.\phantom{[}\HH_{k_{\upsilon_1}+1}X_{k_{\upsilon_1+1}},\ldots,\HH_{k_{\upsilon_1}+1}X_{k_{\upsilon_2-1}},
  %     X_{k_{\upsilon_2}}-\HH_{k_{\upsilon_1}+1}X_{k_{\upsilon_2}},\right.\cr
  %   & \left.\phantom{[} \cdots, \right. \cr
  %   & \left.\phantom{[} \HH_{k_{\upsilon_{q}}+1}X_{k_{\upsilon_{q}+1}},\ldots,\HH_{k_{\upsilon_{q}}+1}X_{k_{d-1}},
  %   X_{k_{d}}-\HH_{k_{\upsilon_q}+1}X_{k_{d}}\right] \cr
  % \end{align*}
  \begin{align}\label{eq:tele}
    \gamma({k_1}, {k_2},\ldots,{k_{d}}) = \sum
    \Cum & \left[X_0, X(k_1,1),\ldots,X({k_{\upsilon_1-1}},1), X_{k_{\upsilon_1}}-X({k_{\upsilon_1}},1),\right.\cr
    & \left.\phantom{[}X({k_{\upsilon_1+1}},{k_{\upsilon_1}+1}),\ldots,X({k_{\upsilon_2-1}},{k_{\upsilon_1}+1}),
      X_{k_{\upsilon_2}}-X({k_{\upsilon_2}},{k_{\upsilon_1}+1}),\right.\cr
    & \left.\phantom{[} \cdots, \right. \cr
    & \left.\phantom{[} X({k_{\upsilon_{q}+1}},{k_{\upsilon_{q}}+1}),\ldots,X({k_{d-1}},{k_{\upsilon_{q}}+1}),
    X_{k_{d}}-X({k_{d}},{k_{\upsilon_{q}}+1})\right];
  \end{align}
  where the sum is taken over all the $2^{d-1}$ increasing sequences
  $\{\upsilon_0,\upsilon_1,\ldots,\upsilon_q,\upsilon_{q+1}\}$ such
  that $\upsilon_0=0$, $\upsilon_{q+1}=d$ and
  $\{\upsilon_1,\upsilon_2-\upsilon_1,\ldots,\upsilon_q-\upsilon_{q-1},d-\upsilon_q\}$
  is a composition of the integer $d$.
  % We claim that
  % \begin{align}\label{eq:37}
  %   \gamma({k_1}, {k_1+k_2},\ldots,{k_1+\cdots+k_{d-1}})
  %   \leq C \sum \prod_{j=1}^q \Theta_d\left(\sum_{i=\Upsilon_{j-1}+1}^{\Upsilon_j}k_i\right),
  % \end{align}
  % where $\Upsilon_0=0$, $\Upsilon_j=\sum_{i=1}^j \upsilon_i$ for $1
  % \leq j \leq q$, and the summation is taken over all the $2^{d-2}$
  % compositions $\{\upsilon_1,\upsilon_2,\ldots,\upsilon_q\}$ of the
  % integer $d-1$. Observing that the condition $\sum_{k=0}^\infty
  % k^{d-1}\delta_d(k)<\infty$ is equivalent to $\sum_{k=0}^\infty
  % k^{d-2}\Theta_d(k) <\infty$, the lemma follows in view of
  % \begin{align*}
  %   \sum_{k_1,\ldots,k_{d-1}\in \N} \Theta_d\left(k_1+k_2+\cdots+k_{d-1}\right)
  %   \leq \sum_{k=0}^{\infty}{k+d-2 \choose d-2}\Theta_d(k).
  % \end{align*}
  % Now we shall prove (\ref{eq:37}).
  We first consider the last summand which corresponds to the sequence $\{\upsilon_0=0,\upsilon_1=d\}$,
  \begin{align*}
    \Cum\left[X_0,X(k_1,1),\ldots,X(k_{d-1},1),X_{k_d}-X(k_d,1)\right]
  \end{align*}
  Observe that $X_0$ and $(X(k_1,1),\ldots,X(k_{d-1},1))$ are
  independent. By definition, only partitions for which $X_0$ and
  $X_{k_d}-X(k_d,1)$ are in the same block contribute to the sum in
  (\ref{eq:cumulant}). Suppose $\{\nu_1,\ldots,\nu_p\}$ is a partition
  of the set $\{k_1,k_2,\ldots,k_{d-1}\}$, since
  \begin{align*}
    \left|\E\left[X_0(X_{k_d}-X(k_d,1))\prod_{k\in\nu_1}X(k,1)\right]\right|
    & = \left|\sum_{j=-\infty}^0 \E\left[\PP_j X_0 \PP_jX_{k_d}\prod_{k\in\nu_1}X(k,1)\right]\right|\cr
    & \leq \sum_{j=-\infty}^0 \delta_{d+1}(-j)\delta_{d+1}(k_d-j)\kappa_{d+1}^{|\nu_1|},
  \end{align*}
  it follows that
  \begin{align*}
    \left|\E\left[X_0(X_{k_d}-X(k_d,1))\prod_{k\in\nu_1}X(k,1)\right]
    \cdot\prod_{j=2}^p\left(\E \prod_{k\in\nu_j}X(k,1)\right)\right|
    \leq \sum_{j=0}^\infty \delta_{d+1}(j)\delta_{d+1}(k_d+j)\kappa_{d+1}^{d-1}
  \end{align*}
  and therefore
  \begin{align*}
    & \sum_{0 \leq k_1\leq k_2\leq\cdots\leq k_{d}}
    \left|\Cum\left[X_0,X(k_1,1),\ldots,X(k_{d-1},1),X_{k_d}-X(k_d,1)\right]\right|\cr
    & \leq C_d  \sum_{0 \leq k_1\leq k_2\leq\cdots\leq k_{d}}\sum_{j=0}^\infty \delta_{d+1}(j)\delta_{d+1}(k_d+j)
    \leq C_d \sum_{j=0}^\infty\sum_{k=0}^{\infty}{k+d-1\choose d-1}\delta_{d+1}(j)\delta_{d+1}(k+j) < \infty,
  \end{align*}
  provided that $\sum_{k=0}^\infty k^{d-1}\delta_{d+1}(k)<\infty$.

  The other terms in (\ref{eq:tele}) are easier to deal with. For
  example, for the term corresponding to the sequence
  $\{\upsilon_0=0,\upsilon_1=1,\upsilon_2=d\}$, we have
  \begin{align*}
    & \left|\Cum\left[X_0,X_{k_1}-X(k_1,1),X(k_2,k_1+1),\ldots,X(k_{d-1},k_1+1),X_{k_d}-X(k_d,k_1+1)\right]\right| \cr
    & \leq C_d \kappa_{d+1}^{d-1}\Psi_{d+1}(k_1)\Psi_{d+1}(k_d-k_1).
  \end{align*}
  Since $\sum_{k=0}^\infty k^{d-1}\delta_{d+1}(k)<\infty$ implies
  $\sum_{k=0}^\infty k^{d-2}\Psi_{d+1}(k) \leq \infty$, it follows
  that
  \begin{align*}
    \sum_{0 \leq k_1\leq k_2\leq\cdots\leq k_{d}}
    & \left|\Cum\left[X_0,X_{k_1}-X(k_1,1),X(k_2,k_1+1),\ldots,X(k_{d-1},k_1+1),X_{k_d}-X(k_d,k_1+1)\right]\right| \cr
    & \leq C_d\kappa_{d+1}^{d-1} \sum_{k=0}^\infty \Psi_{d+1}(k) \sum_{k=0}^\infty {k+d-2 \choose d-2}\Psi_{d+1}(k) \leq \infty.
  \end{align*}
We have shown that every cumulant in (\ref{eq:tele}) is absolutely
summable over $0\leq k_1\leq \cdots\leq k_d$, and it remains to
show the claim (\ref{eq:tele}). We shall derive the case $d=3$,
(\ref{eq:tele}) for other values of $d$ are obtained using the
same idea. By multilinearity of cumulants, we have
  \begin{align*}
    \gamma(k_1,k_2,k_3) = & \Cum(X_0,X_{k_1},X_{k_2},X_{k_3})\cr
    = &\Cum\left[X_0,X_{k_1}-X(k_1,1),X_{k_2},X_{k_3}\right] \cr
    & + \Cum\left[X_0,X(k_1,1),X_{k_2}-X(k_2,1),X_{k_3}\right] \cr
    & +
    \Cum\left[X_0,X(k_1,1),X({k_2},1),X_{k_3}-X(k_3,1)\right]\cr
    & + \Cum\left[X_0,X(k_1,1),X({k_2},1),X(k_3,1)\right].
  \end{align*}
  Since $X_0$ and $(X(k_1,1),X({k_2},1),X(k_3,1))$ are independent,
  the last cumulant is 0. Apply the same trick for the first two
  cumulants, we have
  \begin{eqnarray*}
    && \Cum\left[X_0,X_{k_1}-X(k_1,1),X_{k_2},X_{k_3}\right] \cr
    && = \Cum\left[X_0,X_{k_1}-X(k_1,1),
               X_{k_2}-X(k_2,k_1+1),X_{k_3}\right]\cr
    && + \Cum\left[X_0,X_{k_1}-X(k_1,1),
               X(k_2,k_1+1),X_{k_3}-X(k_3,k_1+1)\right]\cr
    && + \Cum\left[X_0,X_{k_1}-X(k_1,1),
               X(k_2,k_1+1),X(k_3,k_1+1)\right]\cr
    && = \Cum\left[X_0,X_{k_1}-X(k_1,1),X_{k_2}-X(k_2,k_1+1),
               X_{k_3}-X(k_3,k_2+1)\right]\cr
    && + \Cum\left[X_0,X_{k_1}-X(k_1,1),X(k_2,k_1+1),
               X_{k_3}-X(k_3,k_1+1)\right]
  \end{eqnarray*}
and
   \begin{align*}
    \Cum\left[X_0,X(k_1,1),X_{k_2}-X(k_2,1),X_{k_3}\right] =
    \Cum\left[X_0,X(k_1,1),X_{k_2}-X(k_2,1),X_{k_3}-X(k_3,k_2+1)\right].
  \end{align*}
Then the proof is complete.
\end{proof}

\begin{remark}
When $d = 1$, (\ref{eq:kddcm}) reduces to the {\it short-range
dependence} or {\it short-memory} condition $\Theta_2 =
\sum_{k=0}^\infty \delta_2(k) < \infty$. If $\Theta_2 = \infty$,
then the process $(X_i)$ may be long-memory in that the
covariances are not summable. When $d \geq 2$, we conjecture that
(\ref{eq:kddcm}) can be weakened to $\Theta_{d+1} < \infty$. It
holds for linear processes. Let $X_k = \sum_{i=0}^\infty a_{i}
\epsilon_{k-i}$. Assume $\epsilon_k \in \LLL^{d+1}$ and
$\sum_{k=0}^\infty|a_k| < \infty$, then $\delta_{d+1}(k) = |a_k|
\|\epsilon_0\|_{d+1}$. Let $\Cum_{d+1}(\epsilon_0)$ be the
$(d+1)$-th cumulant of $\epsilon_0$. Set $k_0=0$, by
multilinearity of cumulants, we have
\begin{eqnarray*}
\gamma(k_1,\ldots,k_d)
 &=&\sum_{t_0,t_1,\ldots,t_d\geq 0}
    \left(\prod_{j=0}^{d}a_{t_j}\right)
    \Cum(\epsilon_{-t_0},\epsilon_{k_1-t_1},
     \ldots,\epsilon_{k_d-t_d})\cr
  &=& \sum_{t=0}^\infty \prod_{j=0}^{d}a_{k_j+t}
  \Cum_{d+1}(\epsilon_0).
\end{eqnarray*}
Therefore, the condition $\Theta_{d+1} < \infty$ suffices for
(\ref{eq:smcmD25}). For a class of functionals of Gaussian
processes, \cite{rosenblatt:1985} showed that (\ref{eq:smcmD25})
holds if $\sum_{k=0}^\infty|\gamma_k|<\infty$, which in turn is
implied by $\Theta_{d+1}<\infty$ under our setting. It is unclear
whether in general the weaker condition $\Theta_{d+1} < \infty$
implies (\ref{eq:smcmD25}).
\end{remark}

\section{Some Auxiliary Lemmas}
\label{sec:aux} Suppose that $\h{X}$ is a $d$-dimensional random
vector, and $\h{X} \sim \mathcal{N}(0,\Sigma)$. If $\Sigma=I_d$,
then by (\ref{eq:normbdd}), it is easily seen that the ratio of
$P\left(z_n-c_n\leq|\h{X}|_\bullet\leq z_n\right)$ over
$P\left(|\h{X}|_\bullet\geq z_n\right)$ tends to zero provided
that $c_n\to 0$, $z_n\to \infty$ and $c_nz_n\to 0$. It is a
similar situation when $\Sigma$ is not an identity matrix, as
shown in the following lemma, which will be used in the proof of
Lemma~\ref{thm:vec_mod_dev}.
\begin{lemma}
  \label{thm:lemma1}
Let $\h{X} \sim \mathcal{N}(0,\Sigma)$ be a $d$-dimensional normal
random vector. Assume $\Sigma$ is nonsingular. Let $\lambda_0^2$
and $\lambda_1^2$ be the smallest and largest eigenvalue of
$\Sigma$ respectively. Then for $0<c<\delta<1/2$ such that $A :=
(2\pi\lambda_1^2)^{(d-1)/2}\lambda_0^2c^2\delta^{-2} + d \delta
\exp\{(\sqrt{6}d\lambda_1+\lambda_0)/\lambda_0^3\} < 1$, then for
any $z \in [1, \delta/c]$,
  \begin{equation}\label{eq:9}
    P\left(z-c \leq \|\h{X}\|_{\bullet} \leq z\right)
    \leq (1-A)^{-1} A \,P\left(\|\h{X}\|_{\bullet} \geq z\right).
  \end{equation}
\end{lemma}

\begin{proof}
Let $C_d = (6d)^{1/2}\lambda_1/\lambda_0$. Since $\lambda_0^2$ is
the smallest eigenvalue of $\Sigma$,
\begin{eqnarray*}
  P(\|\h{X}\|_{\bullet}\geq z-c)
   &\geq& {(2\pi\det(\Sigma))^{-d/2}}
    \exp\left\{-\frac{d(z+1)^2}{2\lambda_0^2}\right\} \cr
   &\geq& (2\pi\lambda_1^2)^{-d/2}
  \exp\left\{-\frac{4d\delta^2}{2\lambda_0^2c^2}\right\}.
\end{eqnarray*}
Since $P(\|\h{X}\|_\infty \geq C_d \delta/c) \leq d
(2\pi\lambda_1^2)^{-1/2} \exp\{6d\delta^2/(2\lambda_0^2c^2)\}$, we
have
  \begin{equation}
    \label{eq:1}
    P(\|\h{X}\|_\infty \geq C_d\delta/c) \leq
    (2\pi\lambda_1^2)^{(d-1)/2}\lambda_0^2c^2\delta^{-2} \,
    P(\|\h{X}\|_{\bullet}\geq z-c).
  \end{equation}
  For $0\leq k \leq
  \floor{1/\delta}$, define the orthotopes
  $R_k=[z+(k-1)c,z+kc]\times[z-c,C_d\delta/c]^{d-1}$. For two points
  $\h{x}=(x_1,\ldots,x_d) \in R_0$, $\h{x}_k=(x_1+kc,x_2,\ldots,x_d)
  \in R_k$, we have $\h{x}_k^\top\Sigma^{-1}\h{x}_k
  -\h{x}^\top\Sigma^{-1}\h{x} \leq
  (2\sqrt{d}C_d+1)/\lambda_0^2$, and hence
  $P(\h{X} \in R_k) \geq
  \exp\{-(\sqrt{d}C_d+1)/\lambda_0^2\}P(\h{X}\in
  R_0)$ for any $1\leq k \leq \floor{1/\delta}$. Since the same
  inequality holds for every coordinate, we have
  \begin{equation}
    \label{eq:8}
    P\left(z-c \leq \|\h{X}\|_{\bullet} \leq
      z,\,\|\h{X}\|_\infty \leq C_d\delta/c\right) \leq
    d\delta\exp\{(\sqrt{d}C_d+1)/\lambda_0^2\}\,
    P\left(\|\h{X}\|_{\bullet} \geq z-c\right)
  \end{equation}
Combine (\ref{eq:1}) and (\ref{eq:8}), we know $P\left(z-c \le
\|\h{X}\|_{\bullet} \leq z\right) \le A \cdot
P\left(\|\h{X}\|_{\bullet} \geq z-c\right)$. So (\ref{eq:9})
follows.
\end{proof}

The preceding lemma requires the eigenvalues of $\Sigma$ to be
bounded both from above and away from zero. In our application,
$\Sigma$ is taken as the covariance matrix of $(G_{k_1}, G_{k_2},
\ldots, G_{k_d})^\top$, where $(G_k)$ is defined in (\ref{eq:gp}).
Furthermore, we need such bounds be uniform over all choices of
$k_1<k_2<\cdots<k_d$. Let $f(\omega) = (2\pi)^{-1} \sum_{h \in \Z}
\sigma_h \cos(h\omega)$ be the spectral density of $(G_k)$. A
sufficient condition would be that there exists $0<m<M$ such that
\begin{equation}
  \label{eq:32}
  m \leq f(\omega) \leq M, \quad\hbox{for } \omega \in [0,2\pi],
\end{equation}
because the eigenvalues of the autocovariance matrix are bounded from
above and below by the maximum and minimum values that $f$ takes
respectively. For the proof see Section~5.2 of
\cite{grenander:1958}. Clearly the upper bound in (\ref{eq:32}) is
satisfied in our situation, because $\sum_{h\in\Z} |\sigma_h| <
\infty$. However, the existence of lower bound in (\ref{eq:32}) rules
out some classical times series models. For example, if $(G_k)$ is the
moving average of the form $G_k = (\eta_k + \eta_{k-1}) / \sqrt{2}$,
then $f(\omega) = (1+\cos(\omega)) / 2\pi$, and
$f(\pi)=0$. Nevertheless, although the minimum eigenvalue of the
autocovariance matrix converges to $\inf_{\omega\in[0,2\pi]}
f(\omega)$ as the dimension of the matrix goes to infinity, there does
exist a positive lower bound for the smallest eigenvalues of all the
principal sub-matrices with a fixed dimension.
\begin{lemma}
  \label{thm:eigenbdd}
  If $\sum_{h \in \Z}\sigma_h^2<\infty$, then for each $d\geq 1$, there
  exists a constant $C_d>0$ such that
  \begin{align*}
    \inf_{k_1<k_2<\cdots<k_d} \lambda_{\min}
    \left\{\Cov\left[(G_{k_1},G_{k_2},
             \ldots,G_{k_d})^\top\right]\right\}
    \geq C_d.
  \end{align*}
\end{lemma}
\begin{proof}
We use induction. It is clear that we can choose $(C_d)$ to be a
non-increasing sequence. Without loss of generality, let us assume
$k_1=1$. The statement is trivially true when $d=1$. Suppose it is
true for all dimensions up to $d$, we now consider the dimension
$(d+1)$ case. There exist an integer $N_{d}$ such that
  $\sum_{h=N_d}\sigma_h^2<2C_d^2/(d+1)$. If all the differences
  $k_{i+1}-k_i\leq N_{d}$ for $1\leq i\leq d-1$, there are $N_d^{d-1}$
  possible choices of $k_1=1<k_2<\cdots<k_d$. Since the process
  $(G_k)$ is non-deterministic, for all these choices, the
  corresponding covariance matrices are non-singular. Pick $C_d'>0$ to
  be the smallest eigenvalue of all these matrices. If there is one
  difference $k_{l+1}-k_l> N_{d}$, set $\Sigma_1=\Cov[(G_{k_i})_{1\leq
    i\leq l}]$ and $\Sigma_2=\Cov[(G_{k_i})_{l< i\leq d}]$, then
  $\lambda_{\min}(\Sigma_1)\geq C_d$ and $\lambda_{\min}(\Sigma_2)
  \geq C_d$.  It follows that for any real numbers
  $c_1,c_2,\ldots,c_d$ such that $\sum_{i=1}^dc_i^2=1$,
  \begin{eqnarray*}
    \sum_{1\leq i,j\leq d}c_ic_j\Cov(G_{k_i},G_{k_j})
    & = & (c_1,\ldots,c_i)^\top\Sigma_J(c_1,\ldots,c_i) \cr
    & & \quad +
    (c_{i+1},\ldots,c_d)^\top\Sigma_J(c_{i+1},\ldots,c_d)\cr
    & & \quad  + 2\sum_{i\leq l,j>l}c_ic_j\sigma_{k_j-k_i} \cr
    & \geq & C_d - 2\left(\sum_{i\leq l,j>l}\sigma_{k_j-k_i}^2\right)^{1/2}
      \left(\sum_{i\leq l,j>l}c_i^2c_j^2\right)^{1/2}\cr
    & \geq & C_d - \frac{1}{2}\left( {{d+1}\over 2}\cdot
     \sum_{h=N_d}\sigma_h^2\right)^{1/2} \geq {{C_d} \over 2}.
  \end{eqnarray*}
  Setting $C_{d+1} = \min\{C_d/2,C_d'\}$, the proof is complete.
\end{proof}

The following lemma is used in the proof of Lemma~\ref{thm:cov_struc}.
\begin{lemma}
  \label{thm:covconvergence}
  Assume $X_i \in \mathcal{L}^4$, $\E X_0=0$, and
  $\Theta_{4}<\infty$. Assume $l_n\to\infty$, $k_n \to
  \infty$, ${\check m_n} < \floor{k_n/3}$ and $h\geq 0$. Define
  ${S}_{n,k}=\sum_{i=1}^{l_n} (X_{i-k}X_i-\gamma_k)$. Then
  \begin{align}
    \label{eq:19}
    % \left|\left\|({S}_{n,k_n}-n\gamma_{k_n})/\sqrt{l_n}\right\|^2 -\sigma_{0} \right|
    % & \leq \Theta_{4}^3\left(12\Delta_4({\check m_n}+1) + 6\Theta_4\sqrt{{\check m_n}/l_n} + 4\Psi_4({\check m_n}+1)\right); \\ \label{eq:20}
    \left|\E\left({S}_{n,k_n}{S}_{n,k_n+h}\right)/l_n - \sigma_h\right|
    & \leq \Theta_{4}^3\left(16\Delta_4({\check m_n}+1) + 6\Theta_4\sqrt{{\check m_n}/l_n} + 4\Psi_4({\check m_n}+1)\right).
  \end{align}
\end{lemma}
\begin{proof}
  Let $\check X_{i} = \HH_{i-{\check m_n}}^i X_i$, then $\check X_i$ and
  $\check X_{i-k_n}$ are independent, because ${\check m_n} \leq
  \floor{k_n/3}$. Define $\check {S}_{n,k}=\sum_{i=1}^{l_n} \check
  X_{i-k}\check X_i$. By (\ref{eq:fact9}), we have for any $k \geq 0$,
  \begin{equation}\label{eq:27}
    \left\|({S}_{n,k}-\check {S}_{n,k})/\sqrt{l_n}\right\|
    % \leq \sum_{i=0}^{\infty}
    % \min\{2\kappa_4\Psi_4({\check m_n}+1),2\kappa_4[\delta_4(i)+\delta(i-k)]\}
    \leq 4\kappa_4\Delta_{4}({\check m_n}+1).
  \end{equation}
  By (\ref{eq:fact5}), $\left\|{S}_{n,k}/\sqrt{l_n}\right\| \leq
  2\kappa_4\Theta_4$ for any $k\geq 0$, and it follows that
  \begin{equation}
    \label{eq:24}
   \begin{aligned}
     & \left|\E({S}_{n,k_n},{S}_{n,k_n+h})
      - \E(\check {S}_{n,k_n}\check {S}_{n,k_n+h})\right| \cr
     & \quad \leq \left\| {S}_{n,k_n} - \check {S}_{n,k_n} \right\|
     \cdot \left\|{S}_{n,k_n+h}\right\|
       +  \left\|\check {S}_{n,k_n}\right\|
        \cdot \left\|{S}_{n,k_n+h} - \check {S}_{n,k_n+h}
        \right\|\cr
     &\quad\le 16 l_n\kappa_4^2\Theta_4\Delta_4({\check m_n}+1).
   \end{aligned}
   \end{equation}
   For any $k > 3{\check m_n}$, define $M_{n,k}=\sum_{j=1}^{l_n}D_j$, where
   $D_j=\sum_{i=j}^{j+{\check m_n}}\check X_{i-k}\PP^j\check X_i =
   \sum_{q=0}^{{\check m_n}} X_{j+q-k}\PP^jX_{j+q}$. Observe that $\PP^j\check
   X_{j+q}$ and $\check X_{j+q-k}$ are independent, we have
  \begin{align}
    \left\|\check {S}_{n,k} - M_{n,k}\right\|
    & = \left\|\sum_{i=1}^{l_n}\sum_{j=i-{\check m_n}}^i \check X_{i-k}\PP^j\check X_i
      - \sum_{j=1}^{l_n}\sum_{i=j}^{j+{\check m_n}}\check X_{i-k}\PP^j\check X_i\right\| \cr
    & \leq \left\|\sum_{j=1-{\check m_n}}^{0}\sum_{i=1}^{j+{\check m_n}} \check X_{i-k}\PP^j\check X_i\right\|
      + \left\|\sum_{j=l_n-{\check m_n}+1}^{l_n} \sum_{i=l_n+1}^{j+{\check m_n}} \check X_{i-k}\PP^j\check X_i\right\| \cr \label{eq:21}
    & \leq 2\left(\sum_{j=1}^{{\check m_n}} \kappa_2^2\Theta_2(j)^2\right)^{1/2} \leq 2\kappa_2\Theta_2\sqrt{{\check m_n}}
  \end{align}
According to the proof of Theorem~2 of \cite{wu:2009}, when $k >
3{\check m_n}$ $\|M_{n,k}/\sqrt{n}\|^2 = \sum_{k \in \Z} \check
\gamma_k^2$, where $\check\gamma_k=\E\check X_{0}\check X_k$.  By
(\ref{eq:fact4}) and (\ref{eq:fact6}), $|\check\gamma_k| \leq
\zeta_{k}$; and hence
\begin{align}\label{eq:22}
    \left\|M_{n,k}/\sqrt{n}\right\|^2
    &\leq \sum_{k \in \Z} \zeta_k^2
    =\sum_{j,j'=0}^\infty \left(\delta_2(j)\delta_2(j')
    \sum_{k \in \Z} \delta_2(j+k)\delta_2(j'+k)\right)\cr
    & \leq \sum_{j,j'=0}^\infty \delta_2(j)\delta_2(j')
     \Psi_2^2 \leq \Theta_2^2\Psi_2^2.
\end{align}
By (\ref{eq:fact5}) and (\ref{eq:fact6}),  $\left\|\check
{S}_{n,k}/\sqrt{l_n}\right\| \leq 2\kappa_4\Theta_4$ for any
$k\geq 0$. Combining (\ref{eq:21}) and (\ref{eq:22}), we have
  \begin{equation}\label{eq:25}
    \left|\E(\check {S}_{n,k_n}\check {S}_{n,k_n+h}) - \E(M_{n,k_n}M_{n,k_n+h})\right|
    \leq (2\kappa_4\Theta_4+\Theta_2\Psi_2)\sqrt{l_n} \cdot 2\kappa_2\Theta_2\sqrt{{\check m_n}}.
  \end{equation}
  % $$\left|\left\|({S}_{n,k_n}-n\gamma_{k_n})/\sqrt{n}\right\|^2 - \left\|M_{n,k_n}/\sqrt{n}\right\|^2\right|
  % \leq (2\kappa_4\Theta_4+\Theta_2\Psi_2) \left(\sum_{i=0}^{\infty}
  %   4\kappa_4\min\{\Psi_4({\check m_n}+1),\delta_4(i)\} +
  %   2\kappa_2\Theta_2\sqrt{\frac{{\check m_n}}{n}}\right)$$
Observe that when $k_n > 3{\check m_n}$, $X_{q-k_n} X_{q'-k_n-h}$
and $\PP^0 X_q \PP^0 X_{q'}$ are independent for $0 \leq q,q' \leq
{\check m_n}$. Therefore,
\begin{align}
\label{eq:23}
 \E(M_{n,k_n}M_{n,k_n+h})
 & = l_n \E\left(\sum_{q,q'=0}^{{\check m_n}}
     X_{q-k_n}X_{q'-k_n-h}
    \PP^0\check X_q\PP^0\check X_{q'}\right) \cr
    & = l_n \sum_{q,q'=0}^{{\check m_n}}
     \check\gamma_{q-q'+h}
     \E\left[(\PP^0\check X_q)(\PP^0\check X_{q'})\right] \cr
    & = l_n \sum_{k \in \Z} \check\gamma_{k+h}
     \sum_{q'\in \Z}
      \E\left[(\PP^0\check X_{q'+k})(\PP^0\check X_{q'})\right]\cr
    &= l_n \sum_{k \in \Z} \check\gamma_{k+h}
     \sum_{q'\in \Z}
     \E\left[(\PP^{q'}\check X_{k})(\PP^{q'}\check
     X_{0})\right]\cr
    & = l_n \sum_{k \in \Z} \check\gamma_{k+h}\check\gamma_{k}.
  \end{align}
By (\ref{eq:fact7}), $|\gamma_k-\check \gamma_k| \leq 2\kappa_2
\Psi_2(m+1)$. Since $|\gamma_k|\leq\zeta_k$ and $|\check \gamma_k|
\leq \zeta_k$, we have
\begin{eqnarray}\label{eq:26}
    \left|\sigma_h
       -\sum_{k \in \Z} \check\gamma_{k+h}\check\gamma_{k}\right|
    &=& \left|\sum_{k \in \Z}
     (\gamma_k\gamma_{k+h}
        - \check\gamma_k\check\gamma_{k+h})\right| \cr
    &\leq& 4\kappa_2\Psi_2(m+1) \sum_{k \in \Z} \zeta_k
    \leq 4\kappa_2\Psi_2(m+1)\Theta_2^2.
\end{eqnarray}
Combining (\ref{eq:24}), (\ref{eq:25}) and (\ref{eq:26}), the
lemma follows by noting that $\kappa_2$, $\kappa_4$ are dominated
by $\Theta_4$; and $\Theta_2(\cdot)$, $\Psi_2(\cdot)$ and
$\Psi_4(\cdot)$ are all dominated by $\Theta_4(\cdot)$.
\end{proof}

\renewcommand{\baselinestretch}{1}
\bibliographystyle{imsart-nameyear}
\bibliography{max_auto_corr}

% \nocite{einmahl:1989} \nocite{wu:2005} \nocite{wu:2009}
% \nocite{wiener:1958} \nocite{box:1970} \nocite{ljung:1978}
% \nocite{romano:1996} \nocite{schott:2005} \nocite{rosenblatt:1985}
% \nocite{anderson:1971} \nocite{fan:2009} \nocite{bickel:2008b}
% \nocite{wu:2004} \nocite{slepian:1962} \nocite{hall:1979}
% \nocite{petrov:1975} \nocite{amosova:1972} \nocite{wu:2008}

\end{document}